COSTAS J. PAPACHRISTOU

# ASPECTS OF INTEGRABILITY OF DIFFERENTIAL SYSTEMS AND FIELDS

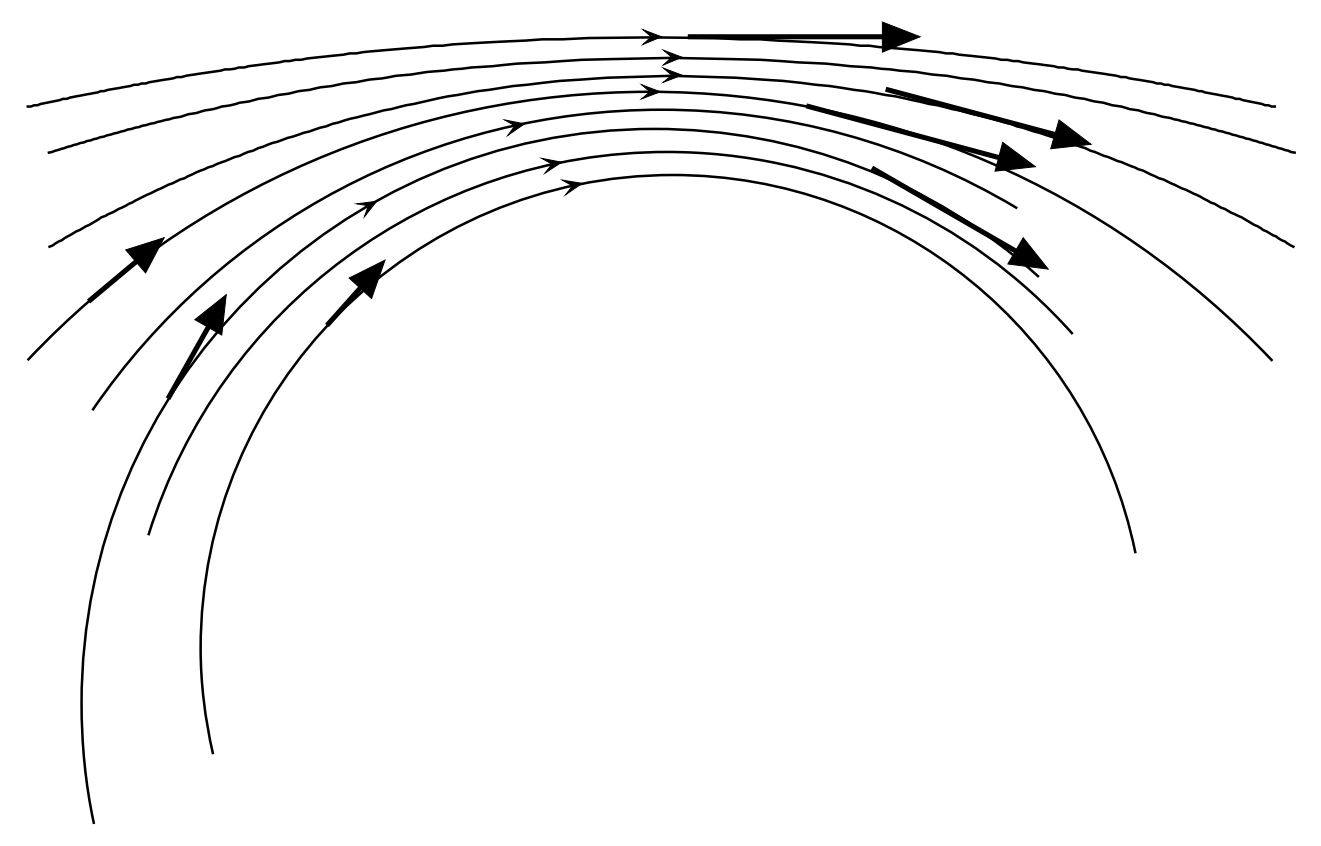

**Aspects of**

# INTEGRABILITY

## of Differential Systems and Fields

**Costas J. Papachristou**

*Department of Physical Sciences*
*Hellenic Naval Academy*

papachristou@hna.gr

# PREFACE

This monograph, written at an intermediate level for educational purposes, serves as an introduction to the concept of integrability as it applies to systems of differential equations (both ordinary and partial) as well as to vector-valued fields. We stress from the outset that this is *not* a treatise on the theory or the methods of solution of differential equations! Instead, we have chosen to focus on specific aspects of integrability that are often encountered in a variety of problems in Applied Mathematics, Physics and Engineering.

With regard to Physics, in particular, integrability is a subject of major importance given that most physical principles are expressed mathematically as systems of differential equations. In Classical Mechanics, certain mathematical techniques are employed in order to integrate the equations of a Newtonian or a Hamiltonian system. These methods involve concepts such as conservation laws, which furnish a number of constants of the motion for the system. In Electrodynamics, on the other hand, the integrability (self-consistency) of the Maxwell system of equations is seen to be intimately related to the wavelike behavior of the electromagnetic field. In the static case, the integrability (in the sense of path-independence) of the electric field leads to the concept of the electrostatic potential. Finally, a number of methods have been developed for finding solutions of nonlinear partial differential equations that are of interest in Mathematical Physics. Before embarking on the study of advanced Physics problems, therefore, the student will benefit by being exposed to some fundamental ideas regarding the mathematical concept of integrability in its various forms.

The following cases of integrability are examined in this book: (*a*) path-independence of line integrals of vector fields on the plane and in space; (*b*) integration of a system of ordinary differential equations (ODEs) by using first integrals; and (*c*) integrable systems of partial differential equations (PDEs). Special topics include the integration of analytic functions and some elements from the geometric theory of differential systems. Certain more advanced subjects, such as Lax pairs and Bäcklund transformations, are also discussed. The presentation sacrifices mathematical rigor in favor of simplicity, as dictated by pedagogical logic. For a deeper study of the subject the reader is referred to the literature cited at the end.

A vector field is said to be integrable in a region of space if its line integral is independent of the path connecting any two points in this region. As will be seen in Chapter 1, this type of integrability is related to the integrability of an associated system of PDEs. Similar remarks apply to the case of analytic functions on the complex plane, examined in Chapter 2. In this case the integrable system of PDEs is represented by the familiar Cauchy-Riemann relations.

In Chapter 3 we introduce the concept of first integrals of ODEs and we demonstrate how these quantities can be used to integrate those equations. As a characteristic example, the principle of conservation of mechanical energy is used to integrate the ODE expressing Newton's second law of motion in one dimension.

This discussion is generalized in Chapter 4 for systems of first-order ODEs, where the solution to the problem is again sought by using first integrals. The method finds an important application in first-order PDEs, the solution process of which is briefly described. Finally, we study the case of a linear system of ODEs, the solution of which reduces to an eigenvalue problem.

Chapter 5 examines systems of ODEs from the geometric point of view. Concepts of Differential Geometry such as the integral and phase curves of a differential system, the differential-operator representation of vector fields, the Lie derivative, etc., are introduced at a fundamental level. The geometric significance of first-order PDEs is also studied, revealing a close connection of these equations with systems of ODEs and vector fields.

Two notions of importance in the theory of integrable nonlinear PDEs are Bäcklund transformations and Lax pairs. In both cases a PDE is expressed as an integrability condition for solution of an associated system of PDEs. These ideas are briefly discussed in Chapter 6. A familiar system of PDEs in four dimensions, namely, the Maxwell equations for the electromagnetic field, is shown to constitute a Bäcklund transformation connecting solutions of the wave equations satisfied by the electric and the magnetic field. The solution of the Maxwell system for the case of a monochromatic plane electromagnetic wave is derived in detail. Finally, the use of Bäcklund transformations as recursion operators for producing symmetries of PDEs is described.

I would like to thank my colleague and friend Aristidis N. Magoulas for an excellent job in drawing a number of figures, as well as for several fruitful discussions on the issue of integrability in Electromagnetism!

Costas J. Papachristou
Piraeus, September 2019

## CHAPTER 1

# INTEGRABILITY ON THE PLANE AND IN SPACE

## 1.1 Simply and Multiply Connected Domains

We begin with a few basic concepts from Topology that will be needed in the sequel.

A domain $D$ on the plane is said to be *simply connected* if, for every closed curve $C$ within this domain, every point of the plane in the interior of $C$ is also a point of $D$. Alternatively, the domain $D$ is simply connected if every closed curve in $D$ can be shrunk to a point without ever leaving this domain. If this condition is not fulfilled, the domain is called *multiply connected*.

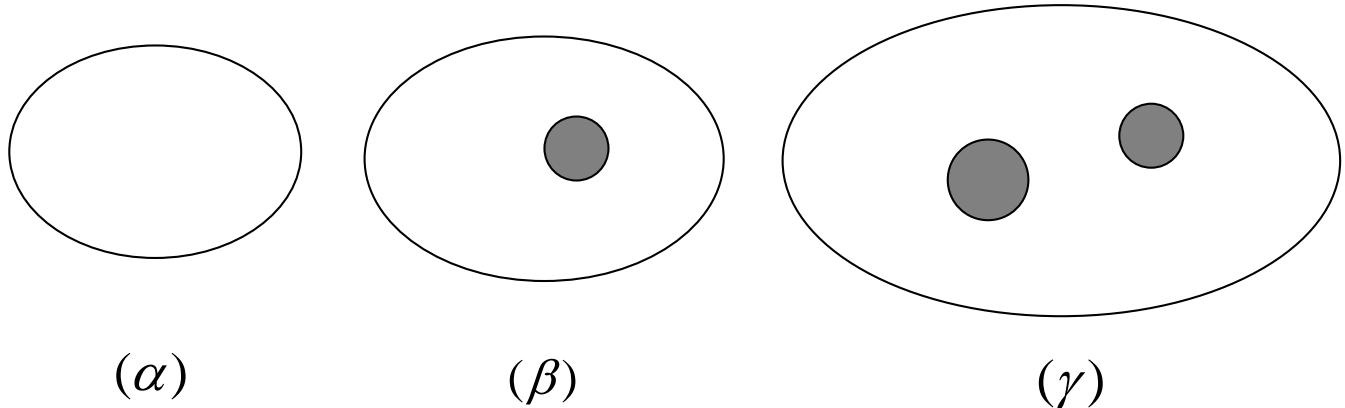

Fig. 1.1. Three domains on the plane, having different types of connectedness.

In Fig. 1.1, the region ($\alpha$) is simply connected, the region ($\beta$) is *doubly* connected while the region ($\gamma$) is *triply* connected. Notice that there are two kinds of closed curves in region ($\beta$): those that do not encircle the "hole" and those that encircle it. (We note that the hole could even consist of a single point subtracted from the plane.) By a similar reasoning, the triple connectedness of region ($\gamma$) is due to the fact that there are three kinds of closed curves in this region: those that do not encircle any hole, those that encircle only one hole (no matter which one!) and those that encircle two holes.

A domain $\Omega$ in space is *simply connected* if, for every closed curve $C$ inside $\Omega$, there is always an open surface bounded by $C$ and located entirely within $\Omega$. This means that every closed curve in $\Omega$ can be shrunk to a point without ever leaving the domain. If this is not the case, the domain is *multiply connected*.

***Examples:***

***1.*** The interior, the exterior as well as the surface of a *sphere* are *simply* connected domains in space. The same is true for a *spherical shell* (the space between two concentric spherical surfaces).

***2.*** The space in the interior of a *torus* (Fig. 1.2) is *doubly* connected (explain why!).

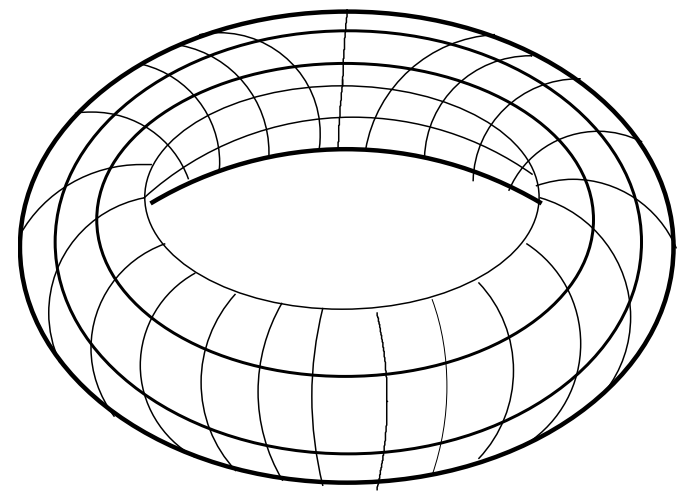

Fig. 1.2. A torus.

## 1.2 Exact Differentials and Integrability

We consider the plane $R^2$ with coordinates $(x,y)$. Let $D \subseteq R^2$ be a domain on the plane and let $P(x,y)$ and $Q(x,y)$ be functions differentiable at every point of $D$. The expression

$$P(x, y)\, dx + Q(x, y)\, dy$$

is an *exact* (or *total*, or *perfect*) *differential* if there exists a function $u(x,y)$, differentiable everywhere in $D$, such that

$$du = P(x, y)\, dx + Q(x, y)\, dy \tag{1}$$

Let us find the *necessary condition* for the existence of $u(x,y)$. In general,

$$du = \frac{\partial u}{\partial x} dx + \frac{\partial u}{\partial y} dy \tag{2}$$

By comparing (1) and (2) and by taking into account that the differentials $dx$ and $dy$ are independent of each other, we find the following system of partial differential equations (PDEs):

$$\frac{\partial u}{\partial x} = P(x, y) \ , \quad \frac{\partial u}{\partial y} = Q(x, y) \tag{3}$$

In order for the system (3) to have a solution for $u$ (that is, to be *integrable*), its two equations must be *compatible* with each other. The *compatibility condition* or *integrability condition* of the system is found as follows: We differentiate the first equation with respect to $y$ and the second one with respect to $x$. By equating the mixed derivatives of $u$ with respect to $x$ and $y$, we find the PDE

$$\frac{\partial P}{\partial y} = \frac{\partial Q}{\partial x} \tag{4}$$

If condition (4) is not satisfied, the system (3) [or, equivalently, the differential relation (1)] does not have a solution for $u$ and the expression $Pdx+Qdy$ is not an exact differential.

***Example:*** $ydx - xdy \neq du$, since $P=y$, $Q=-x$, and $\partial P/\partial y=1$ while $\partial Q/\partial x=-1$.

***Note:*** Relation (4) is a *necessary* condition for the existence of a solution $u$ of the system (3) or, equivalently, of the differential relation (1). This condition will also be *sufficient* if the domain $D \subseteq R^2$ is *simply connected* [1] (by assumption, this is a domain where the functions $P$ and $Q$ are differentiable).

***Examples:***

***1.*** We consider the differential relation

$$du = ydx + xdy.$$

We have $P=y$, $Q=x$, so that $\partial P/\partial y=\partial Q/\partial x=1$. Moreover, the functions $P$ and $Q$ are differentiable everywhere on the plane $R^2$, which is a simply connected space. Thus, the conditions for existence of $u$ are fulfilled. Relations (3) are written

$$\partial u/\partial x = y\ , \quad \partial u/\partial y = x\ .$$

The first one yields

$$u = xy + C(y),$$

where $C$ is an arbitrary function of $y$. Substituting this into the second relation, we find

$$C'(y)=0 \;\;\Rightarrow\;\; C=constant.$$

Thus, finally,

$$u(x,y) = xy + C.$$

***2.*** We consider the relation

$$du = (x+e^y)\, dx + (xe^y - 2y)\, dy.$$

The functions $P=x+e^y$ and $Q=xe^y-2y$ are differentiable on the entire plane $R^2$, which is a simply connected space. Furthermore, $\partial P/\partial y=\partial Q/\partial x=e^y$. Relations (3) are written

$$\partial u/\partial x = x + e^y\ , \quad \partial u/\partial y = xe^y - 2y.$$

By the first one we get

$$u = (x^2/2) + xe^y + \varphi(y) \quad \text{(arbitrary } \varphi\text{)}.$$

Then, the second relation yields

$$\varphi'(y) = -\,2y \;\Rightarrow\; \varphi(y) = -\,y^2 + C.$$

Thus, finally,

$$u(x,y) = (x^2/2) + xe^y - y^2 + C.$$

Consider, now, a domain $\Omega \subseteq R^3$ in a space with coordinates $(x,y,z)$. Also, consider the functions $P(x,y,z)$, $Q(x,y,z)$ and $R(x,y,z)$, differentiable at each point $(x,y,z)$ of $\Omega$. The expression

$$P(x, y, z)\, dx + Q(x, y, z)\, dy + R(x, y, z)\, dz$$

is an *exact differential* if there exists a function $u(x,y,z)$, differentiable in $\Omega$, such that

$$du = P(x, y, z)\, dx + Q(x, y, z)\, dy + R(x, y, z)\, dz \tag{5}$$

Equivalently, since

$$du = (\partial u/\partial x)\, dx + (\partial u/\partial y)\, dy + (\partial u/\partial z)\, dz,$$

the function $u$ will be a solution of the system of PDEs

$$\frac{\partial u}{\partial x} = P(x, y, z) \ , \quad \frac{\partial u}{\partial y} = Q(x, y, z) \ , \quad \frac{\partial u}{\partial z} = R(x, y, z) \tag{6}$$

The *integrability* (*compatibility*) *conditions* of the system (*necessary* conditions for existence of solution for $u$) are

$$\frac{\partial P}{\partial y} = \frac{\partial Q}{\partial x} \ , \quad \frac{\partial P}{\partial z} = \frac{\partial R}{\partial x} \ , \quad \frac{\partial Q}{\partial z} = \frac{\partial R}{\partial y} \tag{7}$$

Conditions (7) are also *sufficient* for solution if the domain $\Omega$, within which the functions $P$, $Q$, $R$ are differentiable, is *simply connected* [1].

***Example:*** Consider the differential relation

$$du = (x+y+z)(dx+dy+dz),$$

with $P=Q=R=x+y+z$. We notice that relations (7) are satisfied, as well as that the functions $P$, $Q$, $R$ are differentiable in the entire $R^3$, which is a simply connected space. Thus, the given differential relation admits a solution for $u$. Relations (6) are written

$$\partial u/\partial x = x+y+z \ , \quad \partial u/\partial y = x+y+z \ , \quad \partial u/\partial z = x+y+z.$$

The first one yields

$$u = (x^2/2) + xy + xz + \varphi(y,z) \quad \text{(arbitrary } \varphi\text{)}.$$

Substituting this into the second relation, we find

$$\partial\varphi/\partial y = y+z \ \Rightarrow\ \varphi(y,z) = (y^2/2) + yz + \psi(z) \quad \text{(arbitrary } \psi\text{)}.$$

Making the necessary replacements into the third relation, we have:

$$\psi'(z) = z \ \Rightarrow\ \psi(z) = (z^2/2) + C.$$

Finally,

$$u = (x^2+y^2+z^2)\,/\,2 + xy + xz + yz + C.$$

## 1.3 Line Integrals and Path Independence

Consider the plane $R^2$ with coordinates $(x,y)$. Let $L$ be an *oriented curve* (*path*) on the plane, with initial point $A$ and final point $B$ (Fig. 1.3). The curve $L$ may be described by parametric equations of the form

$$\{\ x=x(t)\ ,\ \ y=y(t)\ \} \tag{1}$$

Eliminating $t$ between these equations, we get a relation of the form $F(x,y)=0$ which, in certain cases, may be written in the form of a function $y=y(x)$.

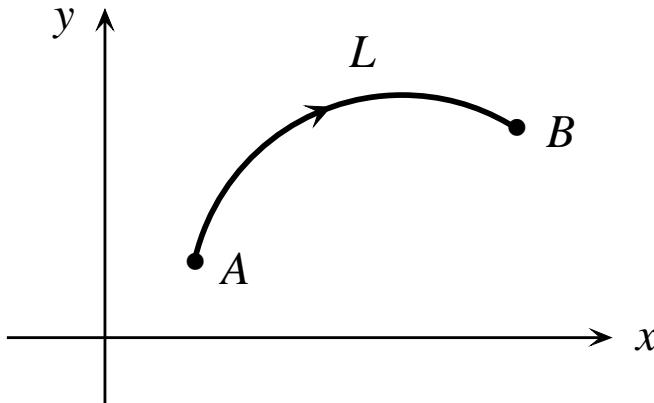

Fig. 1.3. An oriented curve on the $xy$-plane.

***Example:*** Consider the parametric curve of Fig. 1.4:

$$\{\ x = R\cos t\ ,\quad y = R\sin t\ \}\ ,\quad 0 \le t \le \pi\ .$$

The orientation of the curve depends on whether $t$ increases ("counterclockwise") or decreases ("clockwise") between 0 and $\pi$. By eliminating $t$, we get

$$x^2 + y^2 - R^2 = 0 \ \ \Rightarrow\ \ y = (R^2 - x^2)^{1/2}.$$

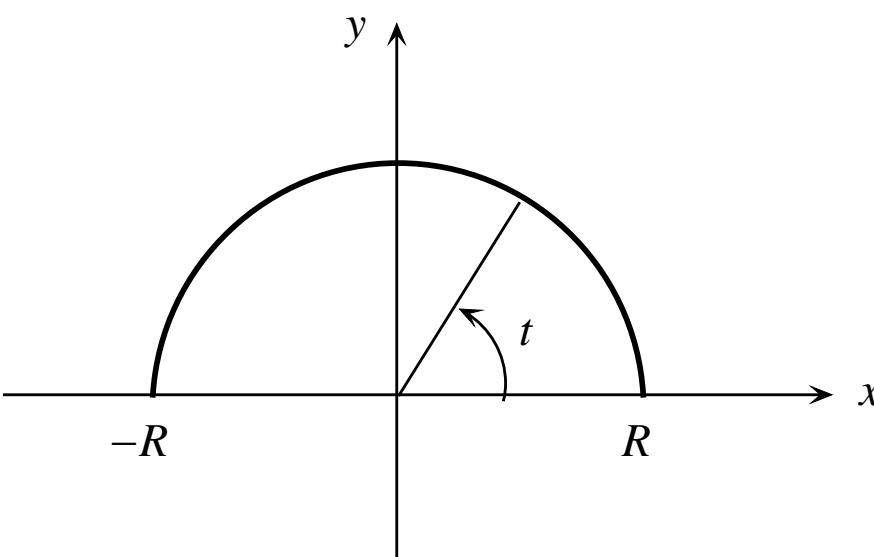

Fig. 1.4. A semicircle on the $xy$-plane.

Given a plane curve $L$ from $A$ to $B$, we now consider a *line integral* of the form

$$I_L = \int_L P(x,y)\,dx + Q(x,y)\,dy \tag{2}$$

In the parametric form (1) of $L$, we have

$$dx = (dx/dt)\,dt = x'(t)\,dt\ ,\quad dy = y'(t)\,dt,$$

so that

$$I_L = \int_{t_A}^{t_B} \{P[x(t), y(t)]\, x'(t) + Q[x(t), y(t)]\, y'(t)\}\, dt \qquad (3)$$

In the form $y=y(x)$ of $L$, we write $dy=y'(x)dx$ and

$$I_L = \int_{x_A}^{x_B} \{P[x, y(x)] + Q[x, y(x)]\, y'(x)\}\, dx \qquad (4)$$

In general, the value of the integral $I_L$ depends on the path $L$ connecting $A$ and $B$.

For every path $L$: $A \to B$, we can define the path $-L$: $B \to A$, with opposite orientation. From (3) it follows that, if

$$I_L = \int_{t_A}^{t_B} (\cdots)\, dt \,,$$

then

$$I_{-L} = \int_{t_B}^{t_A} (\cdots)\, dt \,.$$

Thus,

$$I_{-L} = -\, I_L \qquad (5)$$

If the end points $A$ and $B$ of a path coincide, then we have a *closed curve* $C$ and, correspondingly, a *closed line integral* $I_C$ , for which we use the symbol $\oint_C$. We then have:

$$\oint_{-C} (\cdots) = -\oint_C (\cdots) \qquad (6)$$

where the orientation of $-C$ is *opposite* to that of $C$ (e.g., if $C$ is counterclockwise on the plane, then $-C$ is clockwise).

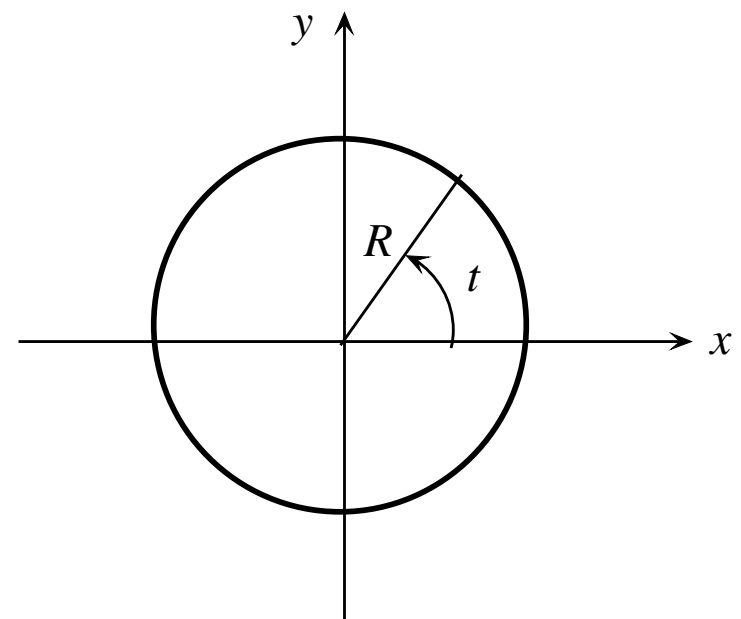

Fig. 1.5. A circle on the *xy*-plane.

***Example:*** The parametric curve

$$\{\ x = R\cos t\,,\quad y = R\sin t\ \}\ ,\quad 0 \le t \le 2\pi$$

represents a circle on the plane (Fig. 1.5). If the *counterclockwise* orientation of the circle (where *t increases* from 0 to $2\pi$) corresponds to the curve $C$, then the *clockwise* orientation (with *t decreasing* from $2\pi$ to 0) corresponds to the curve $-C$.

***Proposition:*** If

$$\oint_C Pdx + Qdy = 0$$

for *every* closed curve $C$ on a plane, then the line integral

$$\int_L Pdx + Qdy$$

is *independent of the path* $L$ connecting any two points $A$ and $B$ on this plane. The converse is also true.

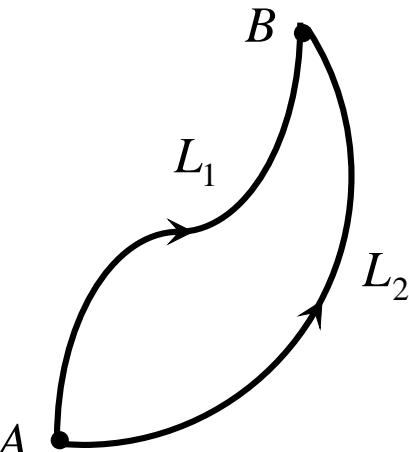

Fig. 1.6. Two paths connecting points $A$ and $B$ on a plane.

***Proof:*** We consider any two points $A$ and $B$ on the plane, as well as two different paths $L_1$ and $L_2$ connecting these points, as seen in Fig. 1.6 (there is an infinite number of such paths). We form the closed path $C=L_1+(-L_2)$ from $A$ to $B$ through $L_1$ and back again to $A$ through $-L_2$ . We then have:

$$\oint_C Pdx + Qdy = 0 \iff \int_{L_1} Pdx + Qdy + \int_{-L_2} Pdx + Qdy = 0 \iff$$

$$\int_{L_1} Pdx + Qdy - \int_{L_2} Pdx + Qdy = 0 \iff \int_{L_1} Pdx + Qdy = \int_{L_2} Pdx + Qdy \,.$$

***Theorem 1:*** Consider two functions $P(x,y)$ and $Q(x,y)$, differentiable in a *simply connected* domain $D$ of the plane. Then, the following 4 conditions are equivalent to one another (if any one is true, then the rest are true as well):

($a$) $\oint_C Pdx + Qdy = 0$, for any closed curve $C$ within $D$.

($b$) The integral $\int_L Pdx + Qdy$ is independent of the curved path $L$ connecting two fixed points $A$ and $B$ of $D$.

($c$) The expression $Pdx+Qdy$ is an *exact differential*. That is, there exists a function $u(x,y)$ such that

$$du = Pdx + Qdy \iff \partial u/\partial x = P \,, \quad \partial u/\partial y = Q \,.$$

($d$) At every point of $D$,

$$\frac{\partial P}{\partial y}=\frac{\partial Q}{\partial x} \quad .$$

(For a proof of this theorem, see, e.g., [1].)

***Comment:*** In the case where the domain $D$ is *not* simply connected, condition ($d$) does not guarantee the validity of the remaining three conditions. However, conditions ($a$), ($b$), ($c$) are still equivalent to one another and each of them separately guarantees ($d$). Note that ($d$) is the *integrability condition* for the validity of ($c$). (Remember that the former condition is *necessary* but not sufficient in the case where the domain $D$, in which the functions $P$ and $Q$ are differentiable, is not simply connected.)

***Example:*** Consider the differential expression

$$\omega=\frac{-y dx+x dy}{x^2+y^2} \quad .$$

Here, $P= -y/(x^2+y^2)$, $Q= x/(x^2+y^2)$, and condition ($d$) is satisfied (show this!). We notice that the functions $P$ and $Q$ are differentiable everywhere on the plane *except* at the origin $O$ of our coordinate system, at which point $(x,y)\equiv(0,0)$. Now, every domain $D$ of the plane *not containing point O* (left diagram in Fig. 1.7) is simply connected (explain why!). A closed curve $C$ within $D$ will then not contain point $O$ in its interior. For such a curve, $\oint_C \omega=0$ and $\omega=du$, with $u(x, y)=\arctan(y/x)$. A curve $C$, however, *containing O* (right diagram) cannot belong to a simply connected domain (why?). For such a curve, $\oint_C \omega \neq 0$.

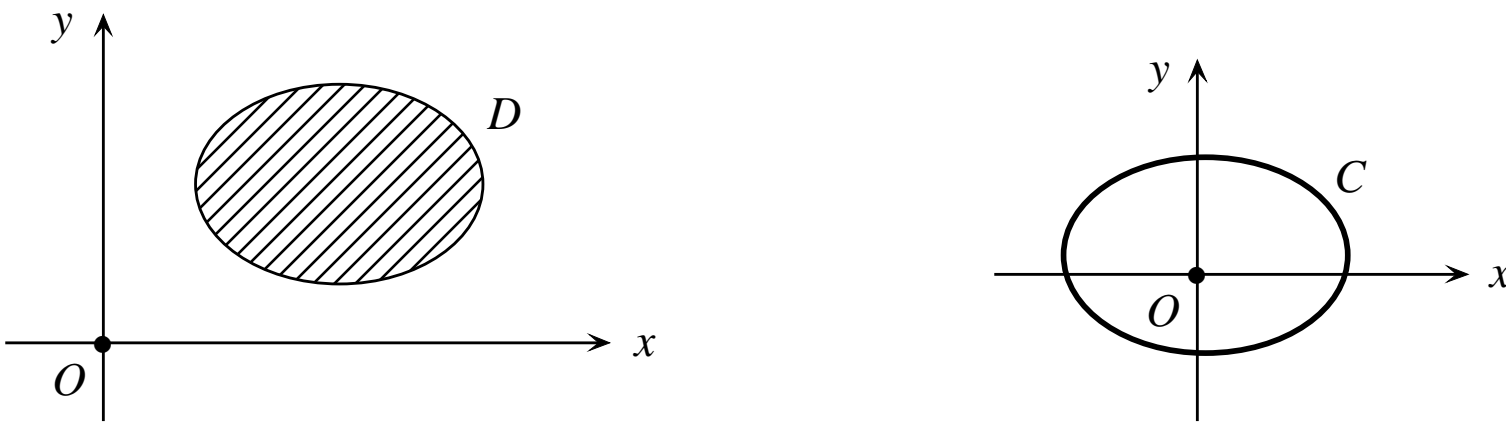

Fig. 1.7. The domain $D$ does not contain the origin $O$. The curve $C$ contains $O$ in its interior.

Let us now consider *line integrals in space*. Let

$$L: \{x=x(t)\ ,\ \ y=y(t)\ ,\ \ z=z(t)\}$$

be a parametric curve from point $A$ of $R^3$ to point $B$. Let $P(x,y,z)$, $Q(x,y,z)$, $R(x,y,z)$ be functions differentiable in the domain $\Omega\subseteq R^3$ in which the curve $L$ is located. We consider the line integral

$$I_L=\int_L P(x,y,z)\,dx+Q(x,y,z)\,dy+R(x,y,z)\,dz \qquad (7)$$

or, in parametric form,

$$I_L = \int_{t_A}^{t_B} dt\,\{P[x(t),y(t),z(t)]\,x'(t) + Q[x(t),y(t),z(t)]\,y'(t) + R[x(t),y(t),z(t)]\,z'(t)\} \tag{8}$$

***Theorem 2:*** If the domain *Ω* is *simply connected*, then the following 4 conditions are equivalent to one another:

(*a*) $\oint_C Pdx + Qdy + Rdz = 0$, for any closed path *C* within *Ω*.

(*b*) The integral $\int_L Pdx + Qdy + Rdz$ is independent of the curved path *L* connecting two fixed points *A* and *B* of *Ω*.

(*c*) The expression $Pdx+Qdy+Rdz$ is an *exact differential*. That is, there exists a function $u(x,y,z)$ such that

$$du = Pdx + Qdy + Rdz \quad \Leftrightarrow \quad \partial u/\partial x = P\,, \quad \partial u/\partial y = Q\,, \quad \partial u/\partial z = R\,.$$

(*d*) At every point of *Ω*,

$$\frac{\partial P}{\partial y} = \frac{\partial Q}{\partial x}\,, \quad \frac{\partial P}{\partial z} = \frac{\partial R}{\partial x}\,, \quad \frac{\partial Q}{\partial z} = \frac{\partial R}{\partial y}\,.$$

(See, e.g., [1] for a proof of the theorem.)

***Comment:*** If the domain *Ω* is *not* simply connected, condition (*d*) does not guarantee the validity of the remaining three conditions. However, conditions (*a*), (*b*), (*c*) are still equivalent to one another and each of them separately guarantees (*d*). Note that (*d*) is the integrability condition for the validity of (*c*). (As we have said, the former condition is *necessary* but not sufficient in the case where the domain *Ω*, in which the functions *P*, *Q* and *R* are differentiable, is not simply connected.)

From (*c*) it follows that, for any open curve *L* limited by two fixed points *A* and *B*,

$$\int_A^B Pdx + Qdy + Rdz = \int_A^B du = u(B) - u(A) \equiv u(x_B, y_B, z_B) - u(x_A, y_A, z_A) \tag{9}$$

Notice, in particular, that this automatically verifies (*b*).

## 1.4 Potential Vector Fields

Consider a vector field in a domain $\Omega \subseteq R^3$:

$$\vec{A}(\vec{r}) = P(x,y,z)\,\hat{u}_x + Q(x,y,z)\,\hat{u}_y + R(x,y,z)\,\hat{u}_z \tag{1}$$

where $\vec{r}$ is the position vector of a point $(x,y,z)$ of the domain, and where $\hat{u}_x, \hat{u}_y, \hat{u}_z$ are the unit vectors on the axes *x*, *y*, *z*, respectively. The functions *P*, *Q*, *R* are as-sumed to be differentiable in the domain *Ω*. We write: $\vec{A} \equiv (P, Q, R)$, $\vec{r} \equiv (x, y, z)$ and $\overrightarrow{dr} \equiv (dx, dy, dz)$.

We say that the field (1) is *potential* if there exists a differentiable function $u(x,y,z)$ such that

$$\vec{A} = \vec{\nabla} u \tag{2}$$

or, in components,

$$P = \frac{\partial u}{\partial x} \quad , \quad Q = \frac{\partial u}{\partial y} \quad , \quad R = \frac{\partial u}{\partial z} \tag{3}$$

The function $u$ is called *potential function* or simply *potential* of the field $\vec{A}$.

If (2) is valid, then

$$\vec{\nabla} \times \vec{A} = \vec{\nabla} \times \vec{\nabla} u = 0 \tag{4}$$

That is, *a potential vector field is necessarily irrotational*. In component form, Eq. (4) is written

$$\frac{\partial P}{\partial y} = \frac{\partial Q}{\partial x} \quad , \quad \frac{\partial P}{\partial z} = \frac{\partial R}{\partial x} \quad , \quad \frac{\partial Q}{\partial z} = \frac{\partial R}{\partial y} \tag{5}$$

Relations (5) are the *integrability conditions* for existence of a solution *u* of the system (3), thus also of the vector relation (2).

Condition (4) is *necessary* in order for the field $\vec{A}$ to be potential. Is it sufficient also? That is, is an irrotational field potential?

***Proposition:*** An *irrotational* field $\vec{A} \equiv (P, Q, R)$ in a *simply connected* domain *Ω* is potential.

***Proof:*** By assumption, the system of PDEs (5) is satisfied at every point of a simply connected domain. Hence, according to Theorem 2 of Sec. 1.3, the expression $Pdx+Qdy+Rdz$ is an exact differential. That is, there exists a function $u(x,y,z)$ such that

$$Pdx + Qdy + Rdz = du \tag{6}$$

Taking into account the independence of the differentials *dx*, *dy*, *dz*, we are thus led to the system (3), thus to the vector equation (2).

Theorem 2 of Sec. 1.3 can be re-expressed in the "language" of vector fields as follows [1,2]:

***Theorem:*** Consider a vector field $\vec{A} \equiv (P, Q, R)$, where the functions $P$, $Q$, $R$ are differentiable in a *simply connected* domain $\Omega \subseteq R^3$. Let $L$ be an open curve and let $C$ be a closed curve, both lying in $\Omega$. Then, the following 4 conditions are equivalent to one another:

($a$) $\oint_C \vec{A} \cdot \overrightarrow{dr} \equiv \oint_C P dx + Q dy + R dz = 0$.

($b$) The integral $\int_L \vec{A} \cdot \overrightarrow{dr} \equiv \int_L P dx + Q dy + R dz$ is independent of the curved path $L$ connecting any two fixed points in $\Omega$.

($c$) There exists a function $u(x, y, z)$ such that, at every point of $\Omega$, $\vec{A} = \vec{\nabla} u$.

($d$) At every point of $\Omega$, $\vec{\nabla} \times \vec{A} = 0$ (i.e., the field $\vec{A}$ is irrotational).

***Comments:***

***1.*** From ($c$) we have

$$\vec{A} \cdot \overrightarrow{dr} = \vec{\nabla} u \cdot \overrightarrow{dr} = du .$$

Thus, if $L$ is a curved path with limit points $a$ and $b$,

$$\int_L \vec{A} \cdot \overrightarrow{dr} = \int_a^b du = u(b) - u(a)$$

independently of the path $a \to b$, in accordance with condition ($b$).

***2.*** Assume that the domain $\Omega$ in which condition ($d$) is valid is simply connected. Then, for every closed curve $C$ in $\Omega$ there exists an open surface $S$ bounded by $C$, as seen in Fig. 1.8. By Stokes' theorem we then have

$$\oint_C \vec{A} \cdot \overrightarrow{dr} = \int_S (\vec{\nabla} \times \vec{A}) \cdot \overrightarrow{da} = 0$$

in accordance with condition ($a$).

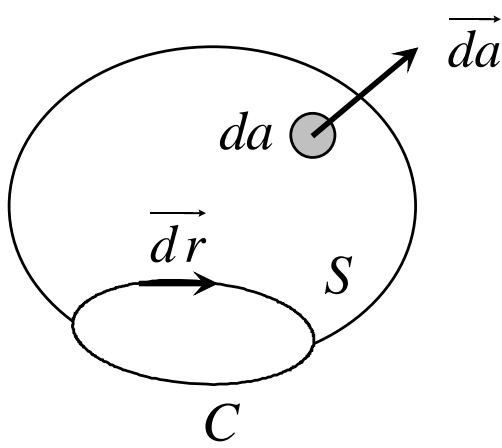

Fig. 1.8. An oriented open surface $S$ bordered by an oriented closed curve $C$. The *relative* orientation of $S$ and $C$ is in accordance with the familiar "right-hand rule".

***3.*** If the domain $\Omega$ is *not* simply connected, condition ($d$) does not guarantee the validity of the remaining three conditions. However, the first three conditions are still equivalent to one another and each of them separately guarantees ($d$).

## 1.5 Conservative Force Fields

In Physics, a *static* (time-independent) force field $\vec{F}(\vec{r})$ is called *conservative* if its work $W_{AB}$ on a test particle moving from point $A$ to point $B$ is independent of the path connecting these points. Equivalently, the work on the particle along a *closed* path $C$ is zero:

$$W_{AB}=\int_A^B \vec{F}\cdot\overrightarrow{dr} \quad \text{is independent of the path } A\to B \quad \Leftrightarrow \quad \oint_C \vec{F}\cdot\overrightarrow{dr}=0 \tag{1}$$

Let $S$ be an open surface inside the field, bounded by the closed curve $C$ (cf. Fig. 1.8). By Stokes' theorem, relation (1) yields

$$\oint_C \vec{F}\cdot\overrightarrow{dr}=\int_S (\vec{\nabla}\times\vec{F})\cdot\overrightarrow{da}=0\,.$$

In order for this to be valid for every open surface bounded by $C$, we must have

$$\vec{\nabla}\times\vec{F}=0 \tag{2}$$

That is, *a conservative force field is irrotational.* (The validity of the converse requires that the domain of space in which the field is defined be *simply connected*.)

From (1) it also follows that, according to the Theorem of Sec. 1.4, there exists a function such that $\vec{F}(\vec{r})$ is the *grad* of this function. We write

$$\vec{F}=-\vec{\nabla}U \tag{3}$$

The function $U(\vec{r})=U(x,y,z)$ is called the *potential energy* of the test particle at the point $\vec{r}\equiv(x,y,z)$ of the field. [The negative sign in (3) is only a matter of convention and has no special physical meaning. One may eliminate it by putting $-U$ in place of $U$. Note also that $U$ is arbitrary to within an additive constant, given that $U$ and $(U+c)$ correspond to the same force $\vec{F}$ in (3).]

The work $W_{AB}$ of $\vec{F}$ (where the latter is assumed to be the total force on the particle) is written

$$W_{AB}=\int_A^B \vec{F}\cdot\overrightarrow{dr}=-\int_A^B(\vec{\nabla}U)\cdot\overrightarrow{dr}=-\int_A^B dU \quad\Rightarrow$$

$$W_{AB}=U(\vec{r}_A)-U(\vec{r}_B)\equiv U_A-U_B \tag{4}$$

Now, by the *work-energy theorem* [3-5]

$$W_{AB}=E_{k,B}-E_{k,A} \tag{5}$$

where $E_k=mv^2/2$ is the *kinetic energy* of the particle ($m$ and $v$ are the particle's mass and speed, respectively). By combining (4) with (5), we have:

$$E_{k,A} + U_A = E_{k,B} + U_B \tag{6}$$

The sum $E=E_k+U$ represents the *total mechanical energy* of the particle. Relation (6), then, expresses the *principle of conservation of mechanical energy*; namely, the total mechanical energy of a particle moving inside a conservative force field assumes a constant value during the motion of the particle.

***Example:*** We consider the *electrostatic Coulomb field* due to a point charge $Q$ located at the origin $O$ of our coordinate system (see Fig. 1.9). Let $q$ be a test charge at a field point with position vector $\vec{r} = r\hat{r}$, where $r$ is the distance of $q$ from $O$ and where $\hat{r}$ is the unit vector in the direction of $\vec{r}$ [in this problem it is convenient to use spherical coordinates ($r$, $\theta$, $\varphi$)].

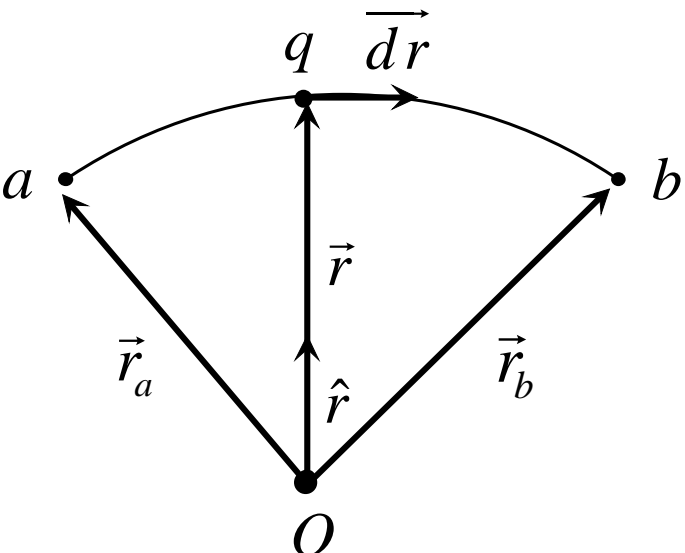

Fig. 1.9. A test charge $q$ inside the Coulomb field produced by a point charge $Q$ (not shown) located at $O$.

The force exerted on $q$ by the field, at a momentary position $\vec{r}$ of the charge, is

$$\vec{F} = \frac{kqQ}{r^2}\,\hat{r}$$

(where $k$ is a constant that depends on the system of units). As can be shown, $\vec{\nabla}\times\vec{F} = 0$. That is, the force field $\vec{F}$ is *irrotational*. This field is defined in a simply connected domain of space (the whole space with the exception of the single point $O$ where the charge $Q$ producing the electrostatic field is located). Hence, the considered irrotational force field will also be *conservative*. Indeed, putting $F(r)=kqQ/r^2$ , we write

$$\vec{F} = F(r)\,\hat{r} = \frac{F(r)}{r}\,\vec{r} \quad\Rightarrow\quad \vec{F}\cdot\overrightarrow{dr} = \frac{F(r)}{r}\,\vec{r}\cdot\overrightarrow{dr} \ .$$

But,

$$\vec{r}\cdot\overrightarrow{dr} = \frac{1}{2}d(\vec{r}\cdot\vec{r}) = \frac{1}{2}d(r^2) = r\,dr \ ,$$

so that $\vec{F}\cdot\overrightarrow{dr} = F(r)\,dr$. Thus, the work produced during the motion of $q$ from a field point $a$ to a field point $b$ is

$$W_{ab} = \int_a^b \vec{F}\cdot\overrightarrow{dr} = \int_a^b F(r)\,dr = kqQ\left(\frac{1}{r_a} - \frac{1}{r_b}\right) \ .$$

This expression allows us to define the potential energy $U(r)$ of $q$ at a given field point by using (4): $W_{ab}=U_a-U_b$ . As is easy to see,

$$U(r) = \frac{kqQ}{r} \quad (+\,const.) \ .$$

We notice that

$$-\vec{\nabla}U = -\frac{\partial U}{\partial r}\,\hat{r} = \frac{kqQ}{r^2}\,\hat{r} = \vec{F} \ .$$

The total mechanical energy of the charge $q$ remains fixed during the motion of the charge inside the field, and is equal to

$$E = E_k + U(r) = mv^2/2 + kqQ/r = const.$$

***Comment:*** In the beginning of this section we stated that the conservative force field $\vec{F}$ is static. Indeed, a time-dependent force *cannot* be conservative! This is explained in Appendix A.

## CHAPTER 2

# INTEGRABILITY ON THE COMPLEX PLANE

## 2.1 Analytic Functions

We consider complex functions of the form

$$w = f(z) = u(x, y) + i\, v(x, y) \qquad (1)$$

where $z=x+iy \equiv (x, y)$ is a point on the complex plane. Let $\Delta z$ be a change of $z$ and let $\Delta w = f(z+\Delta z) - f(z)$ be the corresponding change of the value of $f(z)$. We say that the function (1) is *differentiable* at point $z$ if we can write

$$\frac{\Delta w}{\Delta z} = f'(z) + \varepsilon(z, \Delta z) \quad \text{with} \quad \lim_{\Delta z \to 0} \varepsilon(z, \Delta z) = 0 \qquad (2)$$

Then, the function

$$f'(z) = \lim_{\Delta z \to 0} \frac{\Delta w}{\Delta z} = \lim_{\Delta z \to 0} \frac{f(z+\Delta z) - f(z)}{\Delta z} \qquad (3)$$

is the *derivative* of $f(z)$ at point $z$. Evidently, in order for $f(z)$ to be differentiable at $z$, this function must be *defined* at that point. We also note that a function differentiable at a point $z_0$ is necessarily *continuous* at $z_0$ (the converse is not always true) [1-3]. That is, $\lim_{z \to z_0} f(z) = f(z_0)$ (assuming that the limit exists).

A function $f(z)$ differentiable at every point of a domain $G$ of the complex plane is said to be *analytic* (or *holomorphic*) in the domain $G$. The criterion for analyticity is the validity of a pair of partial differential equations (PDEs) called the *Cauchy-Riemann relations*.

***Theorem:*** Consider a complex function $f(z)$ of the form (1), continuous at every point $z \equiv (x, y)$ of a domain $G$ of the complex plane. The real functions $u(x,y)$ and $v(x,y)$ are differentiable at every point of $G$ and, moreover, their partial derivatives with respect to $x$ and $y$ are continuous functions in $G$. Then, the function $f(z)$ is analytic in the domain $G$ if and only if the following system of PDEs is satisfied [1-3]:

$$\frac{\partial u}{\partial x} = \frac{\partial v}{\partial y} \quad , \quad \frac{\partial u}{\partial y} = -\frac{\partial v}{\partial x} \qquad (4)$$

It is convenient to use the following notation for partial derivatives:

$$\frac{\partial \phi}{\partial x} \equiv \phi_x \ , \quad \frac{\partial \phi}{\partial y} \equiv \phi_y \ , \quad \frac{\partial^2 \phi}{\partial x^2} \equiv \phi_{xx} \ , \quad \frac{\partial^2 \phi}{\partial y^2} \equiv \phi_{yy} \ , \quad \frac{\partial^2 \phi}{\partial x \partial y} \equiv \phi_{xy} \ , \quad \text{etc.}$$

The Cauchy-Riemann relations (4) then read

$$u_x = v_y \quad , \qquad u_y = -\, v_x \tag{4´}$$

The derivative of the function (1) may now be expressed in the following alternate forms:

$$f\,'(z) = u_x + i v_x = v_y - i u_y = u_x - i u_y = v_y + i v_x \tag{5}$$

***Comments:***

***1.*** Relations (4) allow us to find $v$ when we know $u$, and vice versa. Let us put $u_x=P$, $u_y=Q$, so that $\{v_x= -\,Q \;,\; v_y=P\}$. The *integrability* (*compatibility*) *condition* of this system for solution for $v$, for a given $u$, is

$$\partial P/\partial x = -\,\partial Q/\partial y \;\Rightarrow\; u_{xx}+u_{yy}= 0.$$

Similarly, the integrability condition of system (4) for solution for $u$, for a given $v$, is $v_{xx}+v_{yy}= 0$. We notice that both the real and the imaginary part of an analytic function are *harmonic functions*, i.e., they satisfy the *Laplace equation*

$$w_{xx} + w_{yy} = 0 \tag{6}$$

Harmonic functions related to each other by means of the Cauchy-Riemann relations (4) are called *conjugate harmonic*.

***2.*** Let $z^*=x-iy$ be the complex conjugate of $z=x+iy$. Then,

$$x = (z+z^*)\,/\,2 \quad , \qquad y = (z-z^*)\,/\,2i \tag{7}$$

By using relations (7) we can express $u(x,y)$ and $v(x,y)$, thus also the sum $w=u+iv$, as functions of $z$ and $z^*$. The real Cauchy-Riemann relations (4), then, are rewritten in the form of a single complex equation [1-3]

$$\partial w \,/\, \partial z^* = 0 \tag{8}$$

One way to interpret this result is the following: The analytic function (1) is *literally* a function of the complex variable $z=x+iy$, not just some complex function of two real variables $x$ and $y$!

***Examples:***

***1.*** We seek an analytic function of the form (1), with $v(x,y)=xy$. Note first that $v$ satisfies the PDE (6): $v_{xx}+v_{yy}=0$ (harmonic function). Thus, the integrability condition of the system (4) for solution for $u$ is satisfied. The system is written

$$\partial u/\partial x=x \;,\;\; \partial u/\partial y= -y \;.$$

The first relation yields

$$u = x^2/2 + \varphi(y) .$$

From the second one we then get

$$\varphi'(y) = -y \quad \Rightarrow \quad \varphi(y) = -y^2/2 + C$$

so that

$$u = (x^2 - y^2)/2 + C .$$

Putting $C$=0, we finally have

$$w = u + iv = (x^2 - y^2)/2 + ixy .$$

***Exercise:*** Using relations (7), show that $w=f(z)=z^2/2$, thus verifying condition (8).

***2.*** Consider the function $w=f(z)=|z|^2$ defined on the entire complex plane. Here, $u(x,y)=x^2+y^2$, $v(x,y)=0$. As is easy to verify, the Cauchy-Riemann relations (4) are not satisfied anywhere on the plane, except at the single point $z$=0 where $(x,y)\equiv(0,0)$. Alternatively, we may write $w=zz^*$, so that $\partial w/\partial z^*=z \neq 0$ (except at $z$=0). We conclude that the given function is not analytic on the complex plane.

## 2.2 Integrals of Complex Functions

Let $L$ be an oriented curve on the complex plane (Fig. 2.1) the points of which plane are represented as $z=x+iy \equiv (x, y)$. The points $z$ of $L$ are determined by some parametric relation of the form

$$z = \lambda(t) = x(t) + i\, y(t) \quad , \qquad \alpha \leq t \leq \beta \tag{1}$$

As $t$ increases from $\alpha$ to $\beta$, the curve $L$ is traced from $A$ to $B$, while the opposite curve $-L$ is traced from $B$ to $A$ with $t$ decreasing from $\beta$ to $\alpha$.

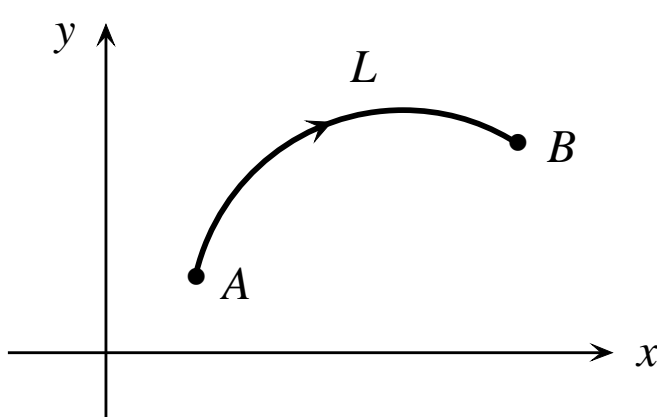

Fig. 2.1. An oriented curve on the complex plane.

We now consider integrals of the form $\int_L f(z)\,dz$, where $f(z)$ is a complex function. We write $dz=\lambda'(t)\,dt$, so that

$$\int_L f(z)\,dz = \int_\alpha^\beta f[\lambda(t)]\,\lambda'(t)\,dt \tag{2}$$

Also,

$$\int_{-L} f(z)\,dz \;=\; \int_{\beta}^{\alpha}(\cdots)\,dt \;=\; -\int_{\alpha}^{\beta}(\cdots)\,dt \quad \Rightarrow$$

$$\int_{-L} f(z)\,dz \;=\; -\int_{L} f(z)\,dz \tag{3}$$

A *closed* curve $C$ will be conventionally regarded as *positively* oriented if it is traced *counterclockwise*. Then, the opposite curve $-C$ will be *negatively* oriented and will be traced *clockwise*. Moreover,

$$\oint_{-C} f(z)\,dz \;=\; -\oint_{C} f(z)\,dz \tag{4}$$

***Examples:***

***1.*** We want to evaluate the integral

$$I = \oint_{|z-a|=\rho} \frac{dz}{z-a} \;,$$

where the circle $|z-a|=\rho$ is traced ($a$) counterclockwise, ($b$) clockwise.

($a$) The circle $|z-a|=\rho$ is described parametrically by the relation

$$z=a+\rho e^{it}\,,\;\; 0\leq t\leq 2\pi\;.$$

Then,

$$dz=(a+\rho e^{it})'\,dt = i\rho e^{it}\,dt\;.$$

Integrating from 0 to $2\pi$ (for *counterclockwise* tracing) we have:

$$I=\int_0^{2\pi}\frac{i\rho e^{it}dt}{\rho e^{it}} \;=\; i\int_0^{2\pi}dt \quad\Rightarrow\quad \oint_{|z-a|=\rho}\frac{dz}{z-a} \;=\; 2\pi i\;.$$

($b$) For *clockwise* tracing of the circle $|z-a|=\rho$, we write, again,

$$z=a+\rho e^{it}\;\;(0\leq t\leq 2\pi)\;.$$

This time, however, we integrate from $2\pi$ to 0. Then,

$$I = i\int_{2\pi}^{0}dt = -2\pi i\;.$$

Alternatively, we write

$$z=a+\rho e^{-it}\;\;(0\leq t\leq 2\pi)$$

and integrate from 0 to $2\pi$, arriving at the same result.

***2.*** Consider the integral

$$I = \oint_{|z-a|=\rho} \frac{dz}{(z-a)^2} \;,$$

where the circle $|z-a|=\rho$ is traced *counterclockwise*. We write

$$z=a+\rho e^{it} \quad (0\leq t \leq 2\pi)$$

so that

$$I=\int_0^{2\pi}\frac{i\rho e^{it}dt}{\rho^2 e^{2it}} = \frac{i}{\rho}\int_0^{2\pi} e^{-it}dt = 0 \ .$$

In general, for $k = 0, \pm1, \pm2, \ldots$ and for a *positively* (*counterclockwise*) oriented circle $|z-a|=\rho$ , one can show that

$$\oint_{|z-a|=\rho} \frac{dz}{(z-a)^k} = \begin{cases} 2\pi i \ , & \text{if} \quad k=1 \\ 0 \ , & \text{if} \quad k\neq 1 \end{cases} \tag{5}$$

## 2.3 Some Basic Theorems

We now state some important theorems concerning analytic functions [1-3].

***Theorem 1 (Cauchy-Goursat):*** Assume that the function $f(z)=u(x,y)+i\,v(x,y)$ is analytic in a *simply connected* domain *G* of the complex plane. Then, for any closed curve *C* in *G*,

$$\oint_C f(z)\,dz = 0 \tag{1}$$

***Proof:*** Write $dz=dx+i\,dy$ , so that

$$f(z)\,dz = (udx - vdy) + i\,(vdx + udy) \ \Rightarrow$$
$$\oint_C f(z)\,dz = \oint_C (udx - vdy) + i\oint_C (vdx + udy)\ .$$

Now, given that $f(z)$ is analytic in *G*, the Cauchy-Riemann relations will be valid in this domain. Moreover, since *G* is simply connected, the conditions of validity of Theorem 1 of Sec. 1.3 are satisfied. Hence,

$$u_y = (-v)_x \ \Leftrightarrow\ \oint_C udx + (-v)dy = 0\ ,$$
$$v_y = u_x \ \Leftrightarrow\ \oint_C vdx + udy = 0\ ,$$

by which relations the result (1) follows immediately.

***Corollary:*** In a simply connected domain *G*, the line integral of an analytic function $f(z)$ is independent of the path connecting any two points *A* and *B*.

***Proof:*** As in Sec. 1.3, we consider two paths $L_1$ and $L_2$ (see Fig. 2.2) and we form the closed path $C=L_1+(-L_2)$. By (1) we then have:

$$\oint_C f(z)\,dz = \int_{L_1} f(z)\,dz + \int_{-L_2} f(z)\,dz = 0 \ \Leftrightarrow$$
$$\int_{L_1} f(z)\,dz - \int_{L_2} f(z)\,dz = 0 \ \Leftrightarrow\ \int_{L_1} f(z)\,dz = \int_{L_2} f(z)\,dz\ .$$

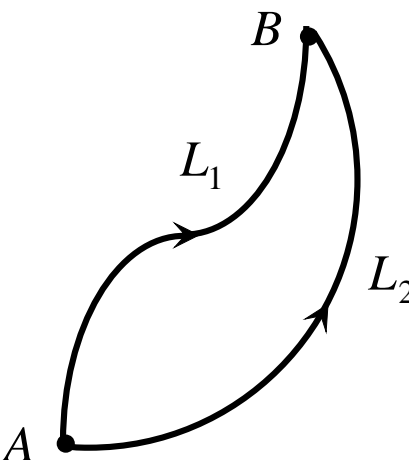

Fig. 2.2. Two paths connecting points $A$ and $B$ on the complex plane.

Let us assume, now, that the function $f(z)$ is analytic in a domain $G$ that is *not* simply connected. (For example, the domain $G$ in Fig. 2.3 is doubly connected.) Let $C$ be a closed curve in $G$. Two possibilities exist:

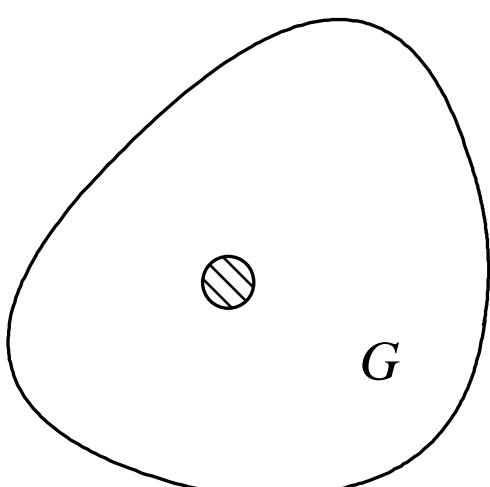

Fig. 2.3. A doubly connected domain $G$ on the complex plane.

($a$) The curve $C$ does not enclose any points not belonging to $G$. Then, $C$ may be considered as the boundary of a simply connected subdomain of $G$ where the conditions of validity of the Cauchy-Goursat theorem are fulfilled. Therefore,

$$\oint_C f(z)\,dz = 0 .$$

($b$) The curve $C$ encloses points not belonging to $G$. Then, $C$ may not belong to some simply connected subdomain of $G$ and the conditions of Theorem 1 are not fulfilled. In such a case, relation (1) may or may not be satisfied.

***Example:*** Let $G$ consist of the complex plane *without* the origin $O$ of its axes (i.e., without the point $z$=0). The function $f(z)=1/z$ is analytic in this domain. Let $C$ be the circle $|z|=\rho$ centered at $O$. Then, as we saw in Sec. 2.2,

$$\oint_C \frac{dz}{z} = 2\pi i \;\; (\neq 0) .$$

On the contrary,

$$\oint_C \frac{dz}{z^k} = 0 \;\text{ for } k\neq 1 .$$

***Theorem 2 (composite contour theorem):*** Consider a multiply (e.g., doubly, triply,...) connected domain $G$ of the complex plane (Fig. 2.4) and let $\Gamma$ be a closed curve in $G$. Let $\gamma_1$ , $\gamma_2$ ,..., $\gamma_n$ be closed curves in the interior of $\Gamma$ (but in the exterior of one another) such that the domain $D$ between the $\gamma_k$ and $\Gamma$ belongs entirely to $G$. Then, for every function $f(z)$ analytic in $G$,

$$\oint_\Gamma f(z)\,dz = \sum_{i=1}^{n} \oint_{\gamma_i} f(z)\,dz = \oint_{\gamma_1} f(z)\,dz + \oint_{\gamma_2} f(z)\,dz + \cdots + \oint_{\gamma_n} f(z)\,dz \qquad (2)$$

where all curves $\Gamma$, $\gamma_1$ ,..., $\gamma_n$ , are traced *in the same direction* (e.g., counterclockwise).

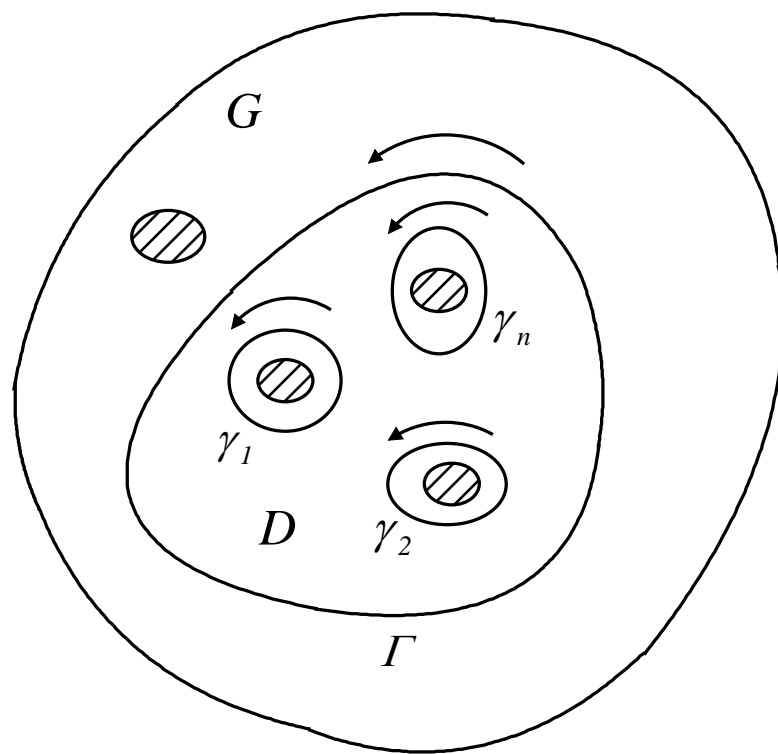

Fig. 2.4. A multiply connected domain $G$ on the complex plane.

***Problem:*** Show that

$$\oint_\Gamma \frac{dz}{z} = 2\pi i$$

where $\Gamma$ is any positively oriented, closed curve enclosing the origin $O$ ($z$=0) of the complex plane. (*Hint:* Consider a circle $\gamma : |z|=\rho$ centered at $O$ and lying in the interior of $\Gamma$.)

***Theorem 3 (Cauchy integral formula):*** Consider a function $f(z)$ analytic in a domain $G$ (Fig. 2.5). Let $C$ be a closed curve in $G$, such that the interior $D$ of $C$ belongs entirely to $G$. Consider a point $z_0 \in D$. Then,

$$f(z_0) = \frac{1}{2\pi i} \oint_C \frac{f(z)\,dz}{z - z_0} \qquad (3)$$

where $C$ is traced in the *positive* direction (i.e., *counterclockwise*).

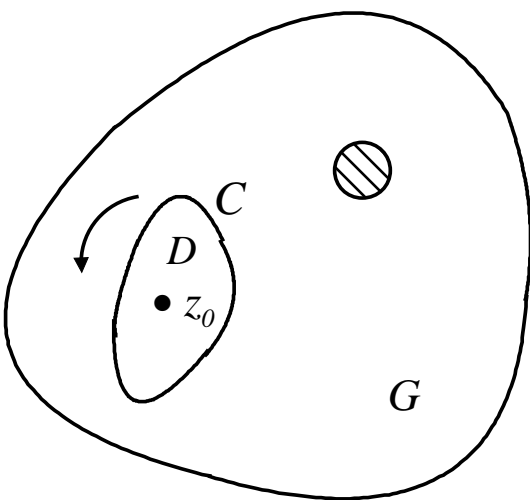

Fig. 2.5. The curve $C$ encloses a simply connected subdomain $D$ of $G$.

***Comments:***

***1.*** The value of the integral in (3) is independent of the choice of the curve $C$ that encloses $z_0$ and satisfies the conditions of the theorem. (This follows from the composite contour theorem for $n=1$.)

***2.*** More generally, we can write

$$\frac{1}{2\pi i}\oint_C \frac{f(z)\,dz}{z-z_0} = \begin{cases} f(z_0)\ , & \text{if} \quad z_0 \in D \\ 0\ , & \text{if} \quad z_0 \in (G-D) \end{cases} \quad .$$

Indeed: If $z_0 \in (G-D)$ (that is, $z_0 \notin D$), then the function $f(z)/(z-z_0)$ is analytic everywhere inside the simply connected domain $D$ and thus satisfies the Cauchy-Goursat theorem.

***Application:*** Putting $f(z)=1$ and considering a positively oriented path $C$ around a point $z_0$, we find

$$\oint_C \frac{dz}{z-z_0} = 2\pi i \quad .$$

More generally, for $k = 0, \pm 1, \pm 2, \ldots$ ,

$$\oint_C \frac{dz}{(z-z_0)^k} = \begin{cases} 2\pi i\ , & \text{if} \quad k=1 \\ 0\ , & \text{if} \quad k \neq 1 \end{cases} \tag{4}$$

where the point $z_0$ is located in the *interior* of $C$.

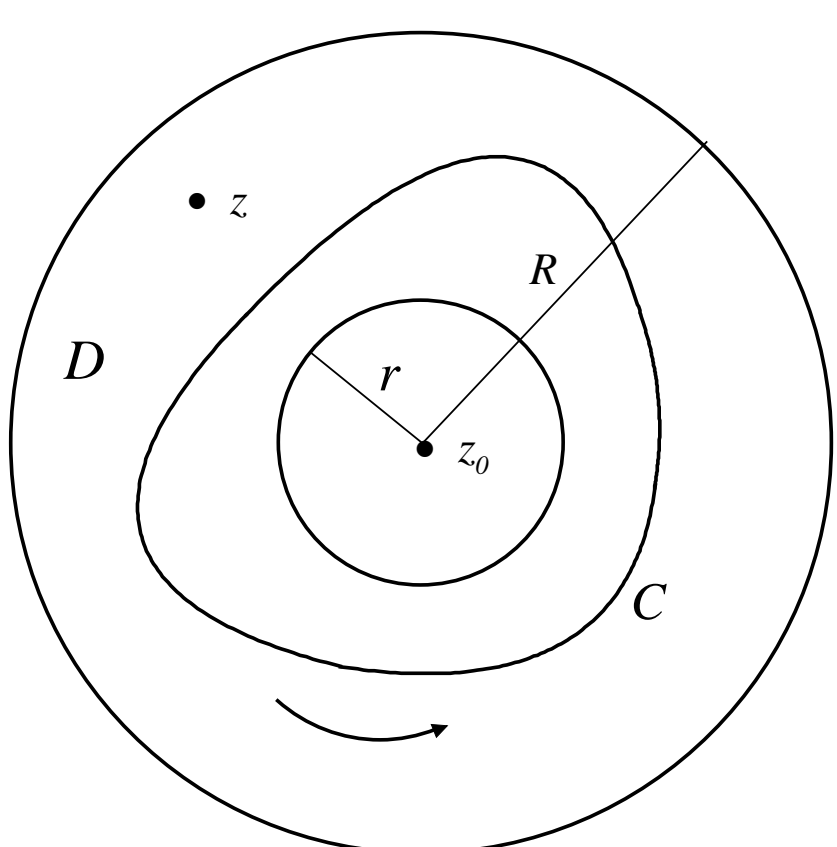

Fig. 2.6. An annular domain $D$ centered at $z_0$ .

***Theorem 4 (Laurent series):*** Consider a function $f(z)$, analytic in an annular domain $D$: $r < |z-z_0| < R$, centered at $z_0$ (Fig. 2.6). Let $C$ be a positively oriented, closed path around $z_0$ and inside the annulus $D$. Then, at every point $z \in D$ the function $f(z)$ may be represented by a convergent series of the form

$$f(z)=\sum_{n=-\infty}^{+\infty} a_n\,(z-z_0)^n \tag{5}$$

where the coefficients $a_n$ are given by

$$a_n=\frac{1}{2\pi i}\oint_C \frac{f(z)\,dz}{(z-z_0)^{n+1}} \tag{6}$$

and where the value of the integral in (6) is independent of the choice of the curve *C*.

***Proof of the coefficient formula:*** Assuming the validity of (5), we have:

$$\frac{f(z)}{(z-z_0)^{k+1}}=\sum_n a_n\,(z-z_0)^{n-k-1}\;\;\Rightarrow$$

$$\oint_C \frac{f(z)\,dz}{(z-z_0)^{k+1}}=\sum_n a_n \oint_C (z-z_0)^{n-k-1}\,dz\equiv\sum_n a_n\,I_{nk}\;.$$

But, by (4),

$$I_{nk}=\begin{Bmatrix}2\pi i\;, & \text{if}\;\; n=k\\ 0\;, & \text{if}\;\; n\neq k\end{Bmatrix}=2\pi i\,\delta_{nk}$$

where $\delta_{nk}$ is the "Kronecker delta", assuming the values 1 and 0 for $n=k$ and $n\neq k$, respectively. Hence,

$$\oint_C \frac{f(z)\,dz}{(z-z_0)^{k+1}}=2\pi i\sum_n a_n\,\delta_{nk}=2\pi i\,a_k\;\;\Rightarrow\;\; a_k=\frac{1}{2\pi i}\oint_C \frac{f(z)\,dz}{(z-z_0)^{k+1}}\;.$$

***Comment:*** The annulus $D$: $r<|z-z_0|<R$ may be

- the region between two concentric circles $(0<r<R)$;
- a circle with its center $z_0$ deleted $(r=0,\; R>0)$;
- the exterior of a circle $(r>0,\; R=\infty)$; or
- the entire complex plane with point $z_0$ deleted $(r=0,\; R=\infty)$.

## 2.4 Antiderivative and Indefinite Integral of an Analytic Function

Let $z_0$ and $z$ be two points in a simply connected domain *G* of the complex plane. We regard $z_0$ as constant while $z$ is assumed to be variable. According to the Cauchy-Goursat theorem, the line integral from $z_0$ to $z$, of a function $f(z)$ analytic in *G*, depends only on the two limit points and is independent of the curved path connecting them. Hence, such an integral may be denoted by

$$\int_{z_0}^{z} f(z')\,dz'$$

or, for simplicity,

$$\int_{z_0}^{z} f(z)\,dz\,.$$

For variable upper limit $z$, this integral is a function of its upper limit. We write

$$\int_{z_0}^{z} f(z)\,dz = I(z) \tag{1}$$

As can be shown [1] $I(z)$ is an analytic function. Moreover, it is an *antiderivative* of $f(z)$; that is, $I'(z)=f(z)$ . Analytically,

$$I'(z) = \frac{d}{dz}\int_{z_0}^{z} f(z)\,dz = f(z) \tag{2}$$

Any antiderivative $F(z)$ of $f(z)$ [$F'(z)=f(z)$] is equal to $F(z)=I(z)+C$, where $C=F(z_0)$ is a constant [note that $I(z_0)=0$]. We observe that $I(z)=F(z)-F(z_0)$ $\Rightarrow$

$$\int_{z_0}^{z} f(z)\,dz = F(z)-F(z_0) \tag{3}$$

In general, for given $z_1$ , $z_2$ and for an *arbitrary* antiderivative $F(z)$ of $f(z)$, we may write

$$\int_{z_1}^{z_2} f(z)\,dz = F(z_2)-F(z_1) \tag{4}$$

Now, if we also allow the lower limit $z_0$ of the integral in Eq. (1) to vary, then this relation yields an *infinite set of antiderivatives* of $f(z)$, which set represents the *indefinite integral* of $f(z)$ and is denoted by $\int f(z)\,dz$ . If $F(z)$ is any antiderivative of $f(z)$, then, by relation (3) and by putting $-F(z_0)=C$ ,

$$\int f(z)\,dz = \left\{ F(z)+C \;/\; F'(z)=f(z),\;\; C=const. \right\}.$$

To simplify our notation, we write

$$\int f(z)\,dz = F(z)+C \tag{5}$$

where the right-hand side represents an *infinite set* of functions, not just any specific antiderivative of $f(z)$!

***Examples:***

***1.*** The function $f(z)=z^2$ is analytic on the entire complex plane and one of its antiderivatives is $F(z)=z^3/3$ . Thus,

$$\int z^2 dz = \frac{z^3}{3} + C \qquad \text{and} \qquad \int_{-1}^{i} z^2 dz = \frac{1}{3}\,(1-i)\ .$$

***2.*** The function $f(z)=1/z^2$ is differentiable everywhere except at the origin $O$ of the complex plane, where $z=0$. An antiderivative, for $z\neq 0$, is $F(z)=-1/z$ . Hence,

$$\int \frac{dz}{z^2} = -\frac{1}{z} + C \qquad \text{and} \qquad \int_{z_1}^{z_2} \frac{dz}{z^2} = \frac{1}{z_1} - \frac{1}{z_2}$$

where the path connecting the points $z_1\neq 0$ and $z_2\neq 0$ does not pass through $O$.

## CHAPTER 3

# ORDINARY DIFFERENTIAL EQUATIONS

### 3.1 The Concept of the First Integral

An *ordinary differential equation* (ODE) is often easier to solve if we can find one or more *first integrals*. In simple terms, a first integral is a relation (algebraic or differential) that gives us the information that some mathematical quantity retains a constant value as a consequence of the given ODE. This quantity may contain the dependent variable $y$, the independent variable $x$, as well as derivatives $y^{(k)}(x) = d^k y/dx^k$.

When derivatives are contained, a first integral leads to an ODE of lower order than the initial ODE. Thus, by using a first integral one may *reduce the order* of a given ODE. If the ODE is of the first order, a first integral is an *algebraic* relation expressing the solution directly. In general, an ODE of order $n$ is completely solved if one manages to find *n independent first integrals*.

In Classical Mechanics one often needs to find the solution of a system of ODEs expressing Newton's second law of motion. With the exception of some simple cases, this system is hard to integrate; for this reason one seeks as many first integrals as possible. These quantities are called *constants of the motion* and they express corresponding *conservation laws*, such as, e.g., conservation of total mechanical energy, of momentum, or of angular momentum [1-3].

### 3.2 Exact Equations

Consider the first-order ODE

$$\frac{dy}{dx} = -\frac{M(x,y)}{N(x,y)} \quad (N \neq 0),$$

which is written more symmetrically as

$$M(x,y)\,dx + N(x,y)\,dy = 0 \tag{1}$$

Equation (1) is said to be *exact* if there exists a function $u(x,y)$ such that

$$M(x,y)\,dx + N(x,y)\,dy = du \tag{2}$$

(that is, if the expression $Mdx+Ndy$ is an exact differential). Then, by (1) and (2), $du=0 \Rightarrow$

$$u(x,y) = C \tag{3}$$

Equation (3) is an algebraic relation connecting $x$ and $y$ and containing an arbitrary constant. Thus it can be regarded as the general solution of (1). Relation (3) is a *first integral* of the ODE (1) and directly determines the general solution of this equation.

According to (2), the function $u(x,y)$ satisfies the following system of first-order partial differential equations (PDEs):

$$\frac{\partial u}{\partial x} = M(x, y) \ , \quad \frac{\partial u}{\partial y} = N(x, y) \tag{4}$$

The *integrability condition* of this system for existence of a solution for $u$, is

$$\frac{\partial M}{\partial y} = \frac{\partial N}{\partial x} \tag{5}$$

If condition (5) is valid at all points of a *simply connected* domain $D$ of the $xy$-plane, then it guarantees the existence of a solution for the system (4) or, equivalently, for the differential relation (2).

The constant $C$ in the solution (3) is determined by the *initial condition* of the problem: if the specific value $x=x_0$ corresponds to the value $y=y_0$ , then $C=C_0=u(x_0, y_0)$. We thus get the particular solution $u(x, y)=C_0$ .

***Example:*** We consider the ODE

$$(x+y+1)\, dx + (x-y^2+3)\, dy = 0 \ , \ \text{ with initial condition } \ y=1 \ \text{ for } \ x=0.$$

Here, $M=x+y+1$, $N=x-y^2+3$ and $\partial M/\partial y=\partial N/\partial x$ (=1), at all points of the $xy$-plane. The system (4) is written

$$\partial u/\partial x = x+y+1, \quad \partial u/\partial y = x-y^2+3.$$

The first equation yields

$$u = x^2/2+xy+x+\varphi(y),$$

while by the second one we get

$$\varphi'(y) = -y^2+3 \ \Rightarrow \ \varphi(y) = -y^3/3+3y+C_1 \, .$$

Thus,

$$u = x^2/2-y^3/3+xy+x+3y+C_1 \, .$$

The general solution (3) is $u(x,y)=C_2$ . Putting $C_2-C_1\equiv C$, we have:

$$x^2/2-y^3/3+xy+x+3y = C \ \text{ (general solution).}$$

Making the substitutions $x=0$ and $y=1$ (as required by the initial condition), we find $C=8/3$ and

$$x^2/2-y^3/3+xy+x+3y = 8/3 \ \text{ (particular solution).}$$

## 3.3 Integrating Factor

Assume that the ODE

$$M(x, y)\, dx + N(x, y)\, dy = 0 \qquad (1)$$

is not exact [i.e., the left-hand side is not a total differential of some function $u(x,y)$]. We say that this equation admits an *integrating factor* $\mu(x,y)$ if there exists a function $\mu(x,y)$ such that the ODE $\mu(Mdx+Ndy)=0$ *is* exact; that is, such that the expression $\mu(Mdx+Ndy)$ is a total differential of a function $u(x,y)$:

$$\mu(x, y)\, [M(x, y)\, dx + N(x, y)\, dy] = du \qquad (2)$$

Then the initial ODE (1) reduces to the differential relation $du=0 \Rightarrow$

$$u(x, y) = C \qquad (3)$$

*on the condition that* the function $\mu(x,y)$ does not vanish identically when $x$ and $y$ are related by (3). Relation (3) is a first integral of the ODE (1) and expresses the general solution of this equation.

***Example:*** The ODE $ydx–xdy=0$ is not exact, since $M=y$, $N= –x$, so that $\partial M/\partial y=1$, $\partial N/\partial x= –1$. However, the equation

$$\frac{1}{y^2}(ydx - xdy) = 0$$

*is* exact, given that the left-hand side is equal to $d(x/y)$. Thus,

$$d(x/y) = 0 \;\Rightarrow\; y = Cx\,.$$

The solution is acceptable since the integrating factor $\mu=1/y^2$ does not vanish identically for $y= Cx$.

## 3.4 Higher-Order Differential Equations

In the case of an ODE of second order or higher, a first integral leads to a *reduction of order* of the ODE.

Consider an ODE of order $n$:

$$F[x, y, y', y'', \ldots, y^{(n)}] = 0 \qquad (1)$$

(where $y^{(n)} \equiv d^n y/dx^n$). We assume that the left-hand side of (1) can be written as the derivative of an expression $\Phi$ of order ($n$–1):

$$F[x, y, y', y'', \cdots, y^{(n)}] = \frac{d}{dx}\Phi[x, y, y', \cdots, y^{(n-1)}] \qquad (2)$$

Then (1) reduces to $d\Phi/dx=0 \Rightarrow$

$$\Phi\,[\,x, y, y', \ldots, y^{(n-1)}\,] = C \tag{3}$$

Relation (3) is a *first integral* of (1); it is an ODE of order ($n$–1).

***Example:*** Consider the second-order ODE

$$yy''+(y')^2 = 0.$$

We notice that the left-hand side is equal to $d\,(yy')/dx$ . Hence, the given equation is written $yy' = C$ (which is a first-order ODE), from which we have

$$y^2 = C_1\,x+C_2\,.$$

Sometimes the left-hand side of (1) is not in itself a total derivative but can transform into one upon multiplication by a suitable *integrating factor*

$$\mu\,[x, y, y', \ldots, y^{(n-1)}]\;.$$

That is,

$$\mu[x, y, y',\cdots, y^{(n-1)}]\; F\,[x, y, y', y'',\cdots, y^{(n)}] = \frac{d}{dx}\Phi\,[x, y, y',\cdots, y^{(n-1)}] \tag{4}$$

Then $d\Phi/dx = 0$, so we are led again to a first integral of the form (3).

***Example:*** Consider the ODE

$$yy'' - (y')^2 = 0\,.$$

Multiplied by $\mu = 1/y^2$ , the left-hand side becomes $(y'/y)'$. The given equation is then written

$$(y'/y)' = 0 \;\Rightarrow\; y'/y = C \;\text{ (a first-order ODE),}$$

by which we get

$$y = C_1\,e^{Cx}.$$

## 3.5 Application: Newton's Second Law in One Dimension

In this section, one or more dots will be used to denote derivatives of various orders with respect to the variable $t$: $\dot{x} = dx/dt$ , $\ddot{x} = d^2x/dt^2$ , etc.

We consider the second-order ODE

$$m\ddot{x} = F(x) \tag{1}$$

with initial conditions $x(t_0)=x_0$ and $v(t_0)=v_0$ , where $v=dx/dt$. Physically, relation (1) expresses *Newton's Second Law* for a particle of mass $m$ moving with instantaneous velocity $v(t)$ along the $x$ axis, under the action of a force (or, more correctly, a force field) $F(x)$. By solving (1) we find the position $x=x(t)$ of the particle as a function of time $t$.

We define an auxiliary function $U(x)$ (*potential energy* of the particle) by

$$U(x) = -\int^x F(x')\,dx' \quad \Leftrightarrow \quad F(x) = -\frac{d}{dx}U(x) \tag{2}$$

(with arbitrary lower limit of integration). The function $U$ may always be defined in a one-dimensional problem, which is not the case in higher dimensions since the integral corresponding to the one in (2) will generally depend on the path of integration and thus will not be uniquely defined (except in the case of *conservative fields*; see Sec. 1.5). We also note that the function $U$ depends on time *t through x only*; i.e., $U$ is not an explicit function of $t$ ($\partial U/\partial t=0$). This means that the value of $U$ changes exclusively due to the motion of the particle along the $x$-axis, while at any fixed point $x$ the value of $U$ is constant in time (see discussion in Appendix A).

Equation (1) is now written

$$m\ddot{x} + \frac{dU}{dx} = 0 \ .$$

The left-hand side is not a perfect derivative *with respect to t*. Let us try the integrating factor $\mu = \dot{x}$:

$$\dot{x}\left(m\ddot{x} + \frac{dU}{dx}\right) = 0 \ \Rightarrow \ m\dot{x}\ddot{x} + \dot{x}\frac{dU}{dx} = 0 \ .$$

But,

$$\dot{x}\ddot{x} = \frac{1}{2}\frac{d}{dt}(\dot{x}^2) \qquad \text{and} \qquad \dot{x}\frac{dU}{dx} = \frac{dU}{dx}\frac{dx}{dt} = \frac{dU}{dt} \ .$$

We thus have

$$\frac{d}{dt}\left[\frac{1}{2}m\dot{x}^2 + U(x)\right] = 0 \ \Rightarrow$$

$$E \equiv \frac{1}{2}m\dot{x}^2 + U(x) \ = \frac{1}{2}mv^2 + U(x) \ = \text{constant} \tag{3}$$

Equation (3) expresses the *conservation of mechanical energy* of the particle. Notice that this result is an immediate consequence of Newton's law. (In higher dimensions this principle is valid only in the case where the force field $\vec{F}(\vec{r})$ is *conservative*; see Sec. 1.5 and Appendix A.)

Relation (3), which constitutes a first integral of the ODE (1), is a first-order ODE that is easy to integrate. We have:

$$\left(\frac{dx}{dt}\right)^2 = \frac{2}{m}[E - U(x)] \ .$$

Taking the case where $v=dx/dt>0$, we write

$$\frac{dx}{dt}=\left\{\frac{2}{m}[E-U(x)]\right\}^{1/2} \Rightarrow \int_{x_0}^{x}\frac{dx}{\{\cdots\}^{1/2}}=\int_{t_0}^{t}dt$$

(where we have taken into account that $x=x_0$ for $t=t_0$). Finally,

$$\int_{x_0}^{x}\frac{dx}{\left\{\frac{2}{m}[E-U(x)]\right\}^{1/2}}=t-t_0 \tag{4}$$

As derived, relation (4) is valid for $v>0$. For $v<0$ one must put $-dx$ in place of $dx$ inside the integral. In general, in cases where the velocity $v$ is positive in a part of the motion and negative in another part, it may be necessary to perform the integration *separately* for each part of the motion.

Relation (4) represents a particular solution of (1) for the given initial conditions. By putting $x=x_0$ and $v=v_0$ in (3) and by taking into account that $E$ is constant, we can determine the value of the parameter $E$ that appears in the solution (4):

$$E=\frac{1}{2}m{v_0}^2+U(x_0) \tag{5}$$

***Comment:*** It is evident from (2) that the function $U(x)$ is arbitrary to within an additive constant whose value will depend on the choice of the lower limit in the integral defining $U$. Through (3), the same arbitrariness is passed on to the value of the constant $E$; it disappears, however, upon taking the difference $E-U(x)$. Thus, this arbitrariness does not affect the result of the integration in (4).

***Example:*** Rectilinear motion under the action of a *constant* force $F$. We take $t_0=0$, $x_0=x(0)=0$, $v(0)=v_0$, and we assume that $v=dx/dt>0$ (in particular, $v_0>0$) for the part of the motion that interests us. From (2) we have (making the arbitrary assumption that $U=0$ for $x=0$):

$$\frac{dU}{dx}=-F \Rightarrow \int_0^U dU=-F\int_0^x dx \Rightarrow U(x)=-Fx \ .$$

Equation (4) then yields

$$\int_0^x\frac{dx}{(E+Fx)^{1/2}}=\left(\frac{2}{m}\right)^{1/2}t \Rightarrow (E+Fx)^{1/2}=\frac{F}{2}\left(\frac{2}{m}\right)^{1/2}t+E^{1/2} \ .$$

Squaring this, we find

$$x=\frac{F}{2m}t^2+\left(\frac{2E}{m}\right)^{1/2}t \ .$$

We set $F/m=a=const.$ (acceleration of the particle). Also, from (5) we have that $E=mv_0^2/2$ [since $U(0)=0$]. Thus, finally (taking into account that $v_0>0$),

$$x=\frac{1}{2}a\,t^2+v_0\,t\ ,$$

which is the familiar formula for uniformly accelerated rectilinear motion.

***Problem:*** Show that a result of the same form will ensue in the case where $v<0$. [*Hint:* Use (4) with $-dx$ in place of $dx$; put $v(0)=v_0$ with $v_0<0$.]

***Note:*** For higher-dimensional conservative systems, conservation of mechanical energy alone is not sufficient in order to obtain a complete solution of the problem; additional conservation laws are needed. For example, motion under a central force is essentially a two-dimensional problem, given that this motion takes place on a single plane. Reduction of order of Newton's equation of motion thus requires *two* first integrals, corresponding to conservation of mechanical energy and angular momentum [1-3].

## CHAPTER 4

# SYSTEMS OF ORDINARY DIFFERENTIAL EQUATIONS

### 4.1 Solution by Seeking First Integrals

We consider a system of $n$ ordinary differential equations (ODEs) of the first order, for $n$ unknown functions $x_1(t), x_2(t), \ldots , x_n(t)$:

$$\frac{dx_i}{dt} = f_i(x_1, x_2, \cdots, x_n, t) \qquad (i = 1, 2, \cdots, n) \tag{1}$$

If the functions $f_i$ are not *explicitly* dependent on $t$ (i.e., if $\partial f_i / \partial t = 0$ for $i=1,2,...,n$) the system (1) is called *autonomous*:

$$\frac{dx_i}{dt} = f_i(x_1, x_2, \cdots, x_n) \qquad (i = 1, 2, \cdots, n) \tag{2}$$

A *conservation law* for the system (1) is an ODE of the form

$$\frac{d}{dt} \Phi(x_1, x_2, \cdots, x_n, t) = 0 \tag{3}$$

which is valid *as a consequence of the system* (i.e., is not satisfied identically). Equation (3) is immediately integrable:

$$\Phi(x_1, x_2, \cdots, x_n, t) = C \tag{4}$$

The function $\Phi$ is a *first integral* of the system (1). It retains a constant value *when the* $x_1, x_2, \ldots , x_n$ *satisfy the system* (that is, it is not identically constant but reduces to a constant on solutions of the system).

If one or more first integrals of the system are known, one can trivially produce an infinity of first integrals by taking sums, multiples, products, powers, etc., of them. We are only interested, however, in first integrals that are *independent* of one another, since it is in this case that we obtain the most useful information for solution of the problem.

Let us now assume that we manage to find $k$ independent first integrals of the system (1) (where $k \leq n$):

$$\begin{aligned} &\Phi_1(x_1, x_2, \cdots, x_n, t) = C_1 \\ &\Phi_2(x_1, x_2, \cdots, x_n, t) = C_2 \\ &\vdots \\ &\Phi_k(x_1, x_2, \cdots, x_n, t) = C_k \end{aligned} \tag{5}$$

Relations (5) allow us to express $k$ of the variables $x_1, \ldots, x_n$ in terms of the remaining $(n–k)$ variables and $t$. We thus eliminate $k$ unknown functions from the problem, so that the system (1) reduces to one with *fewer unknowns*, that is, $(n–k)$. If $k=n$, then *all* unknown functions $x_1, \ldots, x_n$ can be determined algebraically from system (5) without the necessity of integrating the differential system (1) itself.

The autonomous system (2) is written

$$\frac{dx_i}{f_i(x_1, \cdots, x_n)} = dt \qquad (i = 1, 2, \cdots, n) \tag{6}$$

Since the $f_i$ do not contain $t$ directly, this variable can be eliminated from the system. Indeed, since all left-hand sides in (6) are equal to $dt$, they will be equal to one another. Hence,

$$\frac{dx_1}{f_1(x_1, \cdots, x_n)} = \frac{dx_2}{f_2(x_1, \cdots, x_n)} = \cdots = \frac{dx_n}{f_n(x_1, \cdots, x_n)} \tag{7}$$

Relation (7) represents a system of $(n–1)$ equations in $n$ variables $x_1, x_2, \ldots, x_n$. To solve it we seek $(n–1)$ independent first integrals of the form

$$\Phi_j(x_1, x_2, \cdots, x_n) = C_j \qquad (j = 1, 2, \cdots, n-1) \tag{8}$$

We also seek a first integral $\Phi_n$ of the complete system (6):

$$\Phi_n(x_1, x_2, \cdots, x_n, t) = C_n \tag{9}$$

Relations (8) and (9) constitute a system of $n$ algebraic equations in $(n+1)$ variables. By solving this system for the $x_1, x_2, \ldots, x_n$, we can express these variables as functions of $t$.

Analytically, one way of integrating the system (6) is the following: With the aid of relations (8), we express $(n–1)$ of the $n$ variables $x_1, \ldots, x_n$ as functions of the remaining variable. Assume, for example, that the $x_1, x_2, \ldots, x_{n-1}$ are expressed as functions of $x_n$. Taking (6) with $i=n$, we have:

$$\frac{dx_n}{f_n(x_1, \cdots, x_n)} = dt \;\Rightarrow\; \int \frac{dx_n}{f_n(x_1, \cdots, x_n)} \equiv F(x_n) + c = t + c' \;\Rightarrow$$

$$\Phi_n(x_n, t) \equiv F(x_n) - t = C_n \tag{10}$$

Equation (10) allows us to express the variable $x_n$ as a function of $t$. Given that the $x_1, x_2, \ldots, x_{n-1}$ are known functions of $x_n$, the above $(n–1)$ variables can, in turn, also be expressed as functions of $t$.

***Examples:***

***1.*** Consider the system

$$\frac{dx}{dt} = y \quad (a) \qquad \frac{dy}{dt} = x \quad (b)$$

(Here, $x_1 \equiv x$, $x_2 \equiv y$.) We seek first integrals of this system. Two are sufficient for a complete solution of the problem. Adding ($a$) and ($b$), we have

$$d(x+y)/dt = x+y.$$

Putting $x+y=u$, we write $du/dt=u$, which yields

$$u = C_1 e^t \;\Rightarrow\; (x+y)\, e^{-t} = C_1 .$$

Similarly, subtracting ($b$) from ($a$) and putting $x-y=u$, we find $du/dt = -u \;\Rightarrow$

$$u = C_2 e^{-t} \;\Rightarrow\; (x-y)\, e^t = C_2 .$$

We have thus found two independent first integrals of the system:

$$\Phi_1(x, y, t) \equiv (x+y)\, e^{-t} = C_1 \;\;, \quad \Phi_2(x, y, t) \equiv (x-y)\, e^t = C_2 \;.$$

[*Exercise:* Verify that $d\Phi_1/dt=0$ and $d\Phi_2/dt=0$ *when* $x$ and $y$ are solutions of the system of ($a$) and ($b$). Note that $\Phi_1$ and $\Phi_2$ are not *identically* constant!]

By using the first integrals $\Phi_1$ and $\Phi_2$ we can now express $x$ and $y$ as functions of $t$. Putting $C_1$ and $C_2$ in place of $C_1/2$ and $C_2/2$, respectively, we find

$$x = C_1 e^t + C_2 e^{-t} \;\;, \quad y = C_1 e^t - C_2 e^{-t} \;.$$

***Comment:*** One can easily find more first integrals of the system ($a$), ($b$). For example, by eliminating $dt$ we have

$$dx/y = dy/x \;\Rightarrow\; xdx = ydy \;\Rightarrow\; d(x^2-y^2)=0,$$

so that $\Phi_3(x,y) \equiv x^2-y^2=C_3$. Let us note, however, that $\Phi_3=\Phi_1\Phi_2$. Thus, the relation $\Phi_3=const.$ is a trivial consequence of $\Phi_1=const.$ and $\Phi_2=const.$ In other words, $\Phi_3$ is not an *independent*, new first integral of the system; therefore, it does not furnish any useful new information for solution of the problem.

***2.*** Consider the system

$$\frac{dx}{dt} = y \quad (a) \qquad \frac{dy}{dt} = -x \quad (b)$$

We seek two first integrals. In this case, we get no useful information by adding or subtracting the two equations of the system. However, since this system is autonomous, we can eliminate $dt$ :

$$dx/y = -\,dy/x \;\Rightarrow\; xdx+ydy = 0 \;\Rightarrow\; d(x^2+y^2) = 0 \;\Rightarrow$$

$$\Phi_1(x, y) \equiv x^2+y^2 = C_1^2 \;.$$

To solve the problem completely, we need another first integral of the system; this time, one that contains $t$ explicitly. From ($a$) and ($b$) we have

$$x\,(dy/dt) - y\,(dx/dt) = -\,(x^2+y^2) \;\Rightarrow\; d(y/x)/dt = -\,[1 + (y/x)^2\,] \;.$$

Putting $y/x=u$, we write

$$du/(1+u^2) = -\,dt \;\Rightarrow\; d(t + \arctan u) = 0 \,,$$

from which we find

$$\Phi_2(x, y, t) \equiv t + \arctan(y/x) = C_2 \;.$$

We now use the first integrals $\Phi_1$ and $\Phi_2$ to solve the system (*a*), (*b*) algebraically. The relation $\Phi_2 = C_2$ yields

$$y = -x \tan(t - C_2) \;.$$

Then, by the relation $\Phi_1 = C_1^2$ we get

$$x^2\,[1+\tan^2(t-C_2)] = C_1^2 \,.$$

By using the identity $\cos^2 a = 1/1+\tan^2 a$ , it is not hard to show that

$$x = C_1 \cos(t - C_2) \quad \text{so that} \quad y = -\,C_1 \sin(t - C_2) \;.$$

***Comment:*** An alternative way to solve the problem is by *transformation of coordinates* from *Cartesian* $(x,y)$ to *polar* $(r,\theta)$, where $r \geq 0$ and $0 \leq \theta < 2\pi$. The transformation equations are

$$x = r\cos\theta \,, \quad y = r\sin\theta \quad \Leftrightarrow \quad r = (x^2+y^2)^{1/2} \,, \quad \theta = \arctan(y/x) \;.$$

The system (*a*), (*b*) is written

$$(dr/dt)\cos\theta - r\,(d\theta/dt)\sin\theta = r\sin\theta \,,$$

$$(dr/dt)\sin\theta + r\,(d\theta/dt)\cos\theta = -\,r\cos\theta \;.$$

Solving for the derivatives, we can separate the variables $r$ and $t$, finding a separate equation for each variable:

$$dr/dt = 0 \,, \quad d\theta/dt = -1 \,,$$

with corresponding solutions

$$r = C_1 \,, \quad \theta = -\,t + C_2 \,.$$

Substituting into the transformation equations, we have:

$$x = C_1 \cos(t - C_2) \;, \quad y = -\,C_1 \sin(t - C_2) \;,$$

as before. The first integrals of the system are easily found by solving the above equations for the constants $C_1^2$ and $C_2$ .

***3.*** Consider the system

$$\frac{dx}{dt} = y - z \quad (a) \qquad \frac{dy}{dt} = z - x \quad (b) \qquad \frac{dz}{dt} = x - y \quad (c)$$

Rather than trying to solve it analytically, we will express the solution implicitly with the aid of three independent first integrals (as required for an algebraic solution of the problem), at least one of which will contain the variable $t$ explicitly. Taking the sum $(a)+(b)+(c)$ , we have

$$d(x+y+z)/dt = 0 \quad \Rightarrow \quad \Phi_1(x, y, z) \equiv x+y+z = C_1 \;.$$

On the other hand, the combination $x.(a) + y.(b) + z.(c)$ yields

$$d(x^2+y^2+z^2)/dt = 0 \quad \Rightarrow \quad \Phi_2(x, y, z) \equiv x^2+y^2+z^2 = C_2 \;.$$

Now, by using the equations $\Phi_1=C_1$ and $\Phi_2=C_2$ we can express two of the dependent variables, say $x$ and $y$, in terms of the third variable, $z$. Then, relation ($c$) in the system is written in the form of an equation for a single variable $z$:

$$\frac{dz}{x(z)-y(z)}=dt \;\Rightarrow\; \int\frac{dz}{x(z)-y(z)}\equiv F(z)+c=t+c' \;\Rightarrow$$

$$\Phi_3(z,t)\equiv F(z)-t=C_3 \;.$$

The equation $\Phi_3=C_3$ allows us to express $z$ as a function of $t$. Given that $x$ and $y$ already are functions of $z$ (thus, implicitly, of $t$), the problem has been solved in principle.

***4.*** The system

$$\frac{dx}{x^2-y^2-z^2}=\frac{dy}{2xy}=\frac{dz}{2xz} \qquad (a)$$

contains two equations with three variables $x$, $y$, $z$. The solution of the system will allow us to express two of these variables as functions of the third (that is, the third variable plays here the same role as $t$ in the preceding examples). We seek two independent first integrals of the system. By the second equality in relation ($a$) we get

$$dy/y=dz/z \;\Rightarrow\; d(\ln y-\ln z)\equiv d\ln(y/z)=0 \;\Rightarrow\; \ln(y/z)=c \;\Rightarrow\; y/z=e^c\equiv C_1 \;.$$

Thus,

$$\Phi_1(y,z)\equiv y/z=C_1 \;.$$

We need to find one more first integral of system ($a$), this time containing $x$ explicitly. To this end, we apply a familiar property of proportions:

$$\frac{xdx}{x(x^2-y^2-z^2)}=\frac{ydy}{y(2xy)}=\frac{zdz}{z(2xz)}=\frac{xdx+ydy+zdz}{x(x^2+y^2+z^2)} \qquad (b)$$

Equating the last term with the second, we have:

$$\frac{xdx+ydy+zdz}{x^2+y^2+z^2}=\frac{dy}{2y} \;\Rightarrow\; \frac{d(x^2+y^2+z^2)}{x^2+y^2+z^2}=\frac{dy}{y} \;\Rightarrow\; d\ln\left(\frac{x^2+y^2+z^2}{y}\right)=0 \;\Rightarrow$$

$$\ln\left(\frac{x^2+y^2+z^2}{y}\right)=c \;\Rightarrow\; \frac{x^2+y^2+z^2}{y}=e^c\equiv C_2 \;.$$

Thus,

$$\Phi_2(x,y,z)\equiv(x^2+y^2+z^2)/y=C_2 \;.$$

The relations $\Phi_1=C_1$ and $\Phi_2=C_2$ represent the solution of system ($a$) since, by means of them, we can express two of the variables as functions of the third.

***Comment:*** If we had chosen to equate the last term in ($b$) with the third term, rather than with the second, we would have found, in a similar way,

$$\Phi_3(x,y,z)\equiv(x^2+y^2+z^2)/z=C_3 \;.$$

This, however, is not a new, *independent* first integral since, as is easy to show, $\Phi_3=\Phi_1\Phi_2$. Thus, the constancy of $\Phi_1$ and $\Phi_2$ automatically guarantees the constancy of $\Phi_3$ as well, so that the relation $\Phi_3 = const.$ does not provide any new information for solution of the problem.

## 4.2 Application to First-Order Partial Differential Equations

We now examine the relation between systems of ordinary differential equations (ODEs) and first-order partial differential equations (PDEs). More on PDEs will be said in the next chapter (see also [1,2]). Here we will confine ourselves to PDEs whose solutions are functions of two variables $x$, $y$. We will denote by $z$ the variable representing the unknown function in the PDE. Thus, the solution of the equation will be of the form $z=f(x,y)$.

The general solution of a PDE of order $p$ is dependent on $p$ arbitrary functions. Let us see some examples:

1. $\dfrac{\partial z}{\partial x} = x + y \;\Rightarrow\; z = \dfrac{x^2}{2} + xy + \varphi(y)$ .

2. $\dfrac{\partial z}{\partial y} = xyz$ . We integrate, treating $x$ as a constant:

$$\int \frac{dz}{z} = x \int y\,dy \;\Rightarrow\; \ln z = \frac{xy^2}{2} + \ln \varphi(x) \;\Rightarrow\; z = \varphi(x)\, e^{xy^2/2} \ .$$

3. $\dfrac{\partial^2 z}{\partial x \partial y} = 0 \;\Rightarrow\; \dfrac{\partial}{\partial x}\left(\dfrac{\partial z}{\partial y}\right) = 0 \;\Rightarrow\; \dfrac{\partial z}{\partial y} = \varphi(y) \;\Rightarrow$

$$z = \int \varphi(y)\,dy + \varphi_1(x) = \varphi_1(x) + \varphi_2(y) \,.$$

A first-order PDE is called *quasilinear* if it is linear in the partial derivatives of $z$ (but not necessarily linear in $z$ itself). This PDE has the general form

$$P(x,y,z)\frac{\partial z}{\partial x} + Q(x,y,z)\frac{\partial z}{\partial y} = R(x,y,z) \qquad (1)$$

The solution $z(x,y)$ is found by the following process, stated here without proof (see [1,2]):

1. We form the *characteristic system* of ODEs,

$$\frac{dx}{P(x,y,z)} = \frac{dy}{Q(x,y,z)} = \frac{dz}{R(x,y,z)} \qquad (2)$$

Relation (2) represents a system of two differential equations with three variables $x$, $y$, $z$. By solving it, two of these variables may be written as functions of the third. The solution can also be expressed as an algebraic system of two independent first integrals of the form

$$\Psi_1(x, y, z) = C_1 \quad , \quad \Psi_2(x, y, z) = C_2 \tag{3}$$

2. We consider an *arbitrary* function $\Phi$ of $C_1$ , $C_2$ and we form the equation $\Phi(C_1,C_2)=0$, or, in view of (3),

$$\Phi[\Psi_1(x, y, z), \Psi_2(x, y, z)] = 0 \tag{4}$$

Equation (4) defines a relation of the form $z=f(x,y)$, depending on an arbitrary function. This relation constitutes the solution of the PDE (1).

***Note:*** By making the special choices

$$\Phi(\Psi_1,\Psi_2) = \Psi_1(x, y, z) - C_1 \quad \text{and} \quad \Phi(\Psi_1,\Psi_2) = \Psi_2(x, y, z) - C_2 \ ,$$

and by demanding that $\Phi(\Psi_1,\Psi_2)=0$ in each case, we are led to relations (3). That is, *the first integrals of the characteristic system* (2) *are particular solutions of the PDE* (1).

***Special case:*** If $R(x,y,z)=0$, and if the functions $P$ and $Q$ do not contain $z$, then the PDE (1) is called *homogeneous linear*:

$$P(x, y)\frac{\partial z}{\partial x} + Q(x, y)\frac{\partial z}{\partial y} = 0 \tag{5}$$

The characteristic system (2) is written

$$\frac{dx}{P(x, y)} = \frac{dy}{Q(x, y)} = \frac{dz}{0} \tag{6}$$

In order for the $dx/P$ and $dy/Q$ to be finite, it is necessary that $dz=0 \Leftrightarrow z=C_1$ . We thus have a first integral

$$\Psi_1(z) \equiv z = C_1 \tag{7}$$

Next, we solve the ODE

$$dx/P(x,y) = dy/Q(x,y)$$

and express the solution in the form of a first integral:

$$\Psi_2(x, y) = C_2 \tag{8}$$

Finally, we take an *arbitrary* function $\Phi$ of $C_1$ , $C_2$ and we demand that $\Phi(C_1,C_2)=0$. Making use of (7) and (8), and putting $\Psi_2(x,y) \equiv \Psi(x,y)$, we have:

$$\Phi[z, \Psi(x, y)] = 0 \tag{9}$$

Relation (9) allows us to express $z$ as a function of $x$ and $y$, thus obtaining the solution of the PDE (5). The arbitrariness in the choice of $\Phi$ means that this solution will depend on an arbitrary function.

***Examples:***

***1.*** Consider the PDE

$$\frac{\partial z}{\partial x}+\frac{\partial z}{\partial y}=1 \qquad (a)$$

Here, $P=Q=R=1$. We form the characteristic system (2):

$$dx/P = dy/Q = dz/R \;\Rightarrow\; dx = dy = dz \;.$$

We find two first integrals:

$$dx = dy \;\Rightarrow\; d(x-y) = 0 \;\Rightarrow\; \Psi_1(x, y) \equiv x-y = C_1 \,,$$

$$dx = dz \;\Rightarrow\; d(z-x) = 0 \;\Rightarrow\; \Psi_2(x, z) \equiv z-x = C_2 \,.$$

The general solution of the PDE ($a$) is

$$\Phi(C_1, C_2) = 0 \;\Rightarrow\; \Phi(x-y, z-x) = 0 \;\Rightarrow\; z-x = F(x-y) \;\Rightarrow\; z = x + F(x-y) \,,$$

where the function $\Phi$ is chosen arbitrarily while $F$ is dependent upon the choice of $\Phi$. Alternatively, we could have taken

$$dx = dy \;\Rightarrow\; \Psi_1(x,y) \equiv x-y = C_1 \,,$$

$$dy = dz \;\Rightarrow\; \Psi_3(y,z) \equiv z-y = C_3 \,,$$

with corresponding general solution $z= y+G(x-y)$. The two solutions we found, however, are not independent of each other. Indeed, by putting $G(x-y)=x-y+F(x-y)$, the second solution reduces to the first.

***Exercise:*** Verify that the expression $z=x+F(x-y)$ indeed satisfies the PDE ($a$). [*Hint:* Put $x-y=u$ and notice that $\partial F/\partial x=F'(u)(\partial u/\partial x)=F'(u)$, $\partial F/\partial y= - F'(u)$ .]

***2.*** Consider the PDE

$$x\frac{\partial z}{\partial y}-y\frac{\partial z}{\partial x}=0 \qquad (a)$$

Here, $P= -y$, $Q=x$, $R=0$ (homogeneous linear). The characteristic system (6) is written

$$dx/(-y) = dy/x = dz/0 \,.$$

We have:

$$dz = 0 \;\Rightarrow\; \Psi_1(z) \equiv z = C_1 \;,$$

$$-dx/y = dy/x \;\Rightarrow\; xdx+ydy = 0 \;\Rightarrow\; d(x^2+y^2) = 0 \;\Rightarrow\; \Psi_2(x,y) \equiv x^2+y^2 = C_2 \;.$$

The general solution of ($a$) is (with arbitrary $\Phi$)

$$\Phi(C_1, C_2) = 0 \;\Rightarrow\; \Phi(z, x^2+y^2) = 0 \;\Rightarrow\; z = F(x^2+y^2) \quad \text{(arbitrary } F) \,.$$

***Exercise:*** Verify that the above expression satisfies the PDE ($a$). [*Hint:* Put $x^2+y^2=u$ and notice that $\partial F/\partial x = F'(u)(\partial u/\partial x) = 2xF'(u)$, $\partial F/\partial y = 2yF'(u)$.]

***3.*** Consider the PDE

$$x\frac{\partial z}{\partial x} + y\frac{\partial z}{\partial y} = z \qquad (a)$$

The characteristic system reads

$$dx/x = dy/y = dz/z\,.$$

We have:

$$dx/x = dy/y \;\Rightarrow\; d(\ln x - \ln y) = 0 \;\Rightarrow\; \ln(x/y) = c \;\Rightarrow\; \Psi_1(x, y) \equiv x/y = C_1\,,$$

$$dx/x = dz/z \;\Rightarrow\; d(\ln z - \ln x) = 0 \;\Rightarrow\; \ln(z/x) = c' \;\Rightarrow\; \Psi_2(x, z) \equiv z/x = C_2\,.$$

The general solution of ($a$) is

$$\Phi(C_1, C_2) = 0 \;\Rightarrow\; \Phi(x/y, z/x) = 0 \;\Rightarrow\; z/x = F(x/y) \;\Rightarrow\; z = x\,F(x/y)\,,$$

where the function $\Phi$ is arbitrary while $F$ depends on the choice of $\Phi$. Alternatively, we could have taken

$$dx/x = dy/y \;\Rightarrow\; \Psi_1(x,y) \equiv x/y = C_1\,,$$

$$dy/y = dz/z \;\Rightarrow\; \Psi_3(y,z) \equiv z/y = C_3\,,$$

with corresponding general solution $z = y\,G(x/y)$. However, the two solutions we found are not independent of each other. Indeed, by putting $G(x/y)=(x/y)F(x/y)$, the second solution reduces to the first.

## 4.3 System of Linear Equations

We will now study a method of integration that does not employ the tool of first integrals. We consider a *homogeneous linear* system of ODEs, with constant coefficients:

$$\frac{dx_i}{dt} = \sum_{j=1}^{n} a_{ij}\,x_j \qquad (i = 1, 2, \cdots, n) \tag{1}$$

(where the $a_{ij}$ are constants). Analytically,

$$\begin{aligned}
\frac{dx_1}{dt} &= a_{11}x_1 + a_{12}x_2 + \cdots + a_{1n}x_n \\
\frac{dx_2}{dt} &= a_{21}x_1 + a_{22}x_2 + \cdots + a_{2n}x_n \\
&\;\;\vdots \\
\frac{dx_n}{dt} &= a_{n1}x_1 + a_{n2}x_2 + \cdots + a_{nn}x_n
\end{aligned} \tag{2}$$

We seek a solution of the form

$$x_1 = \psi_1 e^{kt} \ , \quad x_2 = \psi_2 e^{kt} \ , \ \cdots \ , \quad x_n = \psi_n e^{kt} \qquad (3)$$

where the $\psi_1,...,\psi_n$ are constants. Substituting (3) into (2) and eliminating the common factor $e^{kt}$, we find a system of $n$ algebraic equations for the $\psi_1,...,\psi_n$:

$$\begin{aligned} &(a_{11}-k)\psi_1 + a_{12}\psi_2 + \cdots + a_{1n}\psi_n = 0 \\ &a_{21}\psi_1 + (a_{22}-k)\psi_2 + \cdots + a_{2n}\psi_n = 0 \\ &\qquad\vdots \\ &a_{n1}\psi_1 + a_{n2}\psi_2 + \cdots + (a_{nn}-k)\psi_n = 0 \end{aligned} \qquad (4)$$

In order that the homogeneous linear system (4) may have a non-trivial solution $(\psi_1,...,\psi_n)$ [different, that is, from the null solution (0,...,0)] the determinant of the coefficients must vanish:

$$\begin{vmatrix} a_{11}-k & a_{12} & \cdots & a_{1n} \\ a_{21} & a_{22}-k & \cdots & a_{2n} \\ \vdots & \vdots & & \vdots \\ a_{n1} & a_{n2} & \cdots & a_{nn}-k \end{vmatrix} = 0 \qquad (5)$$

(*characteristic equation* of the system). Relation (5) is an $n$th-degree polynomial equation for $k$. By solving it we find the values of the constant $k$ for which the system (4) has non-trivial solutions for the $\psi_i$ $(i=1,...,n)$.

It would be easier (and more elegant also!) to write our equations in matrix form. To this end we define the $(n\times n)$ matrix $A$ and the $(n\times 1)$ matrices $X$ and $\Psi$, as follows:

$$A = \begin{bmatrix} a_{11} & a_{12} & \cdots & a_{1n} \\ a_{21} & a_{22} & \cdots & a_{2n} \\ \vdots & \vdots & & \vdots \\ a_{n1} & a_{n2} & \cdots & a_{nn} \end{bmatrix} \ , \quad X = \begin{bmatrix} x_1 \\ x_2 \\ \vdots \\ x_n \end{bmatrix} \ , \quad \Psi = \begin{bmatrix} \psi_1 \\ \psi_2 \\ \vdots \\ \psi_n \end{bmatrix} \ .$$

The system (2) is written

$$\frac{dX}{dt} = AX \qquad (6)$$

where we have used the facts that

$$(dX/dt)_i = dx_i/dt \ \text{(see Appendix B) and } (AX)_i = \textstyle\sum_j a_{ij} x_j \, .$$

The candidate solution (3) takes on the form

$$X = \Psi e^{kt} \qquad (7)$$

Substituting (7) into (6), we get a matrix equation corresponding to the system (4):

$$A\Psi = k\Psi \quad \Leftrightarrow \quad (A - k\cdot 1_n)\Psi = 0 \tag{8}$$

where $1_n$ denotes the ($n\times n$) unit matrix. Relation (8) has the form of an *eigenvalue equation*. In order for this to have a non-trivial solution $\Psi$, it is necessary that

$$\det(A - k\cdot 1_n) = 0 \tag{9}$$

which is precisely Eq. (5). The values of $k$ that satisfy Eq. (9) are the *eigenvalues* of matrix $A$, while the corresponding non-trivial solutions $\Psi$ of (8) are the *eigenvectors* of $A$.

For each root $k_i$ of the characteristic equation (9), the solution of the eigenvalue equation (8) yields a non-trivial eigenvector $\Psi^{(i)}$. If *all* roots $k_i$ of (9) are *different* from one another, we obtain $n$ linearly independent eigenvectors $\Psi^{(i)}$ and an equal number of linearly independent solutions (7) of the ODE (6):

$$X^{(i)} = \Psi^{(i)} e^{k_i t} \quad (i = 1, 2, \cdots, n) \tag{10}$$

The general solution of the linear equation (6) is then

$$X = \sum_{i=1}^{n} c_i X^{(i)} = \sum_{i=1}^{n} c_i \Psi^{(i)} e^{k_i t} \tag{11}$$

where $c_1, ..., c_n$ are arbitrary constants.

The case of *multiple roots* of the characteristic equation is more complex. Let $k_i$ be a root of (9), of multiplicity $\lambda_i$. Then, the solution $X^{(i)}$ of the ODE (6) is not the one given by relation (10) but has the more general form [1]

$$X^{(i)} = \left(\Psi_0^{(i)} + \Psi_1^{(i)} t + \cdots + \Psi_{\lambda_i - 1}^{(i)} t^{\lambda_i - 1}\right) e^{k_i t} \tag{12}$$

Again, the general solution of (6) is $X = \sum_i c_i X^{(i)}$.

***Examples:***

***1.*** Consider the system

$$\frac{dx}{dt} = x + 2y \ , \quad \frac{dy}{dt} = 4x + 3y \qquad (a)$$

In matrix form,

$$\frac{dX}{dt} = AX \quad \text{where} \quad X = \begin{bmatrix} x \\ y \end{bmatrix}, \quad A = \begin{bmatrix} 1 & 2 \\ 4 & 3 \end{bmatrix}.$$

We seek the eigenvalues $k$ and the eigenvectors $\Psi$ of the matrix $A$, according to (8): $A\Psi=k\Psi$. Relation (5) for the eigenvalues is written

$$\begin{vmatrix} 1-k & 2 \\ 4 & 3-k \end{vmatrix} = 0 \;\Rightarrow\; k^2-4k-5=0 \;\Rightarrow\; k_1=5,\; k_2=-1\,.$$

Let

$$\Psi^{(1)} = \begin{bmatrix} \alpha \\ \beta \end{bmatrix}, \quad \Psi^{(2)} = \begin{bmatrix} \gamma \\ \delta \end{bmatrix}$$

be the eigenvectors corresponding to the eigenvalues $k_1$ , $k_2$ . The relation $A\Psi^{(1)}=k_1\Psi^{(1)}$ leads to a linear system of two equations for $\alpha$ and $\beta$. Since this system is homogeneous, these equations are not independent of each other but yield the same result, $\beta=2\alpha$. Hence,

$$\Psi^{(1)} = \begin{bmatrix} \alpha \\ 2\alpha \end{bmatrix} = \alpha \begin{bmatrix} 1 \\ 2 \end{bmatrix} \qquad \text{with } \textit{arbitrary } \alpha\,.$$

Similarly, the relation $A\Psi^{(2)}=k_2\Psi^{(2)}$ yields $\delta=-\gamma$, so that

$$\Psi^{(2)} = \begin{bmatrix} \gamma \\ -\gamma \end{bmatrix} = \gamma \begin{bmatrix} 1 \\ -1 \end{bmatrix} \qquad \text{with } \textit{arbitrary } \gamma\,.$$

The general solution (11) of the given system ($a$) is written

$$X = c_1 \Psi^{(1)} e^{k_1 t} + c_2 \Psi^{(2)} e^{k_2 t}\,.$$

Making substitutions and putting $c_1$ and $c_2$ in place of $c_1\alpha$ and $c_2\gamma$, respectively, we have:

$$\begin{bmatrix} x \\ y \end{bmatrix} = c_1 \begin{bmatrix} 1 \\ 2 \end{bmatrix} e^{5t} + c_2 \begin{bmatrix} 1 \\ -1 \end{bmatrix} e^{-t} \;\Rightarrow\; x = c_1 e^{5t} + c_2 e^{-t}\,,\quad y = 2c_1 e^{5t} - c_2 e^{-t}\,.$$

***2.*** Consider the system

$$\frac{dx}{dt} = x - y\,, \quad \frac{dy}{dt} = x + 3y \qquad (a)$$

In matrix form,

$$\frac{dX}{dt} = AX \quad \text{where} \quad X = \begin{bmatrix} x \\ y \end{bmatrix}, \quad A = \begin{bmatrix} 1 & -1 \\ 1 & 3 \end{bmatrix}.$$

The eigenvalues $k$ of the matrix $A$ are given by (5):

$$\begin{vmatrix} 1-k & -1 \\ 1 & 3-k \end{vmatrix} = 0 \;\Rightarrow\; k^2-4k+4=0 \;\Rightarrow\; k_1=k_2=2\,.$$

Here, the characteristic equation has a *double* root. Thus, we seek a solution $X$ of the form (12) with $\lambda_i=2$:

$$X = (\Psi_0 + \Psi_1 t)\, e^{2t} \qquad (b)$$

By the relation $dX/dt=AX$ we then get

$$(\Psi_1 + 2\Psi_0) + (2\Psi_1)\,t = A\Psi_0 + (A\Psi_1)\,t\,.$$

In order for this to be valid for all $t$, it is necessary that coefficients of equal powers of $t$ on the two sides of the equation be equal. That is,

$$A\Psi_1 = 2\Psi_1 \;\;,\quad A\Psi_0 = \Psi_1 + 2\Psi_0 \qquad (c)$$

Assume that

$$\Psi_0 = \begin{bmatrix} \alpha \\ \beta \end{bmatrix} \;,\quad \Psi_1 = \begin{bmatrix} \gamma \\ \delta \end{bmatrix}\,.$$

The first of relations ($c$), then, leads to a homogeneous linear system for $\gamma$ and $\delta$. The two equations of the system are not independent of each other but yield the same result, $\delta=-\gamma$. Hence,

$$\Psi_1 = \begin{bmatrix} \gamma \\ -\gamma \end{bmatrix} = \gamma \begin{bmatrix} 1 \\ -1 \end{bmatrix} \qquad \text{with } \textit{arbitrary } \gamma\,.$$

By the second equation in ($c$) we then get $\beta=-(\alpha+\gamma)$, so that

$$\Psi_0 = \begin{bmatrix} \alpha \\ -(\alpha+\gamma) \end{bmatrix} \qquad \text{with } \textit{arbitrary } \alpha\,.$$

The solution ($b$), now, of system ($a$) is written

$$\begin{bmatrix} x \\ y \end{bmatrix} = \begin{bmatrix} \alpha+\gamma t \\ -(\alpha+\gamma+\gamma t) \end{bmatrix} e^{2t} \;\Rightarrow\; \text{(by putting } \alpha=c_1\,,\ \gamma=c_2\text{)}$$

$$x = (c_1 + c_2\,t)\,e^{2t} \;\;,\quad y = -\,(c_1 + c_2 + c_2\,t)\,e^{2t}\,.$$

# CHAPTER 5

# DIFFERENTIAL SYSTEMS: GEOMETRIC VIEWPOINT

## 5.1 Dynamical Systems

We consider the system of first-order ODEs

$$\frac{dx_i}{dt} = f_i\,(x_1, x_2, \cdots, x_n,\, t) \qquad (i = 1, 2, \cdots, n) \tag{1}$$

The initial conditions are relations of the form $x_i\,(t_0) = x_{0i}$ .

Defining the vectors $X \equiv (x_1\,,\, x_2\,,\, \ldots\,,\, x_n\,)$ and $F \equiv (\,f_1\,, f_2\,,\, \ldots\,, f_n\,)$ , we rewrite the system (1) and its initial conditions in compact vector form:

$$\frac{dX}{dt} = F\,(X\,, t) \ , \quad X\,(t_0) = X_0 \tag{2}$$

where $X_0 \equiv (x_{01}\,,\, x_{02}\,,\, \ldots\,,\, x_{0n}\,)$ . The system (2) is called *autonomous* if $F$=$F(X)$, i.e., if the vector function $F$ is not explicitly dependent on $t$ :

$$\frac{dX}{dt} = F\,(X\,) \ , \quad X\,(t_0) = X_0 \tag{3}$$

Let $X(t) \equiv \big(x_1(t),\, x_2(t), \ldots,\, x_n(t)\big)$ be the general solution of system (2). This solution will depend on $n$ parameters that are determined by the initial conditions, thus are expressed in terms of the $x_{0i}$ . The solution defines an *integral curve* in the $(n$+1)-dimensional Euclidean space $R^n{\times}R$ with coordinates $(x_1, x_2, \ldots, x_n, t)$. The *projection* of this curve onto the space $R^n$: $(x_1, x_2, \ldots, x_n)$, i.e., the image of the mapping $(t{\in}R){\to}X(t){\in}R^n$ , defines a *trajectory* in $R^n$ (Fig. 5.1).

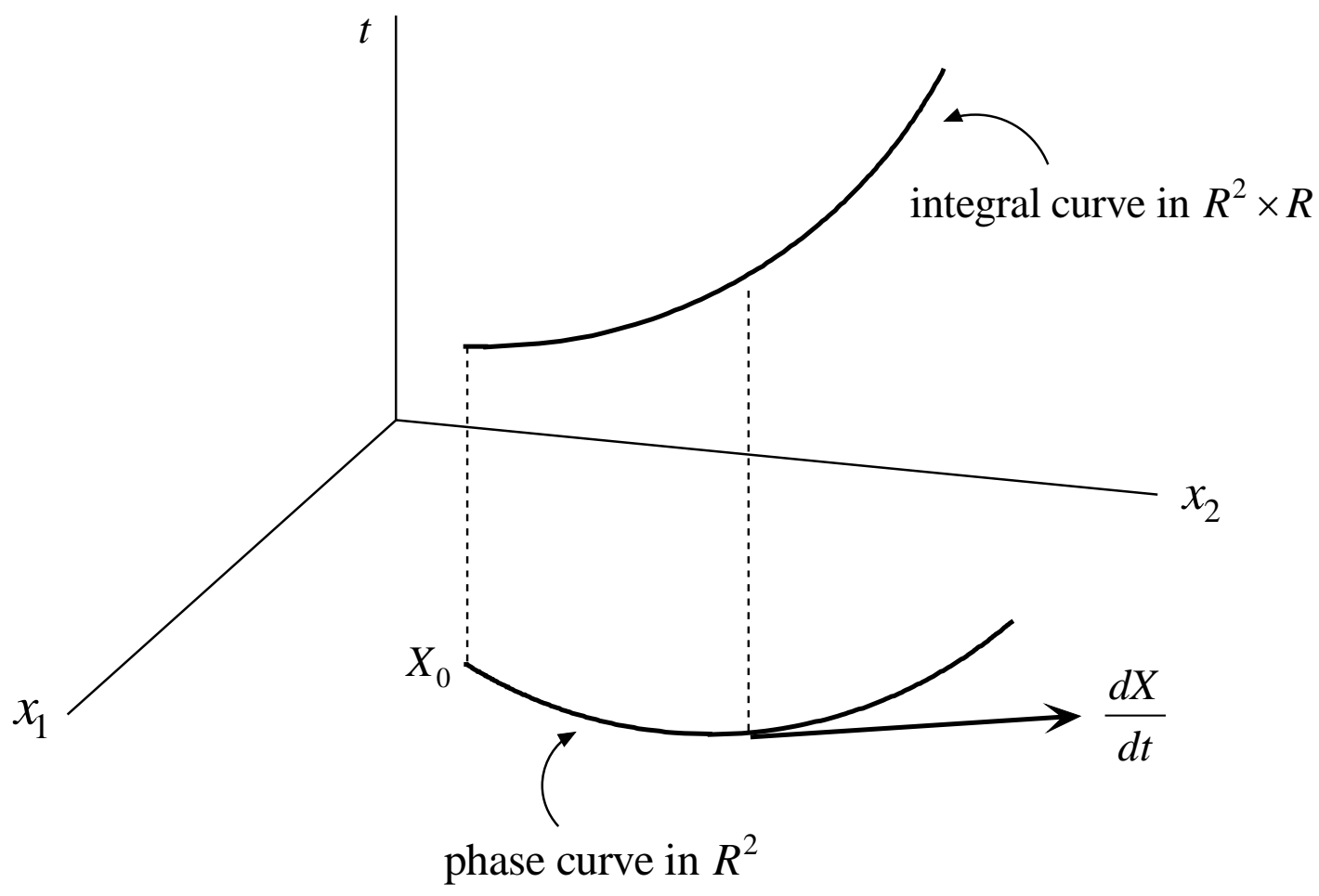

Fig. 5.1. Integral curve, phase curve, and phase velocity for $n$=2.

Under certain conditions on the vector function $F$ [1] the solution of system (2) for given initial conditions is *unique*. This means that *a unique integral curve passes through every point of* $R^n \times R$. Thus, *integral curves in* $R^n \times R$ *do not intersect*.

Adopting the terminology of Classical Mechanics, we call the system (2) a *dynamical system* and the space $R^n$: $(x_1,..., x_n)$ the *phase space* [1,2]. The variable $t$ represents *time* and the trajectory $X(t)$ in phase space describes a *phase curve* in $R^n$ (Fig. 5.1). The vector

$$\frac{dX}{dt} \equiv \left( \frac{dx_1}{dt}, \frac{dx_2}{dt}, \cdots, \frac{dx_n}{dt} \right) \tag{4}$$

is called the *phase velocity*; it represents the velocity of motion at the point $X(t)$ of the phase curve. The direction of the phase velocity is always *tangent to the phase curve* at $X(t)$. Finally, at every instant $t$ the vector function $F(X,t)$ defines a *vector field* in $R^n$, which, in view of the dynamical system (2), is a *velocity field*. In the case of an *autonomous* dynamical system of the form (3), the velocity field $F(X)$ is *static* (time-independent).

A physical analog may help visualize the situation better. Imagine that the entire phase space $R^n$ is filled with a “fluid” consisting of an infinite number of pointlike “particles”. At every moment $t_0$ a fluid particle passes through any point $X_0 \equiv (x_{01}, x_{02},..., x_{0n})$ of space, with velocity $(dX/dt)_0=F(X_0,t_0)$. For $t > t_0$ the particle describes a phase curve in $R^n$. Two particles passing by the same point $X_0$ at different moments $t_0$ and $t_1$ will generally move at this point with different velocities, except if $\partial F/\partial t=0 \Leftrightarrow F=F(X)$, i.e., except if the dynamical system is *autonomous* (thus the velocity field is *static*). In this latter case, *every* particle passing by a given point $X_0$ at any moment will move at this point with the *same* velocity and will describe the *same* phase curve in phase space. This means, in particular, that, in the case of a static velocity field a *unique* phase curve goes through every point of phase space $R^n$; that is, *the phase curves of an autonomous dynamical system do not intersect*. (It should be noted, however, that even in the non-autonomous case where the *phase* curves in $R^n$ can intersect, the *integral* curves in $R^n \times R$ still do *not* intersect as they represent unique solutions of the dynamical system.)

We will now concentrate our attention to the autonomous system (3) where the velocity field is static. In this case the (time-independent) phase curves constitute the *field lines* of $F(X)$. At each point of a field line the phase velocity $dX/dt$ is constant in time and is tangent to the line. Furthermore, a unique field line passes through every point of phase space; that is, *the field lines of* $F(X)$ *do not intersect*.

Analytically, the field lines of $F(X)$ are determined as follows: We consider the autonomous system

$$\frac{dx_i}{dt} = f_i(x_1, x_2, \cdots, x_n) \qquad (i=1,2,\cdots,n) \tag{5}$$

This is written $dx_i / f_i(x_k) = dt$ $(i=1, 2,..., n)$, where by $x_k$ we collectively denote the set of all variables $x_1,..., x_n$ . By eliminating $dt$ we obtain a system of $(n–1)$ equations with $n$ variables $x_1,..., x_n$ :

$$\frac{dx_1}{f_1(x_k)} = \frac{dx_2}{f_2(x_k)} = \cdots = \frac{dx_n}{f_n(x_k)} \tag{6}$$

By solving system (6) we can express $(n–1)$ of the variables as functions of the remaining variable. The solution may be expressed as a set of $(n–1)$ independent first integrals,

$$\Phi_j(x_1, x_2, \dots, x_n) = C_j \quad (j=1, 2,..., n-1) .$$

This solution determines the field lines of $F(X)$, which are curves in the phase space $R^n$: $(x_1,..., x_n)$. We note that these curves are *static* (time-independent) in the case of an autonomous system. As mentioned previously, the field lines of $F(X)$ do not intersect anywhere in $R^n$.

Finally, for a complete solution of the autonomous system (5) we also need a first integral of this system containing the variable $t$ explicitly:

$$\Phi_n(x_1, x_2, \dots, x_n, t) = C_n \quad .$$

The combination of the $n$ first integrals $\Phi_1$, $\Phi_2$,..., $\Phi_n$ allows us to find the functions $x_i(t)$ $(i=1, 2,..., n)$ that satisfy the system (5).

***Example:*** Consider a simplified form of the equation describing harmonic oscillation:

$$\frac{d^2x}{dt^2} + x = 0 .$$

This is written as an autonomous system of first-order ODEs,

$$\frac{dx}{dt} = y \quad , \quad \frac{dy}{dt} = -x \qquad (a)$$

(here, $x_1=x$, $x_2=y$). Eliminating $dt$, we have:

$$\frac{dx}{y} = -\frac{dy}{x} \qquad (b)$$

The solution of $(b)$ can be expressed in the form of a first integral, as follows:

$$xdx + ydy = 0 \;\Rightarrow\; d\,(x^2+y^2) = 0 \;\Rightarrow\; \Phi_1(x, y) \equiv x^2+y^2 = C_1^2 .$$

The first integral $\Phi_1$ determines the field lines of the velocity field

$$(dx/dt,\; dy/dt) \equiv (y, -x),$$

i.e., the phase curves. These are circles centered at the origin $O$ of $R^2$ (Fig. 5.2); obviously, they do not intersect with one another. For $C_1=0$ the field "line" is just a single point $O$, called the *equilibrium point* of the system. At this point the phase velocity vanishes at all $t$ : $(dx/dt, dy/dt) \equiv (y, -x) \equiv (0,0)$.

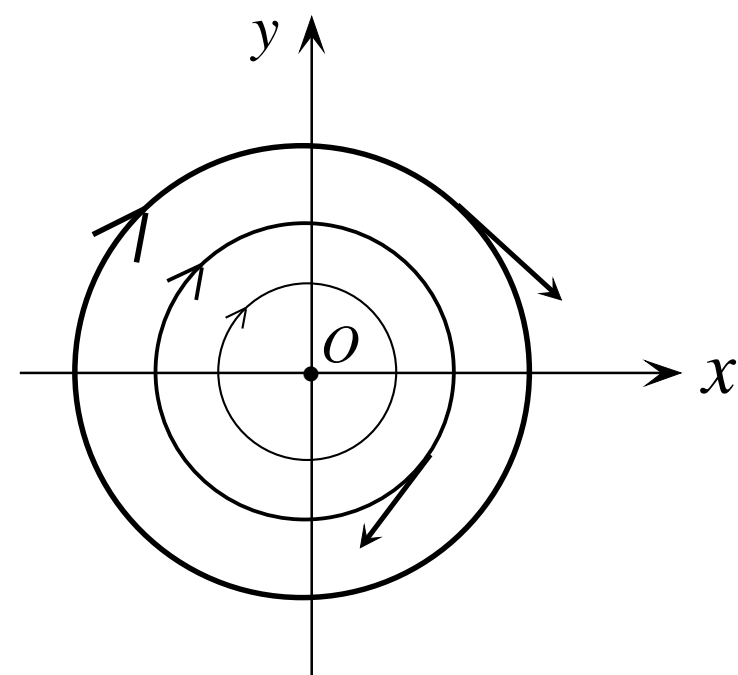

Fig. 5.2. Phase curves of the harmonic oscillator equation.

Now, for a complete solution of the autonomous system (*a*) we also need a first integral of the system directly dependent on *t*. As we have shown (see Example 2 in Sec. 4.1) this first integral is

$$\Phi_2(x, y, t) \equiv t + \arctan(y/x) = C_2 .$$

From $\Phi_1$ and $\Phi_2$ we get the solution of system (*a*):

$$x = C_1 \cos(t - C_2), \quad y = -C_1 \sin(t - C_2) .$$

The above relations describe the integral curves of system (*a*) in the space $R^2 \times R$: $(x, y, t)$. By eliminating the time *t* we find the phase curves of the system, which represent the field lines of the velocity field $F(x, y) \equiv (y, -x)$ and are the projections of the integral curves onto the phase space $R^2$: $(x, y)$. These phase curves are precisely the circles $x^2+y^2 = C_1^{\,2}$ and correspond to the solutions of Eq. (*b*).

Notice that the velocity field $(dx/dt, dy/dt) \equiv (y, -x)$, defined by the dynamical system, endows the phase curves with a sense of direction for increasing *t*. Analytically, let $dt > 0$. Then $dx > 0$ for $y > 0$ and $dx < 0$ for $y < 0$, while $dy > 0$ for $x < 0$ and $dy < 0$ for $x > 0$. That is, *x* increases (decreases) when *y* is positive (negative), while *y* increases (decreases) when *x* is negative (positive). This means that the curves are described *clockwise* for increasing *t*.

## 5.2 Geometric Significance of the First Integral

We consider the autonomous system of $n$ equations,

$$\frac{dx_i}{dt} = f_i(x_1, x_2, \cdots, x_n) \qquad (i = 1, 2, \cdots, n) \tag{1}$$

By eliminating *dt* we obtain a system of $(n-1)$ equations:

$$\frac{dx_1}{f_1(x_k)} = \frac{dx_2}{f_2(x_k)} = \cdots = \frac{dx_n}{f_n(x_k)} \tag{2}$$

Consider a first integral of system (2), of the form

$$\Phi(x_1, \ldots, x_n) = C \tag{3}$$

Relation (3) defines an $(n-1)$-*dimensional surface* in the phase space $R^n$. In fact, we have an infinite family of such surfaces, each surface corresponding to a given value of the constant $C$. For given initial conditions $x_i(t_0)=x_{0i}$ and by taking into account the fact that the value of the function $\Phi$ is constant (i.e., the same for all $t$), we calculate the corresponding value of $C$ as follows:

$$\Phi\big(x_1(t), \dots, x_n(t)\big) = C\,, \ \ \forall t \ \ \Rightarrow \ \ (\text{by putting } \ t=t_0)$$

$$\Phi\big(x_1(t_0), \dots, x_n(t_0)\big) = C \ \ \Rightarrow \ \ C = \Phi(x_{01}, \dots, x_{0n})\,.$$

Consider now a solution $x_i=x_i(t)$ $(i=1,\dots,n)$ [or, in vector form, $X=X(t)$] of system (1), for given initial conditions $x_i(t_0)=x_{0i}$ [or $X(t_0)=X_0$]. The points $X(t)\in R^n$ constitute a phase curve in the phase space $R^n$.

***Proposition:*** If the phase curve $X(t)$ has a common point with the $(n-1)$-dimensional surface (3): $\Phi(X)=C$, then this curve lies entirely on that surface.

***Proof:*** Consider a phase curve corresponding to the solution $X=X(t)$ of system (1) with initial condition $X(t_0)=X_0$. Moreover, assume that $\Phi(X_0)=C$. That is, $X_0$ is a common point of the phase curve and the surface. Given that the function $\Phi(X)$ assumes a constant value for all points $X(t)$ of a phase curve, we have that

$$\Phi\big(X(t)\big)=\Phi\big(X(t_0)\big)=\Phi(X_0)=C,$$

which means that all points of the phase curve $X(t)$ lie on the surface $\Phi(X)=C$. The uniqueness of the solution of the system for given initial conditions guarantees the impossibility of existence of any other phase curve passing through the same point $X_0$ of the surface. Indeed, such a curve, if it existed, would intersect with $X(t)$ at the point $X_0$, which cannot be the case for an *autonomous* system. In conclusion, a unique phase curve passes through the point $X_0$ of the surface $\Phi(X)=C$ and lies entirely on this surface.

Now, if we have $(n-1)$ independent first integrals of system (2), of the form $\Phi_i(X)=C_i$ $(i=1, 2,\dots, n-1)$, these define (for given $C_i$) a set of $(n-1)$-dimensional surfaces in $R^n$. The *intersection* of these surfaces is precisely the phase curve of the system that corresponds to the given initial conditions. [As mentioned earlier, these conditions determine the values of the constants $C_i$ by means of the relation $C_i = \Phi_i(X_0)$.] Indeed, note that this phase curve must simultaneously belong to all surfaces; hence, it must coincide with their intersection.

## 5.3 Vector Fields

At this point we need to introduce a few notational conventions:

1. For the coordinates of $R^n$, as well as for the components of vectors in $R^n$, we will be using *superscripts*. Thus, we will denote by $(x^1, x^2,\dots, x^n) \equiv (x^k)$ the coordinates of a point in space, and by $(V^1, V^2,\dots, V^n)$ the components of a vector $\vec{V}$.

2. We will be using the *summation convention*, according to which, if an expression contains the same index (say, $i$) as both a superscript *and* a subscript, then a *summation* of this expression is implied from $i=1$ to $n$. For example,

$$A^i B_i \equiv \sum_{i=1}^{n} A^i B_i = A^1 B_1 + A^2 B_2 + \cdots + A^n B_n \ .$$

3. The partial derivative with respect to $x^i$,

$$\frac{\partial}{\partial x^i} \equiv \partial_i \ ,$$

will be treated as an expression carrying a *subscript*. Thus, for example,

$$\sum_{i=1}^{n} V^i \frac{\partial \Phi}{\partial x^i} = \sum_{i=1}^{n} V^i \partial_i \Phi \equiv V^i \partial_i \Phi = V^i \frac{\partial \Phi}{\partial x^i} \ .$$

We note that the *name* of the repeated index is immaterial and may change without affecting the result of the summation. For example,

$$A^i B_i = A^j B_j = A^k B_k = \cdots \ ; \quad V^i \, \partial\Phi/\partial x^i = V^j \, \partial\Phi/\partial x^j = \cdots \ ; \quad \text{etc.}$$

Let, now, $\{\hat{e}_1, \cdots, \hat{e}_n\} \equiv \{\hat{e}_k\}$ be a basis of unit vectors in $R^n$. We consider a vector field $\vec{V}$ in this space:

$$\vec{V} = \sum_{i=1}^{n} V^i (x^1, \cdots, x^n) \, \hat{e}_i \equiv V^i (x^k) \, \hat{e}_i \qquad (1)$$

In component form,

$$\vec{V} \equiv \left( V^1(x^k), \cdots, V^n(x^k) \right).$$

We also consider the autonomous system of ODEs,

$$\frac{d}{dt} x^i(t) = V^i(x^k) \ , \quad x^i(0) = x_0{}^i \qquad (2)$$

(In the initial conditions we have put $t_0=0$.) We note that the phase curves $x^i(t)$ ($i=1,\ldots, n$) corresponding to the solutions of system (2) constitute the *field lines* of $\vec{V}$ (this field can be regarded as a velocity field if the variable $t$ represents time). At each point of a field line the vector $\vec{V}$ is *tangent* to this line. In other words, at each point of phase space $R^n$ the field $\vec{V}$ is tangential with respect to the (unique) field line passing through that point.

We now seek functions $\Phi(x^k)$ that *retain constant values along the field lines of* $\vec{V}$. That is, $\Phi(x^k)=C$ for solutions $x^i=x^i(t)$ of system (2). Obviously, every such function is a *first integral* of this system. Let $\Phi(x^k)$ be such a function. Then, $\Phi(x^k)=C \;\Rightarrow$

$$\frac{d}{dt}\Phi(x^k)=0 \;\Rightarrow\; \frac{\partial \Phi(x^k)}{\partial x^i}\frac{dx^i}{dt}=0 \quad (\textit{sum on } i\,!)\,.$$

Substituting for $dx^i/dt$ from system (2), we have:

$$V^i(x^k)\,\frac{\partial \Phi(x^k)}{\partial x^i}=0 \tag{3}$$

(Careful: By $x^k$ we collectively denote the whole *set* of variables $x^1,\dots,x^n$; thus $x^k$ does *not* represent any particular variable. This means that we do *not* sum with respect to $k$ but only with respect to $i$, since it is only the latter index that appears both as a superscript and as a subscript.) We conclude that if $\Phi(x^k)=C$ is a first integral of system (2), then the function $z=\Phi(x^k)$ is a solution of the homogeneous linear PDE

$$V^i(x^k)\,\frac{\partial z}{\partial x^i}=0 \tag{4}$$

***Example:*** Consider the vector field $\vec{V}\equiv(y,-x)$ in $R^2$: $(x, y)$. Here, $(x^1, x^2)\equiv(x, y)$ and $(V^1, V^2)\equiv(y, -x)$. The field lines are determined by the system

$$dx/dt = y\,,\quad dy/dt=-x \qquad (a)$$

We seek a first integral of this system, of the form $\Phi(x, y)=C$. The function $z=\Phi(x, y)$ will then be a solution of the PDE (4), namely, of

$$V^1\frac{\partial z}{\partial x^1}+V^2\frac{\partial z}{\partial x^2}\equiv y\frac{\partial z}{\partial x}-x\frac{\partial z}{\partial y}=0 \qquad (b)$$

The characteristic system of the above PDE is

$$\frac{dx}{y}=\frac{dy}{-x}=\frac{dz}{0} \qquad (c)$$

We seek two independent first integrals of system $(c)$:

$$dz=0 \;\Rightarrow\; z=C_1\,,$$

$$xdx+ydy=0 \;\Rightarrow\; d(x^2+y^2)=0 \;\Rightarrow\; x^2+y^2=C_2\,.$$

The general solution of the PDE $(b)$ is

$$F(C_1, C_2)=0 \;\; (\textit{arbitrary } F) \;\Rightarrow\; F(z, x^2+y^2)=0 \;\Rightarrow\; z=G(x^2+y^2) \;\; (\textit{arbitrary } G)\,.$$

In particular, by choosing $z=x^2+y^2\equiv\Phi(x,y)$ we find a first integral of the system $(a)$ that determines the field lines of $\vec{V}\equiv(y,-x)$:

$$\Phi(x, y)\equiv x^2+y^2=C\,.$$

We notice that $d\Phi/dt=0$ when the functions $x(t)$ and $y(t)$ are solutions of system ($a$) (show this!). This means that the function $\Phi(x, y)$ retains a constant value along any field line of $\vec{V}$ .

## 5.4 Differential Operators and Lie Derivative

Consider a vector field $\vec{V}$ in $R^n$ :

$$\vec{V} = V^i(x^k)\,\hat{e}_i \equiv \left(V^1(x^k), \cdots, V^n(x^k)\right) \tag{1}$$

The field lines are determined by the system of ODEs

$$\frac{dx^i}{dt} = V^i(x^k) \ , \quad x^i(0) = x_0^{\ i} \tag{2}$$

Given a function $f(x^1,..., x^n) \equiv f(x^k)$ in $R^n$, we would like to know the rate of change of $f$ along the field lines of the field (1); that is, the rate at which the value $f(x^k)$ of this function changes when the $x^i(t)$ are solutions of the system (2).

Let $x^i(t)$ $(i=1,..., n)$ be a solution of system (2). This solution corresponds to a field line of $\vec{V}$ . Along this line, the function $f$ takes on the values $f\left(x^k(t)\right)$. The rate of change of $f(x^k)$ along the field line is given by the directional derivative

$$\frac{d}{dt} f(x^k) = \frac{\partial f(x^k)}{\partial x^i} \frac{dx^i}{dt} \quad .$$

Substituting for $dx^i/dt$ from system (2), we have:

$$\frac{d}{dt} f(x^k) = V^i(x^k) \frac{\partial f(x^k)}{\partial x^i} \tag{3}$$

(sum on $i$ only!). Relation (3) is written

$$\frac{d}{dt} f(x^k) = \left( V^i(x^k) \frac{\partial}{\partial x^i} \right) f(x^k) = \left( V^i(x^k)\,\partial_i \right) f(x^k) \tag{3΄}$$

According to (3΄), the rate of change of the function $f(x^k)$ along the field lines of the field (1) is given quantitatively by the result of the action of the *differential operator*

$$V^i(x^k) \frac{\partial}{\partial x^i} = V^i(x^k)\,\partial_i$$

on $f(x^k)$.

We notice a one-to-one correspondence between vector fields and differential operators:

$$\vec{V} = V^i(x^k)\,\hat{e}_i \quad \leftrightarrow \quad V^i(x^k)\,\frac{\partial}{\partial x^i} = V^i(x^k)\,\partial_i \ .$$

Furthermore, the basis vectors $\hat{e}_i$ and the partial derivatives $\partial_i$ obey similar transformation rules under changes $\{x^k\}\rightarrow\{y^k\}$, $\{\hat{e}_k\}\rightarrow\{\hat{h}_k\}$, of the coordinate system [3-5]. Analytically,

$$\hat{h}_j = \frac{\partial x^i}{\partial y^j}\,\hat{e}_i \ , \qquad \frac{\partial}{\partial y^j} = \frac{\partial x^i}{\partial y^j}\,\frac{\partial}{\partial x^i} \ .$$

These observations suggest a new understanding of the vector concept in $R^n$: We no longer distinguish between the vector field $\vec{V}$ and the corresponding differential operator $V^i(x^k)\,\partial_i$ but we regard the two objects as being "identical"! Thus, we *define* the vector $\vec{V}$ as the differential operator

$$\vec{V} \equiv V^i(x^k)\,\frac{\partial}{\partial x^i} = V^i(x^k)\,\partial_i \tag{4}$$

Relation (3΄), then, takes on the new form

$$\frac{d}{dt}\,f(x^k) = \vec{V} f(x^k) \tag{5}$$

Now, in the case where a function $\Phi(x^k)$ is a *first integral* of system (2), the value of $\Phi$ is constant along any field line of $\vec{V}$, so that

$$\Phi(x^k) = C \ \Leftrightarrow\ \frac{d}{dt}\,\Phi(x^k) = 0 \ \Leftrightarrow\ \vec{V}\,\Phi(x^k) = V^i(x^k)\,\frac{\partial\Phi(x^k)}{\partial x^i} = 0 \tag{6}$$

We thus recover the homogeneous linear PDE found in Sec. 5.3.

***Definition:*** The derivative of a function $f(x^k)$ along the field lines of a vector field $\vec{V}$ in $R^n$ is called the *Lie derivative* of $f(x^k)$ with respect to $\vec{V}$, denoted $L_{\vec{V}}\, f(x^k)$.

According to (4) and (5) we can now write

$$\frac{d}{dt}\,f(x^k) \equiv L_{\vec{V}}\, f(x^k) = \vec{V} f(x^k) = V^i(x^k)\,\partial_i\, f(x^k) \tag{7}$$

In particular, according to (6) the Lie derivative with respect to $\vec{V}$, of a first integral $\Phi(x^k)$ of the system (2), is zero: $L_{\vec{V}}\,\Phi(x^k) = 0$.

***Comment:*** The introduction of the symbol $L_{\vec{V}}$ may seem superfluous since this operator appears to do the same job as the operator in (4). This coincidence, however, is valid only in the case of scalar functions of the form $f(x^k)$. The Lie derivative is a much more general concept of differential geometry and its mathematical expression varies in accordance with the tensor character of the function on which this derivative acts [3-5].

***Exercise:*** For the special case $f(x^k)=x^j$ (for given $j$), show that

$$L_{\vec{V}}\, x^j = V^j(x^k) \ .$$

(*Hint:* Notice that $\partial x^j/\partial x^i = \delta_{ij}$ .)

### 5.5 Exponential Solution of an Autonomous System

Consider the autonomous system of ODEs,

$$\frac{dx^i(t)}{dt} = V^i\left(x^1(t),\cdots,x^n(t)\right) \equiv V^i\left(x^k(t)\right) \ , \quad x^i(0) = x_0^{\ i} \tag{1}$$

($i$=1,..., $n$). The solution of this system will depend on $n$ parameters that, in turn, are dependent upon the initial values $x_0^{\ i}$. Hence, this solution will be expressed by a set of functions of the form

$$x^i = \Phi^i(t, x_0^{\ 1}, \ldots , x_0^{\ n}) \equiv \Phi^i(t, x_0^{\ k}) \tag{2}$$

where, by the initial conditions of the problem,

$$\Phi^i(0, x_0^{\ k}) = x_0^{\ i} \tag{3}$$

As a first step toward an analytic expression for the solution of system (1), we now define the differential operator

$$D_V = V^i(x_0^{\ k})\,\frac{\partial}{\partial x_0^{\ i}} \tag{4}$$

This operator acts on functions $f(x_0^{\ k})$ as follows:

$$D_V f(x_0^{\ k}) = V^i(x_0^{\ k})\,\frac{\partial f(x_0^{\ k})}{\partial x_0^{\ i}}$$

(sum on $i$ only!). We also define the *exponential operator*

$$e^{tD_V} \equiv \exp(tD_V) = \sum_{l=0}^{\infty}\frac{1}{l!}(tD_V)^l = 1 + tD_V + \frac{t^2}{2!}D_V^{\ 2} + \frac{t^3}{3!}D_V^{\ 3} + \cdots \tag{5}$$

where $D_V{}^2 f \equiv D_V(D_V f)$, etc. In particular, for $t$=0 we have the unit operator $e^0$=1.

We are now in a position to write the analytic expression for the solution (2) of system (1). As can be proven [5] this solution can be written in power-series form, as follows:

$$x^i = \Phi^i(t, x_0{}^k) = e^{tD_V} x_0{}^i = \left\{ \exp\left( t\, V^j(x_0{}^k) \frac{\partial}{\partial x_0{}^j} \right) \right\} x_0{}^i \tag{6}$$

(here, sum on $j$). Analytically,

$$x^i = \left[ 1 + t\,D_V + (t^2/2!)\,D_V{}^2 + (t^3/3!)\,D_V{}^3 + \cdots \right] x_0{}^i$$

$$= x_0{}^i + t\,D_V\, x_0{}^i + (t^2/2)\,D_V(D_V\, x_0{}^i) + \cdots$$

By taking into account that $\partial x_0{}^i / \partial x_0{}^j = \delta_{ij}$, we have:

$$D_V\, x_0{}^i = (V^j \partial/\partial x_0{}^j)\, x_0{}^i = V^j\, \partial x_0{}^i / \partial x_0{}^j = V^i(x_0{}^k)\,.$$

Thus, finally,

$$x^i = x_0{}^i + t\,V^i(x_0{}^k) + \frac{t^2}{2} V^j(x_0{}^k) \frac{\partial V^i(x_0{}^k)}{\partial x_0{}^j} + \cdots \tag{7}$$

Let us now make a little change of notation in system (1). Specifically, in place of $x_0{}^i$ we simply write $x^i$, while in place of $x^i(t)$ we set $\overline{x}^i(t)$. That is, we have a system of equations for the unknown functions $\overline{x}^i(t)$ with initial values $x^i$:

$$\frac{d\overline{x}^i(t)}{dt} = V^i\left( \overline{x}^k(t) \right) , \quad \overline{x}^i(0) = x^i \tag{8}$$

The solution (6) is then written

$$\overline{x}^i = \Phi^i(t, x^k) = \left\{ \exp\left( t\, V^j(x^k) \frac{\partial}{\partial x^j} \right) \right\} x^i \tag{9}$$

Note that, in this new notation, $\Phi^i(0, x^k) = x^i$, as demanded by the initial conditions.

Relation (9) describes a parametric curve in $R^n$ which starts from the point ($x^1$,..., $x^n$) for $t$=0 and passes through the point $\left( \overline{x}^1(t), \cdots, \overline{x}^n(t) \right)$ for $t$ >0. This relation admits the following geometrical interpretation: The operator $\exp\left( tV^j(x^k) \partial/\partial x^j \right)$ *pushes* the point ($x^1$,..., $x^n$) of the curve to the point $\left( \overline{x}^1(t), \cdots, \overline{x}^n(t) \right)$.

Now, according to what was said in Sec. 5.4 on the equivalence between differential operators and vector fields, the operator $D_V = V^i(x^k)\partial/\partial x^i$ can be identified with the vector field

$$\vec{V} = V^i(x^k)\,\frac{\partial}{\partial x^i} = V^i(x^k)\,\partial_i \tag{10}$$

The field lines of the field (10) are described by the curves (9), at each point of which the field is tangential. These lines represent the phase curves of the autonomous system (8). In vector notation, we write

$$\overline{x}^i = \Phi^i(t, x^k) = e^{t\vec{V}}\,x^i \tag{11}$$

If the variable $t$ is given the physical interpretation of time, then the field (10) is *static*. Indeed, the $V^i$ do not depend explicitly on time but only implicitly, through the $x^k$. Thus, for any given value of the $x^k$ (i.e., at any point of phase space $R^n$) the field is constant in time; its change with respect to time is only due to a displacement along a phase curve within a time interval, resulting in a corresponding change of the coordinates $x^k$ themselves.

We conclude that the field lines of the vector field (10) are static (time-independent) and, moreover, they do not intersect. Indeed, if they did intersect we would have two or more tangent vectors at the same point of phase space. This would mean either that the vector field changes with time (i.e., is not static) or that the static field is not uniquely defined everywhere.

## 5.6 Vector Fields as Generators of Transformations

Consider a vector field

$$\vec{V} = V^i(x^1,\cdots,x^n)\,\frac{\partial}{\partial x^i} \equiv V^i(x^k)\,\partial_i \tag{1}$$

The field line $\overline{x}^i(t)$ $(i=1,\cdots,n)$, starting from the point $(x^1,..., x^n)$ of $R^n$ for $t$=0, is given by

$$\overline{x}^i = \Phi^i(t, x^k) = e^{t\vec{V}}\,x^i = \left\{\exp\left(t\,V^j(x^k)\,\partial_j\right)\right\}\,x^i \tag{2}$$

We may say that the operator $e^{t\vec{V}}$ pushes the point $(x^1,..., x^n)$ of $R^n$ to the point $\left(\overline{x}^1(t),\cdots,\overline{x}^n(t)\right)$ along the (unique) field line passing through $(x^1,..., x^n)$ .

Now, let $F(x^1,..., x^n) \equiv F(x^k)$ be a function in $R^n$. The replacement $x^i \to \overline{x}^i(t)$ $(i$=1,..., $n)$ leads to a new function $F_t$ such that

$$F_t(x^k) = F\left(\overline{x}^k(t)\right) \quad \text{with} \quad F_0(x^k) = F(x^k) \text{ for } t=0$$

[where we have taken into account that $\overline{x}^i(0) = x^i$, according to the initial conditions of the problem]. We say that the field (1) is the *generator of the transformation*

$$x^i \to \overline{x}^i(t) \ \ (i=1,\cdots,n) \ , \quad F(x^k) \to F_t(x^k) = F\left(\overline{x}^k(t)\right) \tag{3}$$

***Example:*** In the two-dimensional space $R^2 : (x^1, x^2) \equiv (x, y)$, we consider the vector field

$$\vec{V} = \alpha x \frac{\partial}{\partial x} + \beta \frac{\partial}{\partial y} \qquad (\alpha, \beta = const.) \ .$$

The field lines are given by the system

$$\frac{d\overline{x}}{dt} = \alpha \overline{x} \ , \quad \frac{d\overline{y}}{dt} = \beta \ , \quad \text{with} \quad (\overline{x}, \overline{y}) \equiv (x, y) \quad \text{for} \quad t=0 \ .$$

The direct solution of this system is easy:

$$\overline{x}(t) = e^{\alpha t} x \ , \quad \overline{y}(t) = y + \beta t \ .$$

Alternatively (but less simply in this case) we can use the general formula (2), according to which

$$\overline{x}(t) = e^{t\vec{V}} x \ , \quad \overline{y}(t) = e^{t\vec{V}} y \quad \text{where} \quad e^{t\vec{V}} = 1 + t\vec{V} + \frac{t^2}{2!}\vec{V}\vec{V} + \frac{t^3}{3!}\vec{V}\vec{V}\vec{V} + \cdots$$

By noting that

$$\vec{V}x = \alpha x \ , \quad \vec{V}\vec{V}x \equiv \vec{V}(\vec{V}x) = \alpha \vec{V}x = \alpha^2 x \ , \quad \vec{V}\vec{V}\vec{V}x \equiv \vec{V}\left(\vec{V}(\vec{V}x)\right) = \alpha^3 x \ , \ \cdots \ ,$$

$$\vec{V}y = \beta \ , \quad \vec{V}\vec{V}y \equiv \vec{V}(\vec{V}y) = 0 \ , \quad \vec{V}\vec{V}\vec{V}y \equiv \vec{V}\left(\vec{V}(\vec{V}y)\right) = 0 \ , \ \cdots \ ,$$

we have:

$$\overline{x}(t) = \left[1 + \alpha t + \frac{(\alpha t)^2}{2!} + \frac{(\alpha t)^3}{3!} + \cdots\right] x = e^{\alpha t} x \ , \qquad \overline{y}(t) = y + \beta t \ ,$$

as before. The transformation (3) of a function $F(x, y)$ in $R^2$ is written

$$F(x,y) \to F_t(x,y) = F\left(\overline{x}(t), \ \overline{y}(t)\right) = F\left(e^{\alpha t} x, \ y + \beta t\right) \ .$$

In general, the variable $t$ in Eq. (2) is called the *parameter of the transformation* (3). For *infinitesimal* $t$ we can make the approximation

$$e^{t\vec{V}} \simeq 1 + t\vec{V} \ .$$

Then, relation (2) yields

$$\bar{x}^i(t) \simeq (1+t\vec{V})\,x^i = \left(1+t\,V^j(x^k)\frac{\partial}{\partial x^j}\right)x^i \;\Rightarrow$$

$$\bar{x}^i(t) \simeq x^i + t\,V^i(x^k) \tag{4}$$

Thus, in the preceding example,

$$\bar{x}(t) \simeq x + t\cdot(\alpha x) = (1+\alpha t)\,x \;, \quad \bar{y}(t) \simeq y + t\cdot\beta = y + \beta t \;.$$

***Exercise:*** Show that, in infinitesimal form,

$$F\left(\bar{x}^k(t)\right) \simeq F(x^k) + t\vec{V}F(x^k) \tag{5}$$

[*Hint:* For infinitesimal changes $dx^k$ of the $x^k$, the change of the value of $F$ is approximately equal to the differential $dF = (\partial F/\partial x^i)\,dx^i$. By Eq. (4), $dx^i = tV^i(x^k)$ .]

## 5.7 Geometric Significance of First-Order PDEs

In the space $R^3$: $(x^1, x^2, x^3) \equiv (x, y, z)$ we consider the vector field (expressed here in standard form)

$$\begin{aligned}\vec{V} &= P(x,y,z)\,\hat{u}_x + Q(x,y,z)\,\hat{u}_y + R(x,y,z)\,\hat{u}_z \\ &\equiv \left(P(x,y,z),\; Q(x,y,z),\; R(x,y,z)\right)\end{aligned} \tag{1}$$

where $\hat{u}_x$, $\hat{u}_y$, $\hat{u}_z$ are the unit vectors on the axes $x$, $y$, $z$, respectively. The field lines of this field, at each point of which $\vec{V}$ is a tangent vector, are given by the solutions of the autonomous system of first-order ODEs,

$$\frac{dx}{dt} = P(x,y,z) \;,\quad \frac{dy}{dt} = Q(x,y,z) \;,\quad \frac{dz}{dt} = R(x,y,z) \tag{2}$$

By eliminating $dt$, we obtain a system of two equations for the field lines:

$$\frac{dx}{P(x,y,z)} = \frac{dy}{Q(x,y,z)} = \frac{dz}{R(x,y,z)} \tag{3}$$

The system (2) describes the curves parametrically: $x=x(t)$, $y=y(t)$, $z=z(t)$. The system (3), on the other hand, describes a curve as a geometric locus of points of $R^3$. These points constitute the image of the mapping $(t\in R) \rightarrow \left(x(t), y(t), z(t)\right)\in R^3$.

We consider, now, a surface $S$ in $R^3$ consisting of field lines of the field (1), these lines being determined by the system (3) (Fig. 5.3). Such a surface can be described mathematically in two ways: by directly expressing one coordinate of space, say $z$, in terms of the other two: $z=f\,(x,y)$, or, equivalently, by a more symmetric equation of the form $u(x,y,z)=C$. The surface $S$ has the following property: every vector $\vec{N}$ normal to this surface at any point of $S$, is normal to the field line passing through this point, thus normal to the vector $\vec{V}$ tangent to the line at this point. Hence, $\vec{N}\cdot\vec{V}=0$.

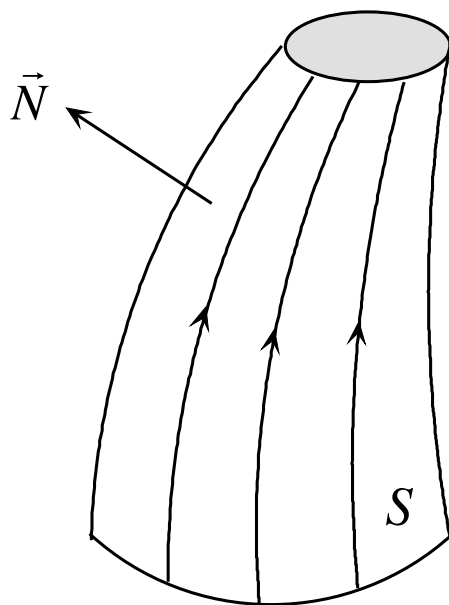

Fig. 5.3. A surface $S$ formed by field lines of a vector field.

As we know from vector analysis [6] a vector normal to the surface $u(x,y,z)=C$ is

$$\vec{N}=\vec{\nabla}u \equiv (\partial u/\partial x,\ \partial u/\partial y,\ \partial u/\partial z).$$

If the surface $S$ is described by the relation $z=f(x,y) \Rightarrow u(x,y,z)\equiv f(x,y)-z=0$, then

$$\vec{N}=\vec{\nabla}u \equiv (\partial f/\partial x,\ \partial f/\partial y,\ -1) \equiv (\partial z/\partial x,\ \partial z/\partial y,\ -1).$$

Given that $\vec{V}\equiv(P,Q,R)$, the orthogonality condition $\vec{N}\cdot\vec{V}=0 \Leftrightarrow \vec{V}\cdot\vec{\nabla}u=0$ can be expressed in the following ways:

1. With the homogeneous linear PDE

$$P(x,y,z)\frac{\partial u}{\partial x}+Q(x,y,z)\frac{\partial u}{\partial y}+R(x,y,z)\frac{\partial u}{\partial z}=0 \tag{4}$$

2. With the quasilinear PDE

$$P(x,y,z)\frac{\partial z}{\partial x}+Q(x,y,z)\frac{\partial z}{\partial y}=R(x,y,z) \tag{5}$$

In conclusion, the solutions of the PDEs (4) and (5) represent surfaces in $R^3$, formed by field lines of the field (1). We note that the system (3) is the *characteristic system* of the quasilinear PDE (5) (cf. Sec. 4.2).

We also note that the solution $u(x,y,z)$ of the PDE (4) is a *first integral* of the autonomous system (2), since $u(x,y,z)=C$ for solutions $x(t)$, $y(t)$, $z(t)$ of this system. Indeed, these solutions correspond to field lines of the field (1), each line lying entirely on some surface $u(x,y,z)=C$.

Alternatively, let us notice that the PDE (4) is written

$$\vec{V}u(x,y,z)=0\,,$$

where $\vec{V}$ now denotes the differential operator corresponding to the vector field (1):

$$\vec{V}=P\,\partial/\partial x+Q\,\partial/\partial y+R\,\partial/\partial z\,.$$

As we showed in Sec. 5.4, the above PDE is satisfied when the function $u(x,y,z)$ is a first integral of system (2), in which case the value of $u$ is constant along any field line of the field (1).

***Example:*** For $P=y$, $Q=-x$, $R=0$, we have the vector field

$$\vec{V}=y\,\partial/\partial x-x\,\partial/\partial y\;(+\,0\cdot\partial/\partial z)\;.$$

The surfaces $S$, consisting of field lines of this field, are given by the solutions of the PDE (5), written here as

$$y\partial z/\partial x-x\partial z/\partial y=0.$$

These solutions are given by the relation $z=F(x^2+y^2)$ with arbitrary $F$ (see Sec. 4.2, Example 2).

Assume now that we seek a surface $S_n$ intersecting the field lines of the vector field (1) in such a way that this field is *normal* to the surface at each point of $S_n$ . If $S_n$ is described by a relation of the form $U(x,y,z)=C$, then, without loss of generality, we may identify the normal vector $\vec{N}=\vec{\nabla}U$ at each point of $S_n$ with the field vector $\vec{V}$ :

$$\vec{\nabla}U(x,y,z)=\vec{V}(x,y,z) \tag{6}$$

In component form, relation (6) corresponds to the system of PDEs

$$\frac{\partial U}{\partial x}=P\;,\quad \frac{\partial U}{\partial y}=Q\;,\quad \frac{\partial U}{\partial z}=R \tag{7}$$

As we know (see Sec. 1.4) the integrability condition of (6) and (7) for solution for $U$ is written

$$\vec{\nabla}\times\vec{V}=0\;\Leftrightarrow\;\frac{\partial P}{\partial y}=\frac{\partial Q}{\partial x}\;,\quad\frac{\partial P}{\partial z}=\frac{\partial R}{\partial x}\;,\quad\frac{\partial Q}{\partial z}=\frac{\partial R}{\partial y} \tag{8}$$

We note that if a surface $S_n$ : $U(x,y,z)=C$, normal to the field lines of the vector field (1), exists, then this field is *potential* (Sec. 1.4). The function $U(x,y,z)$ represents the potential function of the field, according to (6).

***Application:*** The *electrostatic field* $\vec{E}$ is irrotational:

$$\vec{\nabla}\times\vec{E}=0\,;$$

thus it satisfies the necessary condition for existence of an *electrostatic potential* $U(x,y,z)$ such that

$$\vec{E} = -\vec{\nabla} U \ .$$

An *equipotential surface* $S_n$ is the geometric locus of points of space that satisfy the relation $U(x,y,z)=C$. Such a surface intersects normally the electric field lines, since at each point of $S_n$ the field vector is normal to the surface (explain why).

# CHAPTER 6

# INTEGRABLE SYSTEMS OF PARTIAL DIFFERENTIAL EQUATIONS

## 6.1 Notation

Let $u(x,t)$ be a function of two variables. For the partial derivatives of $u$ the following notation will be used:

$$\frac{\partial u}{\partial x} = \partial_x u = u_x \ , \quad \frac{\partial u}{\partial t} = \partial_t u = u_t \ , \quad \frac{\partial^2 u}{\partial x^2} = u_{xx} \ , \quad \frac{\partial^2 u}{\partial t^2} = u_{tt} \ , \quad \frac{\partial^2 u}{\partial x \partial t} = u_{xt} \ ,$$

etc. In general, a subscript will denote partial differentiation with respect to the indicated variable.

Consider now a function $F$ of $x$, $t$, $u$, as well as of a number of partial derivatives of $u$. We will denote this type of dependence by writing

$$F(x,t,u,u_x,u_t,u_{xx},u_{tt},u_{xt},\cdots) \equiv F[u] \ .$$

We also write

$$F_x = \partial_x F = \partial F / \partial x \quad , \quad F_t = \partial_t F = \partial F / \partial t \quad , \quad F_u = \partial_u F = \partial F / \partial u \quad ,$$

etc. Note that in determining $F_x$ and $F_t$ we must take into account both the *explicit* and the *implicit* (through $u$ and its partial derivatives) dependence of $F$ on $x$ and $t$. As an example, for $F[u] = 3xtu^2$ we have $F_x = 3tu^2 + 6xtuu_x$ , $F_t = 3xu^2 + 6xtuu_t$ .

## 6.2 Bäcklund Transformations

Consider two partial differential equations (PDEs) $P[u]=0$ and $Q[v]=0$ for the unknown functions $u$ and $v$, respectively, where the bracket notation introduced in the previous section is adopted. Both $u$ and $v$ are functions of two variables $x$, $t$.

Independently, for the moment, consider also a pair of coupled PDEs for $u$ and $v$,

$$B_1[u,v] = 0 \quad (a) \qquad B_2[u,v] = 0 \quad (b) \tag{1}$$

where the expressions $B_i\,[u,v]$ ($i$=1,2) may contain $u$, $v$ as well as partial derivatives of $u$ and $v$ with respect to $x$ and $t$. We notice that $u$ appears in both equations ($a$) and ($b$). The question then is: if we find an expression for $u$ by integrating ($a$) for a given $v$, will it match the corresponding expression for $u$ found by integrating ($b$) for the same $v$? The answer is that, in order that ($a$) and ($b$) be consistent with each other for solution for $u$, the function $v$ must be properly chosen so as to satisfy a certain *consistency condition* (or *integrability condition* or *compatibility condition*).

By a similar reasoning, in order that (*a*) and (*b*) in (1) be mutually consistent for solution for $v$, for some given $u$, the function $u$ must now itself satisfy a corresponding integrability condition.

If it happens that the two consistency conditions for integrability of the system (1) are precisely the PDEs $P[u]=0$ and $Q[v]=0$, we say that the above system constitutes a *Bäcklund transformation* (BT) connecting solutions of $P[u]=0$ with solutions of $Q[v]=0$. In the special case where $P$ and $Q$ are functionally identical, i.e., if $u$ and $v$ satisfy the *same* PDE, the system (1) is an *auto-Bäcklund* transformation (auto-BT) for this PDE.

Suppose now that we seek solutions of the PDE $P[u]=0$. Assume that we are able to find a BT connecting solutions $u$ of this equation with solutions $v$ of the PDE $Q[v]=0$ (if $P$ and $Q$ are functionally identical, the auto-BT connects solutions $u$ and $v$ of the same PDE). Let $v=v_0(x,t)$ be some known solution of $Q[v]=0$. The BT is then a system of PDEs for the unknown $u$,

$$B_i\,[u, v_0] = 0 \quad (i = 1, 2) \tag{2}$$

The system (2) is integrable for $u$, given that the function $v_0$ satisfies *a priori* the required integrability condition $Q[v]=0$. The solution $u$ then of the system satisfies the PDE $P[u]=0$. Thus a solution $u(x,t)$ of the latter PDE is found without actually solving the equation itself, simply by integrating the BT (2) with respect to $u$. Of course, this method will be useful provided that integrating the system (2) for $u$ is simpler than integrating the PDE $P[u]=0$ itself. If the transformation (2) is an auto-BT for the PDE $P[u]=0$, then, starting with a known solution $v_0(x,t)$ of this equation and integrating the system (2) we find another solution $u(x,t)$ of the same equation.

***Examples:***

***1.*** The *Cauchy-Riemann relations* of Complex Analysis (cf. Sec. 2.1),

$$u_x = v_y \quad (a) \qquad u_y = -v_x \quad (b) \tag{3}$$

(where here the variable $t$ has been renamed $y$) constitute an auto-BT for the *Laplace equation*,

$$P\,[w] \equiv w_{xx} + w_{yy} = 0 \tag{4}$$

Let us explain this: Suppose we want to solve the system (3) for $u$, for a given choice of the function $v(x,y)$. To see if the PDEs (*a*) and (*b*) match for solution for $u$, we must compare them in some way. We thus differentiate (*a*) with respect to $y$ and (*b*) with respect to $x$, and equate the mixed derivatives of $u$. That is, we apply the integrability condition $(u_x)_y= (u_y)_x$ . In this way we eliminate the variable $u$ and find the condition that must be obeyed by $v(x,y)$:

$$P\,[v] \equiv v_{xx} + v_{yy} = 0 \ .$$

Similarly, by using the integrability condition $(v_x)_y= (v_y)_x$ to eliminate $v$ from the system (3), we find the necessary condition in order that this system be integrable for $v$, for a given function $u(x,y)$:

$$P[u] \equiv u_{xx} + u_{yy} = 0 \ .$$

In conclusion, the integrability of system (3) with respect to either variable requires that the other variable satisfy the Laplace equation (4).

Let now $v_0(x,y)$ be a known solution of the Laplace equation (4). Substituting $v=v_0$ in the system (3), we can integrate this system with respect to $u$. It is not hard to show (by eliminating $v_0$ from the system) that the solution $u$ will also satisfy the Laplace equation (4). As an example, by choosing the solution $v_0(x,y)=xy$, we find a new solution $u(x,y)=(x^2-y^2)/2+C$.

***2.*** The *Liouville equation* is written

$$P[u] \equiv u_{xt} - e^u = 0 \quad \Leftrightarrow \quad u_{xt} = e^u \tag{5}$$

Due to its nonlinearity, this PDE is hard to integrate directly. A solution is thus sought by means of a BT. We consider an auxiliary function $v(x,t)$ and an associated PDE,

$$Q[v] \equiv v_{xt} = 0 \tag{6}$$

We also consider the system of first-order PDEs,

$$u_x + v_x = \sqrt{2}\ e^{(u-v)/2} \quad (a) \qquad u_t - v_t = \sqrt{2}\ e^{(u+v)/2} \quad (b) \tag{7}$$

Differentiating the PDE ($a$) with respect to $t$ and the PDE ($b$) with respect to $x$, and eliminating $(u_t - v_t)$ and $(u_x + v_x)$ in the ensuing equations with the aid of ($a$) and ($b$), we find that $u$ and $v$ satisfy the PDEs (5) and (6), respectively. Thus, the system (7) is a BT connecting solutions of (5) and (6). Starting with the trivial solution $v=0$ of (6), and integrating the system

$$u_x = \sqrt{2}\ e^{u/2} \ , \qquad u_t = \sqrt{2}\ e^{u/2} \ ,$$

we find a nontrivial solution of (5):

$$u(x,t) = -2\ln\left( C - \frac{x+t}{\sqrt{2}} \right) .$$

***3.*** The "*sine-Gordon" equation* has applications in various areas of Physics, e.g., in the study of crystalline solids, in the transmission of elastic waves, in magnetism, in elementary-particle models, etc. The equation (whose name is a pun on the related linear Klein-Gordon equation) is written

$$P[u] \equiv u_{xt} - \sin u = 0 \quad \Leftrightarrow \quad u_{xt} = \sin u \tag{8}$$

The following system of equations is an auto-BT for the nonlinear PDE (8):

$$\frac{1}{2}(u+v)_x = a \sin\left( \frac{u-v}{2} \right) \ , \qquad \frac{1}{2}(u-v)_t = \frac{1}{a} \sin\left( \frac{u+v}{2} \right) \tag{9}$$

where $a$ ($\neq 0$) is an arbitrary real constant. [Because of the presence of $a$, the system (9) is called a *parametric* BT.] When $u$ is a solution of (8) the BT (9) is integrable for $v$, which, in turn, also is a solution of (8): $P[v]=0$; and vice versa. Starting with the trivial solution $v=0$ of $v_{xt}=\sin v$, and integrating the system

$$u_x = 2a\sin\frac{u}{2} \quad , \quad u_t = \frac{2}{a}\sin\frac{u}{2} \quad ,$$

we obtain a new solution of (8):

$$u(x,t) = 4\arctan\left\{ C\exp\left( ax + \frac{t}{a} \right) \right\} \ .$$

## 6.3 Lax Pair for a Nonlinear PDE

Let $F[u]=0$ be a nonlinear PDE for $u(x,t)$. Independently, for the moment, consider also a pair of *linear* PDEs for a new variable $\psi$, in which pair the variable $u$ enters as a sort of "parametric" function. We write

$$L_1(\psi\,;u) = 0 \quad , \quad L_2(\psi\,;u) = 0 \tag{1}$$

In order for the system (1) to be integrable for $\psi$ [i.e., for the two PDEs in (1) to be compatible with each other for solution for $\psi$] the function $u(x,t)$ must be properly chosen. We now make the special assumption that the linear system (1) is integrable for $\psi$ on the condition that $u$ satisfies the given nonlinear PDE $F[u]=0$. In this case the system (1) constitutes a *Lax pair* for $F[u]=0$. The construction of a Lax pair is closely related to a method of integration of nonlinear PDEs, called the *inverse scattering method* (see, e.g., [1,2]).

***Examples:***

***1.*** The *Korteweg-de Vries* (KdV) equation describes the propagation of particle-like nonlinear waves called *solitons* [1,2]. One form of this equation is

$$F[u] \equiv u_t - 6u\,u_x + u_{xxx} = 0 \tag{2}$$

The Lax pair for the nonlinear PDE (2) is written

$$\psi_{xx} = (u-\lambda)\psi \quad (a) \qquad \psi_t = 2(u+2\lambda)\psi_x - u_x\psi \quad (b) \tag{3}$$

where $\lambda$ is an arbitrary parameter. For system (3) to be integrable for $\psi$, equations ($a$) and ($b$) must agree with each other for all values of $\lambda$. Thus, in particular, the mixed derivative $(\psi_{xx})_t$ from ($a$) must match the derivative $(\psi_t)_{xx}$ from ($b$). The corresponding integrability condition is, therefore, $(\psi_{xx})_t=(\psi_t)_{xx}$. Performing suitable differentiations of ($a$) and ($b$) and using these same equations to eliminate $\psi_{xx}$ and $\psi_t$, we obtain the relation

$$(u_t - 6u\,u_x + u_{xxx})\psi \equiv F[u]\,\psi = 0 \ .$$

Hence, in order for the system (3) to have a nontrivial solution $\psi \neq 0$, it is necessary that $F[u]=0$; that is, $u$ must satisfy the KdV equation (2).

***2.*** One form of the *chiral field equation* is

$$F[g] \equiv \partial_t(g^{-1}g_x) + \partial_x(g^{-1}g_t) = 0 \tag{4}$$

where $g=g(x,t)$ is a non-singular, complex ($n\times n$) matrix. This equation constitutes a two-dimensional reduction of the four-dimensional self-dual Yang-Mills equation [1]. The Lax pair for the nonlinear PDE (4) is written

$$\psi_t = \frac{\lambda}{1-\lambda}\, g^{-1} g_t \psi \quad (a) \qquad \psi_x = -\frac{\lambda}{1+\lambda}\, g^{-1} g_x \psi \quad (b) \tag{5}$$

where $\psi$ is a complex ($n\times n$) matrix and $\lambda$ is an arbitrary complex parameter. The compatibility of Eqs. ($a$) and ($b$) with each other requires that $(\psi_t)_x=(\psi_x)_t$ . Cross-differentiating ($a$) and ($b$), using these same relations to eliminate $\psi_x$ and $\psi_t$ , and finally eliminating the common factor $\psi$ (assuming $\psi\neq 0$), we find the relation

$$\partial_t (g^{-1} g_x) + \partial_x (g^{-1} g_t) - \lambda \left\{ \partial_t (g^{-1} g_x) - \partial_x (g^{-1} g_t) + [g^{-1} g_t ,\, g^{-1} g_x] \right\} = 0$$

where, in general, by $[A\,,B]\equiv AB-BA$ we denote the *commutator* of two matrices $A$ and $B$. As can be shown (see Appendix B) the quantity inside the curly brackets vanishes identically. Thus, in order for the system (5) to have a nontrivial solution for $\psi$, the matrix function $g$ must satisfy the PDE (4).

## 6.4 The Maxwell Equations as a Bäcklund Transformation

A somewhat different approach to the concept of a Bäcklund transformation (BT) has recently been suggested [3,4]. Specifically, rather than being an auxiliary tool for integrating a given (usually nonlinear) PDE, it is the BT itself (regarded as a system) whose solutions are sought. To this end, it is examined whether the PDEs expressing the integrability conditions of the BT possess known, parameter-dependent solutions. By properly matching the parameters it may then be possible to find *conjugate solutions* of these PDEs; solutions, that is, which jointly satisfy the BT. This method is particularly effective for BTs whose integrability conditions are *linear* PDEs.

A nice example of this scheme is furnished by the *Maxwell equations* of classical Electrodynamics [5,6]. In *empty space* where no sources (charges and/or currents) exist, these equations form a homogeneous linear system:

$$\begin{aligned} &(a) \quad \vec{\nabla}\cdot\vec{E} = 0 \qquad &(c) \quad \vec{\nabla}\times\vec{E} = -\frac{\partial \vec{B}}{\partial t} \\ &(b) \quad \vec{\nabla}\cdot\vec{B} = 0 \qquad &(d) \quad \vec{\nabla}\times\vec{B} = \varepsilon_0 \mu_0 \frac{\partial \vec{E}}{\partial t} \end{aligned} \tag{1}$$

where $\vec{E}, \vec{B}$ are the electric and the magnetic field, respectively, and where the $\varepsilon_0$ , $\mu_0$ are constants associated with the S.I. system of units.

We will show that, by the self-consistency of the Maxwell system (1), each field $\vec{E}$ and $\vec{B}$ satisfies a corresponding wave equation. In other words, the system (1) is a BT connecting two *separate* wave equations, one for the electric field and one of similar form for the magnetic field. Since the two fields are physically different (they have different physical properties and dimensions) this BT is *not* an auto-BT.

As can be checked [3,4] the only nontrivial integrability conditions for system (1) are

$$\vec{\nabla}\times(\vec{\nabla}\times\vec{E})=\vec{\nabla}(\vec{\nabla}\cdot\vec{E})-\nabla^2\vec{E} \tag{2}$$

and

$$\vec{\nabla}\times(\vec{\nabla}\times\vec{B})=\vec{\nabla}(\vec{\nabla}\cdot\vec{B})-\nabla^2\vec{B} \tag{3}$$

Taking the *rot* of (1)(*c*) and using (2) and (1)(*a*),(*d*), we have:

$$\vec{\nabla}\times(\vec{\nabla}\times\vec{E})=-\vec{\nabla}\times\frac{\partial\vec{B}}{\partial t} \;\Rightarrow\; \vec{\nabla}(\vec{\nabla}\cdot\vec{E})-\nabla^2\vec{E}=-\frac{\partial}{\partial t}(\vec{\nabla}\times\vec{B}) \;\Rightarrow$$

$$\nabla^2\vec{E}-\varepsilon_0\mu_0\frac{\partial^2\vec{E}}{\partial t^2}=0 \tag{4}$$

Similarly, taking the *rot* of (1)(*d*) and using (3) and (1)(*b*),(*c*), we find:

$$\nabla^2\vec{B}-\varepsilon_0\mu_0\frac{\partial^2\vec{B}}{\partial t^2}=0 \tag{5}$$

Thus, the consistency conditions of the Maxwell system (1) yield two separate second-order linear PDEs, one for each field $\vec{E}$ and $\vec{B}$. We conclude that the system (1) is a BT relating the *wave equations* (4) and (5) for the electric and the magnetic field, respectively.

Although of different physical content, Eqs. (4) and (5) share the common form

$$\nabla^2\vec{A}-\varepsilon_0\mu_0\frac{\partial^2\vec{A}}{\partial t^2}=0 \tag{6}$$

We set

$$\varepsilon_0\mu_0\equiv\frac{1}{c^2} \;\Leftrightarrow\; c=\frac{1}{\sqrt{\varepsilon_0\mu_0}} \tag{7}$$

(where $c$ is the speed of light in empty space) and we write (6) as

$$\nabla^2\vec{A}-\frac{1}{c^2}\frac{\partial^2\vec{A}}{\partial t^2}=0 \tag{8}$$

where $\vec{A}=\vec{E}$ or $\vec{B}$.

The wave equation (8) admits plane-wave solutions of the form

$$\vec{A}=\vec{F}(\vec{k}\cdot\vec{r}-\omega t) \quad \text{where} \quad \omega/k=c \quad \text{with} \quad k=|\vec{k}\,| \tag{9}$$

The simplest such solution is a *monochromatic plane wave* of angular frequency $\omega$, propagating in the direction of the *wave vector* $\vec{k}$ :

$$\begin{aligned} \vec{E}(\vec{r},t) &= \vec{E}_0\, e^{\,i(\vec{k}\cdot\vec{r}-\omega t)} \qquad (a) \\ \vec{B}(\vec{r},t) &= \vec{B}_0\, e^{\,i(\vec{k}\cdot\vec{r}-\omega t)} \qquad (b) \end{aligned} \tag{10}$$

where the $\vec{E}_0, \vec{B}_0$ are constant *complex* amplitudes. (The term *"monochromatic"* indicates that the electromagnetic wave is a harmonic wave containing a single frequency $\omega$.) All constants appearing in Eqs. (10) (i.e., the amplitudes and the frequency), as well as the direction of the wave vector, can be chosen arbitrarily; thus these choices can be regarded as *parameters* on which the solutions (10) of the wave equations (4) and (5) depend.

It must be emphasized that, whereas every pair $(\vec{E},\vec{B})$ satisfying the Maxwell system (1) also satisfies the wave equation (8), the converse is not true. That is, the solutions (10) of the wave equation are not *a priori* solutions of the Maxwell equations. We must thus substitute the general solutions (10) into the system (1) to find the additional constraints that the latter system imposes on the parameters contained in Eqs. (10). By fixing these parameters, the wave solutions (*a*) and (*b*) in (10) will become *BT-conjugate* with respect to the Maxwell system (1).

To this end, we need two more vector identities: If $\Phi$ is a scalar field and if $\vec{A}$ is a vector field, then

$$\vec{\nabla}\cdot(\Phi\vec{A}) = (\vec{\nabla}\Phi)\cdot\vec{A} + \Phi(\vec{\nabla}\cdot\vec{A}) \;,$$
$$\vec{\nabla}\times(\Phi\vec{A}) = (\vec{\nabla}\Phi)\times\vec{A} + \Phi(\vec{\nabla}\times\vec{A}) \;.$$

In our case we set $\Phi = e^{\,i(\vec{k}\cdot\vec{r}-\omega t)} = e^{\,i\vec{k}\cdot\vec{r}}\, e^{-i\omega t}$ and $\vec{A}=\vec{E}_0$ or $\vec{B}_0$. We also note that

$$\vec{\nabla}\cdot\vec{E}_0 = \vec{\nabla}\cdot\vec{B}_0 = 0 \;, \quad \vec{\nabla}\times\vec{E}_0 = \vec{\nabla}\times\vec{B}_0 = 0 \qquad (\text{since } \vec{E}_0, \vec{B}_0 \text{ are constants}) \;,$$

$$\vec{\nabla}e^{\,i\vec{k}\cdot\vec{r}} = \Big(\hat{u}_x\frac{\partial}{\partial x} + \hat{u}_y\frac{\partial}{\partial y} + \hat{u}_z\frac{\partial}{\partial z}\Big)e^{\,i(k_x x + k_y y + k_z z)} = i\,(k_x\hat{u}_x + k_y\hat{u}_y + k_z\hat{u}_z)\,e^{\,i\vec{k}\cdot\vec{r}} = i\,\vec{k}\,e^{\,i\vec{k}\cdot\vec{r}} \;,$$

$$\frac{\partial}{\partial t}\,e^{-i\omega t} = -i\,\omega\, e^{-i\omega t} \;.$$

Substituting Eqs. (10)(*a*) and (*b*) into Eqs. (1)(*a*) and (*b*), respectively, we have:

$$(\vec{E}_0\, e^{-i\omega t})\cdot\vec{\nabla}e^{\,i\vec{k}\cdot\vec{r}} = 0 \;\Rightarrow\; (\vec{k}\cdot\vec{E}_0)\, e^{\,i(\vec{k}\cdot\vec{r}-\omega t)} = 0 \,,$$
$$(\vec{B}_0\, e^{-i\omega t})\cdot\vec{\nabla}e^{\,i\vec{k}\cdot\vec{r}} = 0 \;\Rightarrow\; (\vec{k}\cdot\vec{B}_0)\, e^{\,i(\vec{k}\cdot\vec{r}-\omega t)} = 0 \,,$$

so that

$$\vec{k}\cdot\vec{E}_0 = 0 \;, \quad \vec{k}\cdot\vec{B}_0 = 0 \,. \tag{11}$$

Multiplying by $e^{i(\vec{k}\cdot\vec{r}-\omega t)}$ and using Eqs. (10), we find

$$\vec{k}\cdot\vec{E}=0\,,\quad \vec{k}\cdot\vec{B}=0 \tag{12}$$

This indicates that in a monochromatic plane electromagnetic wave the oscillating fields $\vec{E}$ and $\vec{B}$ are always normal to the wave vector $\vec{k}$, that is, normal to the direction of propagation of the wave. Thus, this wave is a *transverse* wave.

Next, substituting Eqs. (10)(*a*) and (*b*) into Eqs. (1)(*c*) and (*d*), we have:

$$\begin{aligned}
&e^{-i\omega t}\,(\vec{\nabla}e^{i\vec{k}\cdot\vec{r}})\times\vec{E}_0 = i\,\omega\,\vec{B}_0\,e^{i(\vec{k}\cdot\vec{r}-\omega t)} \Rightarrow\\
&\qquad\qquad (\vec{k}\times\vec{E}_0)\,e^{i(\vec{k}\cdot\vec{r}-\omega t)} = \omega\,\vec{B}_0\,e^{i(\vec{k}\cdot\vec{r}-\omega t)}\,,\\
&e^{-i\omega t}\,(\vec{\nabla}e^{i\vec{k}\cdot\vec{r}})\times\vec{B}_0 = -i\,\omega\,\varepsilon_0\,\mu_0\,\vec{E}_0\,e^{i(\vec{k}\cdot\vec{r}-\omega t)} \Rightarrow\\
&\qquad\qquad (\vec{k}\times\vec{B}_0)\,e^{i(\vec{k}\cdot\vec{r}-\omega t)} = -\frac{\omega}{c^2}\,\vec{E}_0\,e^{i(\vec{k}\cdot\vec{r}-\omega t)}\,,
\end{aligned}$$

so that

$$\vec{k}\times\vec{E}_0=\omega\,\vec{B}_0\,,\quad \vec{k}\times\vec{B}_0=-\frac{\omega}{c^2}\,\vec{E}_0 \tag{13}$$

Multiplying by $e^{i(\vec{k}\cdot\vec{r}-\omega t)}$ and using (10), we find

$$\vec{k}\times\vec{E}=\omega\,\vec{B}\,,\quad \vec{k}\times\vec{B}=-\frac{\omega}{c^2}\,\vec{E} \tag{14}$$

We notice that at each instant the fields $\vec{E}$ and $\vec{B}$ are normal to each other as well as normal to the direction of propagation $\vec{k}$ of the wave.

Let us now assume that the complex amplitudes $\vec{E}_0\,,\vec{B}_0$ can be written as

$$\vec{E}_0=\vec{E}_{0,R}\,e^{i\alpha}\,,\quad \vec{B}_0=\vec{B}_{0,R}\,e^{i\beta}$$

where the $\vec{E}_{0,R}\,,\vec{B}_{0,R}$ are *real* vectors and where $\alpha,\beta$ are real numbers (physically, this indicates a *linearly polarized* wave). As we can show, relations (13) then demand that $\alpha=\beta$ and that

$$\vec{k}\times\vec{E}_{0,R}=\omega\,\vec{B}_{0,R}\,,\quad \vec{k}\times\vec{B}_{0,R}=-\frac{\omega}{c^2}\,\vec{E}_{0,R} \tag{15}$$

The monochromatic waves (10) are written

$$\vec{E}=\vec{E}_{0,R}\,e^{i(\vec{k}\cdot\vec{r}-\omega t+\alpha)}\,,\quad \vec{B}=\vec{B}_{0,R}\,e^{i(\vec{k}\cdot\vec{r}-\omega t+\alpha)} \tag{16}$$

Taking the real parts of Eqs. (16) we find the expressions for the *real* fields $\vec{E}$ and $\vec{B}$:

$$\vec{E} = \vec{E}_{0,R} \cos(\vec{k} \cdot \vec{r} - \omega t + \alpha) \ , \quad \vec{B} = \vec{B}_{0,R} \cos(\vec{k} \cdot \vec{r} - \omega t + \alpha) \tag{17}$$

Note, in particular, that the two fields "oscillate" in phase, acquiring their maximum, minimum and zero values simultaneously.

Taking the magnitudes of the vector relations (15) and using the fact that the $\vec{E}_{0,R}$ and $\vec{B}_{0,R}$ are normal to the wave vector $\vec{k}$, as well as that $\omega/k=c$, we find

$$E_{0,R} = c\, B_{0,R} \tag{18}$$

where $E_{0,R}=|\vec{E}_{0,R}|$ and $B_{0,R}=|\vec{B}_{0,R}|$. Also, taking the magnitudes of Eqs. (17) and using (18), we find a relation for the *instantaneous values* of the electric and the magnetic field:

$$E = c\, B \tag{19}$$

where $E=|\vec{E}|$ and $B=|\vec{B}|$.

## 6.5 Bäcklund Transformations as Recursion Operators

The concept of symmetries of partial differential equations (PDEs) is discussed in a number of books and articles (see, e.g., [7,8]). Let us briefly review the main ideas.

Consider a PDE $F[u]=0$, where, for simplicity, $u=u(x,t)$. A transformation

$$u(x,t) \ \rightarrow \ u'(x,t)$$

from the function $u$ to a new function $u'$ represents a *symmetry* of the given PDE if the following condition is satisfied: $u'(x,t)$ is a solution of $F[u]=0$ *if* $u(x,t)$ is a solution. That is,

$$F[u'] = 0 \quad \textit{when} \quad F[u] = 0 \tag{1}$$

An *infinitesimal symmetry transformation* is written

$$u' = u + \delta u = u + \alpha Q[u] \tag{2}$$

where $\alpha$ is an infinitesimal parameter. The function $Q[u] \equiv Q(x, t, u, u_x, u_t, ...)$ is called the *symmetry characteristic* of the transformation (2).

In order that a function $Q[u]$ be a symmetry characteristic for the PDE $F[u]=0$, it must satisfy a certain PDE that expresses the *symmetry condition* of $F[u]=0$. We write, symbolically,

$$S(Q\,;u)=0 \quad when \quad F[u]=0 \tag{3}$$

where the expression $S$ depends *linearly* on $Q$ and its partial derivatives. Thus, (3) is a linear PDE for $Q$, in which equation the variable $u$ enters as a sort of parametric function that is required to satisfy the PDE $F[u]=0$.

A *recursion operator* $\hat{R}$ [7] is a linear operator which, acting on a symmetry characteristic $Q$, produces a new symmetry characteristic $Q'=\hat{R}Q$. That is,

$$S(\hat{R}Q\,;u)=0 \quad when \quad S(Q\,;u)=0 \tag{4}$$

Obviously, *any power of a recursion operator also is a recursion operator*. This means that, starting with any symmetry characteristic $Q$, one may in principle obtain an infinite set of characteristics (thus, an infinite number of symmetries) by repeated application of the recursion operator.

A new approach to recursion operators was suggested in the early 1990s [9-12] (see also [13,14]). According to this view, a recursion operator is an auto-Bäcklund transformation for the linear PDE (3) that expresses the symmetry condition of the problem; that is, a BT producing new solutions $Q'$ of (3) from old ones, $Q$. Typically, this type of BT produces *nonlocal* symmetries, i.e., symmetry characteristics depending on *integrals* (rather than derivatives) of $u$.

As an example, consider the *chiral field equation* of Sec. 6.3, an alternative form of which is

$$F[g]\equiv (g^{-1}g_x)_x+(g^{-1}g_t)_t\ =0 \tag{5}$$

(as usual, subscripts denote partial differentiations) where $g$ is a $GL(n,C)$-valued function of $x$ and $t$ (i.e., an invertible complex $n\times n$ matrix, differentiable in $x$ and $t$).

Let $Q[g]$ be a symmetry characteristic of the PDE (5). It is convenient to put

$$Q[g]\ =g\,\Phi[g]$$

and write the corresponding infinitesimal symmetry transformation in the form

$$g'=g+\delta g\ =\ g+\alpha g\,\Phi[g] \tag{6}$$

The symmetry condition that $Q$ must satisfy will be a PDE linear in $Q$, thus in $\Phi$ also. As can be shown [8] this PDE is

$$S(\Phi;g)\equiv \Phi_{xx}+\Phi_{tt}+[g^{-1}g_x\,,\,\Phi_x]+[g^{-1}g_t\,,\,\Phi_t]=0 \tag{7}$$

which must be valid when $F[g]=0$ (where, in general, $[A,B]\equiv AB–BA$ denotes the commutator of two matrices $A$ and $B$).

For a given $g$ satisfying $F[g]=0$, consider now the following system of PDEs for the matrix functions $\Phi$ and $\Phi'$:

$$\begin{aligned} \Phi'_x &= \Phi_t + [g^{-1} g_t \,,\, \Phi] \\ -\Phi'_t &= \Phi_x + [g^{-1} g_x \,,\, \Phi] \end{aligned} \tag{8}$$

The integrability condition $(\Phi'_x)_t = (\Phi'_t)_x$, together with the equation $F[g]=0$, require that $\Phi$ be a solution of (7): $S(\Phi\,;g)=0$. Similarly, by the integrability condition $(\Phi_t)_x = (\Phi_x)_t$ one finds, after a lengthy calculation: $S(\Phi';g)=0$.

In conclusion, for any $g$ satisfying the PDE (5), the system (8) is a BT relating solutions $\Phi$ and $\Phi'$ of the symmetry condition (7) of this PDE; that is, relating different symmetries of the chiral field equation (5). Thus, if a symmetry characteristic $Q=g\Phi$ of (5) is known, a new characteristic $Q'=g\Phi'$ may be found by integrating the BT (8); the converse is also true. Since the BT (8) produces new symmetries from old ones, it may be regarded as a *recursion operator* for the PDE (5).

As an example, for any constant matrix $M$ the choice $\Phi=M$ clearly satisfies the symmetry condition (7). This corresponds to the symmetry characteristic $Q=gM$. By integrating the BT (8) for $\Phi'$, we get $\Phi'=[X, M]$ and $Q'=g[X, M]$, where $X$ is the "potential" of the PDE (5), defined by the system of PDEs

$$X_x = g^{-1} g_t \;, \quad -X_t = g^{-1} g_x \tag{9}$$

Note the *nonlocal* character of the BT-produced symmetry $Q'$, due to the presence of the potential $X$. Indeed, as seen from (9), in order to find $X$ one has to *integrate* the chiral field $g$ with respect to the independent variables $x$ and $t$. The above process can be continued indefinitely by repeated application of the recursion operator (8), leading to an infinite sequence of increasingly nonlocal symmetries.

## APPENDIX A

# CONSERVATIVE AND IRROTATIONAL FIELDS

Let $\vec{F}(\vec{r})$ be a static force field. According to the definition given in Sec. 1.5, this field is *conservative* if the work it does on a test particle of mass $m$ is path-independent, or equivalently, if

$$\oint_C \vec{F}(\vec{r}) \cdot \overrightarrow{dr} = 0 \tag{1}$$

for any closed path $C$ within the field. If $S$ is an open surface in the field, bounded by a given closed curve $C$ (see Fig. A.1), then, by Stokes' theorem and Eq. (1),

$$\oint_C \vec{F}(\vec{r}) \cdot \overrightarrow{dr} = \int_S (\vec{\nabla} \times \vec{F}) \cdot \overrightarrow{da} = 0 \tag{2}$$

In order for this to be true for every open surface $S$ bounded by $C$, the field $\vec{F}(\vec{r})$ must be *irrotational*:

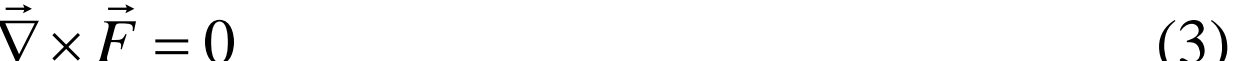

$$\vec{\nabla} \times \vec{F} = 0 \tag{3}$$

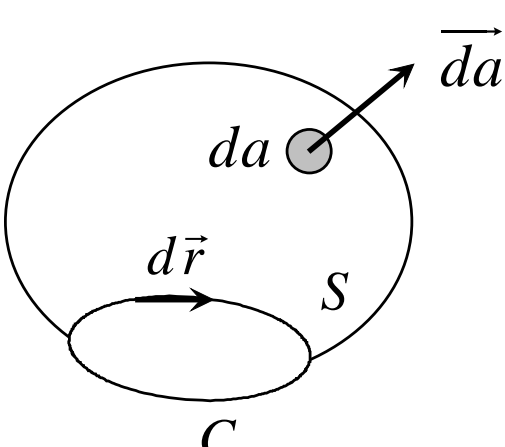

Fig. A.1. An open surface $S$ bounded by a closed curve $C$.

Conversely, an irrotational force field $\vec{F}(\vec{r})$ will also be conservative in a region of space that is *simply connected*. Indeed, given any closed curve $C$ in such a region, it is always possible to find an open surface $S$ having $C$ as its boundary. Then, if (3) is valid, the force is conservative in view of (2).

Given a conservative force field $\vec{F}(\vec{r})$, there exists a function $U(\vec{r})$ (*potential energy* of the test particle $m$) such that

$$\vec{F} = -\vec{\nabla} U \tag{4}$$

The work done on $m$ by $\vec{F}$, when this particle moves along *any* path from point $A$ to point $B$ in the field, is then equal to

$$W = \int_A^B \vec{F}(\vec{r}) \cdot \overrightarrow{dr} = U(\vec{r}_A) - U(\vec{r}_B) \tag{5}$$

Consider now a *time-dependent* force field $\vec{F}(\vec{r},t)$ in a *simply connected* region $\Omega$ of space. This field is assumed to be *irrotational* for all values of $t$:

$$\vec{\nabla}\times\vec{F}(\vec{r},t)=0 \tag{6}$$

Can we conclude that the field $\vec{F}$ is conservative?

It is tempting but *incorrect* (!) to argue as follows: Let $C$ be an arbitrary closed curve in $\Omega$. Since $\Omega$ is simply connected, there is always an open surface $S$ bounded by $C$. By Stokes' theorem,

$$\oint_C \vec{F}(\vec{r},t)\cdot\overrightarrow{dr}=\int_S(\vec{\nabla}\times\vec{F})\cdot\overrightarrow{da}=0 \tag{7}$$

for all values of $t$. This *appears* to imply that $\vec{F}$ is conservative. This is not so, however, for the following reason: For any fixed value of $t$, the integral

$$I(t)=\oint_C \vec{F}(\vec{r},t)\cdot\overrightarrow{dr}$$

does *not* represent work. Indeed, $I(t)$ expresses the integration of a function of two independent variables $\vec{r}$ and $t$ over one of these variables (namely, $\vec{r}$), the other variable ($t$) playing the role of a "parameter" of integration which remains *fixed*. Thus, $I(t)$ is evaluated *for a given instant of time t*, and all values of $\vec{F}$ at the various points of $C$ must be recorded *simultaneously* at this time.

On the other hand, in the integral representation of work, time is assumed to flow as the particle $m$ travels along the closed path $C$. In this case, $\vec{r}$ and $t$ are no longer considered independent of each other but are connected through the equation of motion of $m$. Furthermore, the values of the force $\vec{F}$ at different points of $C$ must now be recorded at *different* times, corresponding to the times the particle passes from the respective points. Equation (7) is therefore not valid if the integral on the left-hand side is interpreted to be the work of $\vec{F}$ on $m$ along $C$. This complication never arises in the case of a static force field, as in that case the time at which the value of $\vec{F}$ is recorded at any given point of the path is immaterial.

We thus conclude the following:

*1. A force field that is both static and irrotational in a simply connected region of space is conservative.*

*2. A time-dependent force field cannot be conservative even if it is irrotational and its region of action is simply connected.*

Finally, let us explain why a time-dependent force field does not lead to conservation of total mechanical energy. Consider again an irrotational force field $\vec{F}(\vec{r},t)$ [as defined according to (6)] in a simply connected region $\Omega$. Then there exists a time-dependent potential energy $U(\vec{r},t)$ of $m$, such that, for any value of $t$,

$$\vec{F}(\vec{r},t) = -\vec{\nabla} U(\vec{r},t) \tag{8}$$

This time we will assume that $\vec{F}(\vec{r},t)$ is the *total* force on *m*. By Newton's 2nd law, then,

$$m\frac{d\vec{v}}{dt} = \vec{F} \quad (\text{where } \vec{v} = \overrightarrow{dr}/dt) \;\Rightarrow\; m\frac{d\vec{v}}{dt} + \vec{\nabla} U = 0 \ .$$

Taking the dot product with $\vec{v}$, this vector being used here as an *integrating factor* (see Sec. 3.3), we have:

$$m\vec{v}\cdot\frac{d\vec{v}}{dt} + \vec{v}\cdot\vec{\nabla} U = 0 \ .$$

Now,

$$\vec{v}\cdot\frac{d\vec{v}}{dt} = \frac{1}{2}\frac{d}{dt}(\vec{v}\cdot\vec{v}) = \frac{1}{2}\frac{d}{dt}(v^2) \quad (v = |\vec{v}\,|)$$

and

$$\vec{v}\cdot\vec{\nabla} U \;=\; \frac{\vec{\nabla} U \cdot \overrightarrow{dr}}{dt} \;=\; \frac{dU - \dfrac{\partial U}{\partial t}dt}{dt} \;=\; \frac{dU}{dt} - \frac{\partial U}{\partial t}$$

where we have used the fact that

$$dU(\vec{r},t) = \vec{\nabla} U \cdot \overrightarrow{dr} + \frac{\partial U}{\partial t}dt \ .$$

Hence, finally,

$$\frac{d}{dt}\left(\frac{1}{2}mv^2\right) + \frac{dU}{dt} - \frac{\partial U}{\partial t} = 0 \;\Rightarrow$$

$$\frac{d}{dt}(E_k + U) = \frac{\partial U}{\partial t} \tag{9}$$

where $E_k = mv^2/2$ is the kinetic energy of *m*. As seen from (9), the total mechanical energy $(E_k+U)$ is not conserved *unless* $\partial U/\partial t=0$, i.e., unless the force field is *static*.

Note that, for a time-dependent irrotational force field [defined according to (6)] the quantity

$$\int_A^B \vec{F}(\vec{r},t)\cdot\overrightarrow{dr} \;=\; U(\vec{r}_A,t) - U(\vec{r}_B,t) \ ,$$

defined for any *fixed t*, does no longer represent the work done by this field on a particle moving from *A* to *B* [comp. Eq. (5) for the case of a static force field].

Let us summarize our main conclusions:

1. A static force field that is irrotational in a simply connected region of space is conservative.

2. A time-dependent force field cannot be conservative even if it is irrotational and its region of action has the proper topology.

3. The work of a time-dependent irrotational force field cannot be expressed as the difference of the values of the time-dependent potential energy at the end points of the trajectory.

4. Time-dependent force fields are incompatible with conservation of total mechanical energy.

**APPENDIX B**

# MATRIX DIFFERENTIAL RELATIONS

Let $A(t)=[a_{ij}(t)]$ be an $(m\times n)$ matrix whose elements are functions of $t$. The *derivative* $dA/dt$ of $A$ is the $(m\times n)$ matrix with elements

$$\left(\frac{dA}{dt}\right)_{ij}=\frac{d}{dt}\,a_{ij}(t) \qquad (1)$$

If $B(t)$ is another $(m\times n)$ matrix, then

$$\frac{d}{dt}\,(A\pm B)=\frac{dA}{dt}\pm\frac{dB}{dt} \qquad (2)$$

For square $(n\times n)$ matrices $A$, $B$, $C$,

$$\frac{d}{dt}\,(AB)=\frac{dA}{dt}B+A\frac{dB}{dt}\ ,$$
$$\frac{d}{dt}\,(ABC)=\frac{dA}{dt}BC+A\frac{dB}{dt}C+AB\frac{dC}{dt}\ , \qquad (3)$$

etc. Similarly, the *integral* of an $(m\times n)$ matrix function $A(t)=[a_{ij}(t)]$ is defined by

$$\left(\int A(t)\,dt\right)_{ij}=\int a_{ij}(t)\;dt \qquad (4)$$

The derivative of the *inverse* $A^{-1}$ of a non-singular $(n\times n)$ matrix $A$ is given by

$$\frac{d}{dt}\,(A^{-1})=-A^{-1}\,\frac{dA}{dt}\,A^{-1} \qquad (5)$$

Indeed, given that $A^{-1}A=1_n$ [where $1_n$ is the unit $(n\times n)$ matrix], we have:

$$\frac{d}{dt}\,(A^{-1}A)=0\;\Rightarrow\;\frac{d(A^{-1})}{dt}A+A^{-1}\frac{dA}{dt}=0\;\Rightarrow\;\frac{d(A^{-1})}{dt}A=-A^{-1}\frac{dA}{dt}\ .$$

Multiplying from the right by $A^{-1}$, we get (5).

As is easy to show with the aid of (2) and (3), for square matrices $A$ and $B$ we have:

$$\frac{d}{dt}\,[A,B]=\left[\frac{dA}{dt}\,,\,B\right]+\left[A\,,\,\frac{dB}{dt}\right] \qquad (6)$$

where by $[A\,,B]\equiv AB-BA$ we denote the *commutator* of two matrices.

Assume now that $A(x,y)= [a_{ij}(x,y)]$ is an invertible square-matrix function. We call $A_x$ and $A_y$ the partial derivatives of $A$ with respect to $x$ and $y$, respectively. The following identities are valid:

$$\begin{aligned} &\partial_x (A^{-1}A_y) - \partial_y (A^{-1}A_x) + [A^{-1}A_x \,,\, A^{-1}A_y] = 0 \\ &\partial_x (A_y A^{-1}) - \partial_y (A_x A^{-1}) - [A_x A^{-1},\, A_y A^{-1}] = 0 \end{aligned} \qquad (7)$$

Moreover,

$$A(A^{-1}A_x)_y A^{-1} = (A_y A^{-1})_x \;\Leftrightarrow\; A^{-1}(A_y A^{-1})_x A = (A^{-1}A_x)_y \qquad (8)$$

Given a *constant* $(n\times n)$ matrix $A$ (where by "constant" we mean independent of $t$) we define the *exponential matrix* $e^{tA}$ by

$$e^{tA} \equiv \exp(tA) = \sum_{k=0}^{\infty} \frac{1}{k!}(tA)^k = 1_n + tA + \frac{t^2}{2!}A^2 + \frac{t^3}{3!}A^3 + \frac{t^4}{4!}A^4 + \cdots \qquad (9)$$

The $(n\times n)$ matrix $e^{tA}$ is a function of $t$. Its derivative is given by

$$\begin{aligned} \frac{d}{dt}e^{tA} &= 0 + A + tA^2 + \frac{t^2}{2!}A^3 + \frac{t^3}{3!}A^4 + \cdots \\ &= A\left(1_n + tA + \frac{t^2}{2!}A^2 + \frac{t^3}{3!}A^3 + \cdots\right) \\ &= \left(1_n + tA + \frac{t^2}{2!}A^2 + \frac{t^3}{3!}A^3 + \cdots\right)A \;\Rightarrow \end{aligned}$$

$$\frac{d}{dt}e^{tA} = A\,e^{tA} = e^{tA}A \qquad (10)$$

Putting $-A$ in place of $A$, we have:

$$\frac{d}{dt}e^{-tA} = -A\,e^{-tA} = -\,e^{-tA}A \qquad (11)$$

***Exercises:***

***1.*** Using property (3) for the derivative of a product of matrices, and taking into account Eqs. (10) and (11), show that, for *constant* $(n\times n)$ matrices $A$ and $B$,

$$\frac{d}{dt}\left(e^{-tA}B\,e^{tA}\right) = \left[e^{-tA}B\,e^{tA},\, A\right] \qquad (12)$$

where the square bracket on the right denotes the commutator

$$(e^{-tA}Be^{tA})A - A(e^{-tA}Be^{tA})\,.$$

***2.*** Consider the ODE and the associated initial condition,

$$\frac{d}{dt}u(t) = Au(t) \quad , \quad u(0) = u_0$$

where $A$ and $u_0$ are constant $(n\times n)$ matrices while $u(t)$ is a matrix function. Show that the solution of this ODE is

$$u(t) = e^{tA}u_0 \ .$$

Similarly, show that the solution of the ODE

$$\frac{d}{dt}u(t) = u(t)A \quad , \quad u(0) = u_0$$

is given by

$$u(t) = u_0\, e^{tA} .$$

***3.*** Show that the solution of the matrix ODE

$$\frac{d}{dt}u(t) = [u(t),\, A] \quad , \quad u(0) = u_0$$

(for constant $A$ ) is

$$u(t) = e^{-tA}u_0\, e^{tA} \ .$$

# GENERAL BIBLIOGRAPHY